\documentclass[opre,sglanonrev]{informs}
\RequirePackage{tgtermes}
\RequirePackage{newtxtext}
\RequirePackage{newtxmath}
\RequirePackage{bm}
\RequirePackage{endnotes}

\usepackage[linesnumbered,ruled,lined,commentsnumbered]{algorithm2e}
\usepackage{tikz}
\usepackage[utf8]{inputenc}   \usepackage[T1]{fontenc}      \usepackage{amsfonts,amscd,amssymb,amsmath}
\usepackage{array,mathtools,bm}
\usepackage{booktabs,caption,longtable,multicol,multirow}
\usepackage{subcaption,lscape,enumitem}
\usepackage[colorlinks,citecolor=purple,linkcolor=blue]{hyperref}
\usepackage{arydshln,xfrac}
\usepackage{csquotes}
\allowdisplaybreaks
\usepackage{tikz}
\usepackage{adjustbox}
\usetikzlibrary{shadings}
\usetikzlibrary{arrows.meta}
\definecolor{ao(english)}{rgb}{0.0, 0.5, 0.0}
\definecolor{cadmiumgreen}{rgb}{0.0, 0.42, 0.24}
\definecolor{darkpastelgreen}{rgb}{0.01, 0.75, 0.24}
\definecolor{rblue}{RGB}{156,166,192}

\usepackage{tcolorbox}
\usepackage{pgfplots}
\usepackage{cleveref}
\tcbuselibrary{listings,theorems,breakable,skins}

\newtheorem{theorem}{Theorem}[section]
\newtheorem{lemma}{Lemma}[section]

\newtheorem{assumption}{Assumption}
\newtheorem{corollary}{Corollary}[section]
\newtheorem{remark}{Remark}
\newtheorem{example}{Example}

\crefname{observation}{Observation}{Observations}
\crefname{theorem}{Theorem}{Theorems}
\crefname{lemma}{Lemma}{Lemmas}
\crefname{definition}{Definition}{Definitions}
\crefname{assumption}{Assumption}{Assumptions}
\crefname{corollary}{Corollary}{Corollaries}
\crefname{remark}{Remark}{Remarks}
\crefname{example}{Example}{Examples}
\crefname{section}{Section}{Sections}
\crefname{subsection}{Section}{Sections}
\crefname{subsubsection}{Section}{Sections} \newcommand{\lump}{\text{LUMP}}
\DeclareMathSymbol{\shortminus}{\mathbin}{AMSa}{"39}
\newcommand{\st}{\mathrm{s.t.}~}
\newcommand{\diag}{\mathrm{diag}}

\newcommand{\real}{\mathbb R}

\newcommand{\cov}{\mathrm{cov}}
\newcommand{\va}{\mathrm{var}}
\newcommand{\inv}{\mathrm{\shortminus1}}
\newcommand{\invhalf}{\mathrm{\shortminus0.5}}
\newcommand{\qlq}{~\Rightarrow~}

\newcommand{\diff}{\mathrm{d}}

\newcommand{\1}{\mathbf{1}}

\newcommand{\inner}[2]{\langle #1,#2\rangle}

\definecolor{rosemary}{RGB}{190,88,69}
\definecolor{stone}{RGB}{42,121,10}
\definecolor{deepblack}{RGB}{90,90,90}

\hypersetup{
    colorlinks,
    linkcolor={deepblack},
    citecolor={rosemary},
    urlcolor={rosemary}
}

\newenvironment{eqalign}[1]
{\setlength{\abovedisplayskip}{6pt}\setlength{\belowdisplayskip}{3pt}\equation\label{#1}\aligned[b]}{\endaligned\endequation}
\newenvironment{subarr}[2]
{\setlength{\abovedisplayskip}{6pt}\setlength{\belowdisplayskip}{3pt}
    \equation\label{#1}\array{#2}}
{\endarray\endequation}
\newenvironment{subalign}[1]
{\setlength{\abovedisplayskip}{6pt}\setlength{\belowdisplayskip}{3pt}\subequations\label{#1}\align}
{\endalign\endsubequations}

\let\oldequation\equation
\let\endoldequation\endequation
\renewenvironment{equation}
{\setlength{\abovedisplayskip}{6pt}\setlength{\belowdisplayskip}{3pt}\oldequation}
{\endoldequation}

\everydisplay{\setlength{\abovedisplayskip}{6pt}\setlength{\belowdisplayskip}{3pt}}

\makeatletter
\newcommand\NoIndentNext{\@afterindentfalse\@afterheading}
\makeatother
\newcommand{\NoIndentAfter}[1]{\AfterEndEnvironment{#1}{\NoIndentNext}}
\forcsvlist{\NoIndentAfter}{theorem,lemma,proposition,corollary,definition,remark,proof}

\def\b{\mathbf{b}}

\def\c{\mathbf{c}}

\def\d{\mathbf{d}}

\def\f{\mathbf{f}}

\def\g{\mathbf{g}}

\def\h{\mathbf{h}}

\def\p{\mathbf{p}}

\def\q{\mathbf{q}}

\def\r{\mathbf{r}}

\def\s{\mathbf{s}}

\def\v{\mathbf{v}}

\def\w{\mathbf{w}}

\def\x{\mathbf{x}}

\def\y{\mathbf{y}}

\def\z{\mathbf{z}}

\def\cA{\mathcal{A}}

\def\bA{\mathbf{A}}

\def\cB{\mathcal{B}}

\def\bB{\mathbf{B}}

\def\cC{\mathcal{C}}

\def\cD{\mathcal{D}}

\def\bD{\mathbf{D}}

\def\vbE{\vmathbb{E}}

\def\cF{\mathcal{F}}

\def\bF{\mathbf{F}}

\def\bG{\mathbf{G}}

\def\cH{\mathcal{H}}

\def\bH{\mathbf{H}}

\def\cI{\mathcal{I}}

\def\bI{\mathbf{I}}

\def\cJ{\mathcal{J}}

\def\cK{\mathcal{K}}

\def\bK{\mathbf{K}}

\def\cL{\mathcal{L}}

\def\bL{\mathbf{L}}

\def\bM{\mathbf{M}}

\def\bbN{\mathbb{N}}

\def\cO{\mathcal{O}}

\def\cP{\mathcal{P}}

\def\vbP{\vmathbb{P}}

\def\bP{\mathbf{P}}

\def\cQ{\mathcal{Q}}

\def\bQ{\mathbf{Q}}

\def\cR{\mathcal{R}}
\def\bbR{\mathbb{R}}

\def\bS{\mathbf{S}}

\def\bT{\mathbf{T}}

\def\cU{\mathcal{U}}

\def\cV{\mathcal{V}}

\def\bW{\mathbf{W}}

\def\cX{\mathcal{X}}

\def\bX{\mathbf{X}}

\def\bGamma{\boldsymbol{\Gamma}}
\def\bgamma{\boldsymbol{\gamma}}

\def\btheta{\boldsymbol{\theta}}

\def\bxi{\boldsymbol{\xi}}

\makeatother

\newcommand{\kc}{\mathbf{k}_{\mathrm{c}}}
\newcommand{\Kc}{\mathbf{K}_{\mathrm{c}}}

\newcommand{\logbar}{\texttt{LogBar}}
\newcommand{\logbarpcg}{\texttt{LogBar-PCG}}

\newcommand{\pathf}{\texttt{PathFol}}
\newcommand{\tat}{\texttt{Tât}}
\newcommand{\pr}{\texttt{PropRes}}

\crefname{algorithm}{Second-order Tâtonnement}{Second-order Tâtonnements}

\usepackage{natbib}
 \bibpunct[, ]{(}{)}{,}{a}{}{,}\def\bibfont{\small}\def\BIBand{and}

\OneAndAHalfSpacedXI
\EquationsNumberedBySection 

\ECRepeatTheorems  

\MANUSCRIPTNO{MOOR-0001-2024.00}

\begin{document}

\RUNAUTHOR{Zhang et al.}

\RUNTITLE{Price Updates by Interior-Point Approaches}

\TITLE{The Second-Order Tâtonnement: Decentralized Interior-Point Methods for Market Equilibrium}

\ARTICLEAUTHORS{\AUTHOR{Chuwen Zhang\textsuperscript{a}, Chang He\textsuperscript{b}, Bo Jiang\textsuperscript{b,c,d,*}, Yinyu Ye\textsuperscript{e,f}
    }
    \AFF{
        \textsuperscript{a} Booth School of Business, University of Chicago,
        \EMAIL{chuwen.zhang@chicagobooth.edu} \\
        \textsuperscript{b} Research Institute for Interdisciplinary Sciences, School of Information Management and Engineering, Shanghai University of Finance and Economics,
        \EMAIL{ischanghe@gmail.com, isyebojiang@gmail.com} \\
        \textsuperscript{c} Key Laboratory of Interdisciplinary Research of Computation and Economics, Shanghai University of Finance and Economics \\
        \textsuperscript{d} Dishui Lake Advanced Institute, Shanghai University of Finance and Economics \\
        \textsuperscript{e} Antai College of Economics and Management, Shanghai Jiao Tong University, \EMAIL{yinyuye@sjtu.edu.cn}\\
        \textsuperscript{f} Department of Management Science and Engineering, Stanford University (Emeritus), \EMAIL{yyye@stanford.edu}\\
        \textsuperscript{*} Corresponding authors
    }
} 

\ABSTRACT{The t\^atonnement process and Smale's process are two classical approaches to compute market equilibrium in exchange economies. While the t\^atonnement process can be seen as a first-order method, Smale's process, being second-order, is less popular due to its reliance on additional information from the players and expensive Newton steps. In this paper, we study Fisher exchange market for a broad class of utility functions, where we show that all high-order information required by Smale's process is readily available from players' best responses.
    Motivated by this observation, we develop two second-order t\^atonnement processes, constructed as decentralized interior-point methods, which are traditionally known to work in a centralized manner. The methods here bear the name ``t\^atonnement'', since, in spirit, they demand no more information than the classical t\^atonnement process.
    To address the Newton systems involved, we introduce an explicitly invertible approximation with high-probability guarantees and a scaling matrix that optimally minimizes the condition number, both of which rely solely on best responses as the methods themselves.
    Using these tools, the first second-order t\^atonnement process has $\cO(\log(\tfrac{1}{\epsilon}))$ complexity rate.
    Under mild conditions, the other method achieves a non-asymptotic superlinear convergence rate.
    Preliminary experiments are presented to justify the capability of the proposed methods for large-scale problems. Extensions of our approach are also discussed.
}


\maketitle

{\small\setlength{\baselineskip}{0.9\baselineskip}
    \tableofcontents
}

\section{Introduction}

In the Arrow-Debreu competitive exchange market problem (e.g., \cite{walrasElementsDeconomiePolitique1874}), a set of divisible goods $j\in \cJ=\{1, ..., n\}$ is allocated to a population of players $i\in \cI=\{1, ..., m\}$. Without loss of generality, we assume that there is a unit amount of goods for each $j \in \cJ$. In the marketplace, the \emph{prices} for the goods $p^{(1)},...,p^{(n)}$ belong to a subset of the nonnegative orthant $\cP\subseteq \real_+^n$.
Each of the $m$ players sells their initial endowment and buys a bundle of goods, $\x_i$, to maximize his or her concave utility $u_i$. The bundle, also referred to as the \emph{allocation vector}, is restricted to a common set $\cX \subseteq \real_+^n$ that is a closed convex cone.
Besides, each player $i$ cannot spend more than the income.

Fisher considers a simpler model \citep{fisherMathematicalInvestigationsTheory1892} where the initial endowment is given as an exogenous budget $w_i \in \real_+$.
An equilibrium of the market is $(\x_1, ..., \x_m, \p)$ at which every player $i$ is optimally satisfied under the competitive budget $w_i$, and the market has no deficiency or surplus such that:
\begin{eqalign}{eq.admarket.equilibrium}
    \p \in \cP, \quad \x_i \in ~& \underset{\x_i \in \cX}{\arg \max} ~ u_i(\x_i), ~\st, \inner{\p}{\x_i} \le w_i; \\
    ~\textrm{and} ~ &\sum_{i\in \cI}  \x_i = \1.
\end{eqalign}
Here we assume that every good is valued by at least one consumer, and so the price vector $\p$ is strictly positive.
At such an equilibrium, we say the market \emph{clears}: the demand matches the supply, i.e., the above equality holds.

The study of market equilibrium has been a central topic in capitalist economic systems \citep{walrasElementsDeconomiePolitique1874,arrowExistenceEquilibriumCompetitive1954,saariEffectivePriceMechanisms1978,smaleConvergentProcessPrice1976}.
The scale of modern internet platforms has reignited such interest due to several appealing features.
For example, the Fisher market model is known for distributive fairness: with equal budgets, it yields the Competitive Equilibrium from Equal Incomes (CEEI) \citep{varianEquityEnvyEfficiency1974}, which is envy-free and Pareto efficient. From a standpoint of mechanism design, the approximate market equilibrium is strategy-proof in large markets \citep{robertsIncentivesPricetakingBehavior1976,budishCombinatorialAssignmentProblem2011,azevedoStrategyproofnessLarge2019}: the players would act as price takers with little incentive to misreport their true preferences.

Except for fairness concerns in public goods, healthcare, or humanitarian applications, a market equilibrium is sometimes desirable even for a profit-seeking marketplace \citep{borgsDynamicsBidOptimization2007,bateniFairResourceAllocation2022a}. For example, in online advertising, the equilibrium approach yields stable and convergent mechanisms \citep{borgsDynamicsBidOptimization2007} and significantly alleviates the congestion issues (e.g., career platforms) typical among click-through rate (CTR) based mechanisms \citep{kroerComputingLargeMarket2022}.

\subsection{Related work and motivation}
Generally, the existence of a market equilibrium is warranted by Brouwer's fixed-point theorem. This viewpoint has motivated Scarf's homotopy method \citep{scarfApproximationFixedPoints1967}. For linear utilities, it may also be formulated as a linear complementarity problem (LCP) \citep{facchineiFiniteDimensionalVariationalInequalities2004}, which admits finite convergence via Lemke's pivotal algorithm. However, the worst-case running time of these algorithms are not polynomial \citep{eavesFiniteAlgorithmLinear1976,yePathArrowDebreu2008}.
Most of the contemporary approaches are based on the Eisenberg-Gale (EG) convex program \citep{eisenbergConsensusSubjectiveProbabilities1959} whenever the players have homogeneous-of-order-one utility functions:
\begin{subalign}{eq.primal.eg}
    \underset{\x_{i} \in \cX, i \in \cI}{\max ~} & \sum_{i\in \cI} w_i \log(u_i(\x_i)) \\
    \label{eq.primal.eg.clear} \st & \sum_{i\in \cI} \x_i \le \1.
\end{subalign}
It is known that the dual solution to \eqref{eq.primal.eg.clear} corresponds to the equilibrium prices $\p$ that clear the market. This convex optimization problem \eqref{eq.primal.eg} allows straightforward applications of the ellipsoid method \citep{khachiyanPolynomialAlgorithmLinear1979,grotschelGeometricAlgorithmsCombinatorial1993}. One can also reformulate it as a linear conic program \citep{ben-talLecturesModernConvex2001,skajaaHomogeneousInteriorpointAlgorithm2015}, allowing for the application of polynomial-time interior-point methods. For the market problem with linear utilities, strongly polynomial-time algorithms have been proposed in the early 21st century \citep{devanurMarketEquilibriumPrimaldualtype2002,gargAuctionAlgorithmsMarket2006,yePathArrowDebreu2008}. However, these approaches decide $\{\x_i\}_{i\in \cI}$ in a \emph{centralized} manner (it assumes the power to allocate the goods for all players). On one hand, it is sometimes not realistic to do so, whereas only the outcomes of market prices (e.g., revealed preferences of the players) are available. On the other hand, centralized approaches do not scale very well in the number of players and goods, which may be inefficient for large-scale problems. For example, the worst-case complexity of the interior-point method \citep{yePathArrowDebreu2008} is $\cO(\sqrt{mn}(m+n)^3\log(\frac{1}{\epsilon}))$ (recall that the number of players and goods are $m$ and $n$, respectively).

In fact, centralization is unnecessary for the purpose of calculating and implementing the equilibrium prices. An alternative class of approaches is the decentralized price updating mechanisms. For example, a prototypical scheme could work as follows:
\begin{itemize}[leftmargin=*]
    \item Collect each player's \emph{best response} $\x_i(\p)$ from his or her Utility Maximization Problem (UMP),
          \begin{subalign}{eq.ump.def}
              \label{eq.ump} (\text{UMP}) \quad \pi_i(\p) = \max_{\x_i \in \cX} u_i(\x_i), ~\st~\inner{\p}{\x_i} \le w_i \\
              \label{eq.ump.x} \text{and}\quad \x_i(\p) \text{ is the solution to}~\text{UMP}.
          \end{subalign}
          The best response is also referred to as the Walrasian demand function. The optimal utility value $\pi_i(\p)$ is denoted as the \emph{indirect utility} \citep{mas-colellMicroeconomicTheory1995}.
    \item Compute the \emph{excess demand function},
          $\z(\p) = \sum_{i\in \cI} \x_i(\p)- \1;$
    \item Update $\p$ based on some rule $G$ that maps $\z(\p)$ to the price space:
          \begin{equation}\label{eq.tatonnement}
              \p_+ \leftarrow \p + G(\z(\p)).
          \end{equation}
\end{itemize}
Strictly speaking, the above scheme is a discrete t\^atonnement process \citep{walrasElementsDeconomiePolitique1874,mas-colellMicroeconomicTheory1995,codenottiComputationMarketEquilibria2007}.
Note that at each round, only the players' responses are needed, so that it is indeed a decentralized approach.
The underlying economic intuition is that when aggregate demand exceeds supply for a good, the market mechanism responds by adjusting prices upward.
To discuss its convergence, we turn to the dual EG problem \citep{coleConvexProgramDuality2017,goktasConsumerTheoreticCharacterizationFisher2022} of \eqref{eq.primal.eg}:
\begin{equation}\label{eq.indirect.dual.eg}
    \min_{\p\in\cP} ~ \varphi(\p) := \inner{\p}{\1} + \sum_{i\in\cI} w_i \log(\pi_i(\p)).
\end{equation}
In this paper, we call $\varphi$ the \emph{potential function}.
The connection between the t\^atonnement \eqref{eq.tatonnement} and dual EG program \eqref{eq.indirect.dual.eg} is established via the following fact:
\begin{equation}\label{eq.dual.grad}
    w_i \nabla \log(\pi_i(\p)) = - \x_i(\p),~ \nabla \varphi(\p) = \1 - \sum_{i\in\cI} \x_i(\p) = -\z(\p).
\end{equation}
Hence, finding the market equilibrium boils down to a task of finding the first-order stationary point of $\varphi$.
If all $\x_i(\p)$ are unique, then $\varphi$ is differentiable; otherwise, the best responses produce the subdifferential.
In either case, it is quite natural to use first-order methods that work with a finite-sum optimization problem.

In this vein, two specific types of first-order methods (price updating mechanisms) are well-studied for the Fisher model.  The first one is the aforementioned t\^atonnement process dating back to \cite{walrasElementsDeconomiePolitique1874}.  The convergence was established for economies satisfying weak gross substitutability (WGS) \citep{arrowStabilityCompetitiveEquilibrium1958,arrowStabilityCompetitiveEquilibrium1959}.
Another one is Proportional Response (PR), originated from \cite{shapleyTradeUsingOne1977} and recently applied in the peer-to-peer networks \citep{wuProportionalResponseDynamics2007,zhangProportionalResponseDynamics2011}.
Strictly speaking, PR is driven by the players rather than the prices: the players submit their bidding results (the proportion of money spent on each good in the bundle), and then the prices are aggregated based on the players' bids. Hence, it has a direct connection to Shmyrev's convex program \citep{shmyrevAlgorithmFindingEquilibrium2009}.
While the two dynamics are motivated differently, both the t\^atonnement and PR dynamics can be understood as non-Euclidean gradient descent methods through the lens of relative smoothness \citep{bauschkeDescentLemmaLipschitz2017,luRelativelySmoothConvex2018}.
The key observation is that the distance to the equilibrium can be measured by Kullback–Leibler divergence (or its variants). For constant-elasticity-of-substitution (CES) utilities, both the t\^atonnement and PR dynamics have been shown to have a convergence rate of $\cO(\log(\tfrac{1}{\epsilon}))$. While in linear and Leontief economies, the convergence rate is $\cO(\tfrac{1}{\epsilon})$ \citep{birnbaumDistributedAlgorithmsGradient2011,cheungTatonnementGrossSubstitutes2013,cheungDynamicsDistributedUpdating2018,gaoFirstOrderMethodsLargeScale2020}\endnote{A linear convergent FOM is given \cite{nanFastInterpretableDynamics2023}, but it requires $\epsilon$ to be not too small.}.
First-order schemes can also tackle a related problem that computes the equilibrium for chores; see, for example, \cite{chenComputingCompetitiveEquilibrium2024,chaudhuryCompetitiveEquilibriumChores2024}.
The per-iteration costs of these methods can be significantly cheaper than centralized approaches. For example, each iteration of the t\^atonnement process requires $\cO(n)$ arithmetic operations to update the prices. Moreover, the dependency on the number of players is $\cO(m)$ due to decentralization.

Apart from the first-order schemes, to accelerate the convergence rate, we may wonder if it is possible to design a second-order price updating mechanism that leverages $\nabla\z(\p)$, the Jacobian of the excess demand function. One such scheme, named \enquote{Global Newton method}, was proposed in a seminal work \citep{smaleConvergentProcessPrice1976}, and is now known as Smale's process (see also \cite{heringsEquilibriumAdjustmentDisequilibrium1997,kamiyaGloballyStablePrice1990} for its variants).
In history, the motivation at that time was to resolve the non-convergent example in \cite{scarfExamplesGlobalInstability1960} in the Arrow-Debreu setting.
However, it turns out that Smale's process
has a few limitations. In particular, it
only has local convergence guarantees \citep{keenanFurtherRemarksGlobal1981} although bearing the name ``Global Newton''. Furthermore, no worst-case complexity has been established for either the Arrow-Debreu or the Fisher model.
Compared to the t\^atonnement process, Smale's process and the variants have no natural economic interpretation and are more informationally demanding \citep{jordanCompetitiveAllocationProcess1982}: except for the best responses, one also needs their marginal changes of the excess demand function. Acquiring such information is highly nontrivial even if the full preferences (e.g., the utilities) of the market are provided.
It was later shown that this informational requirement cannot be substantially reduced for an effective mechanism if one wants guaranteed convergence to equilibrium \citep{saariEffectivePriceMechanisms1978,saariIterativePriceMechanisms1985}.  In practice, Smale's process is also unrealistic: each update step requires computing the Newton steps, which seems to be prohibitively costly for large-scale problems in today's economic systems.

In this paper, we resolve the above drawbacks of Smale's process under a ubiquitous class of utility functions for the Fisher model.
Based on a few technical results well-known in the interior-point methods \citep{nesterovSelfscaledBarriersInteriorpoint1997}, we show that the second-order information needed for a price updating mechanism is free: the best responses $\x_1(\p),...,\x_m(\p)$ are all we need. Conceptually, in terms of information dependency, a second-order price updating mechanism has no difference from a t\^atonnement process.
Practically, we show that the computational burden of deriving the necessary Newton steps can also be significantly reduced by exploring lightweight strategies based on exactly the same information.
These nice properties enable us to design two second-order t\^atonnement processes, specifically, decentralized interior-point methods (IPMs) with worst-case complexity guarantees.

\subsection{Contributions}

This paper lays the mathematical apparatus for the decentralized interior-point methods to compute the Fisher market equilibrium. To facilitate the analysis, we introduce the Logarithmic Utility Maximization problem ($\lump{}$) whose best-response mapping $\x_i(\p)$ is equivalent to that of a classical UMP \eqref{eq.ump}.   Specifically, for a ubiquitous class that is commonly used in the literature, we show that the Jacobian $\nabla\x_i(\p)$ can be constructed from $\x_i(\p)$ itself. Hence, our IPMs consume no extra information from the players and the dependency on $m$ ties to the discrete t\^atonnement process. Thereby, the Scaled Lipschitz Continuity (SLC) of the best response are discussed to facilitate the worst-case complexity analysis.

To tackle the Newton systems arising from step computations in the IPMs, we introduce two lightweight approaches, again based on the best responses, to compute inexact Newton steps. Firstly, we introduce a simple Diagonal plus Rank-1 (DR1) approximation of the Jacobian operator with \emph{explicit} inverse, and hence the step computations require $\cO(n)$ arithmetic operations. The quality of this approximate operator can be guaranteed, with high probability, in large-market regimes with reasonably many players and goods (e.g., $m \ge 10^2, n\ge 10^2$). Alternatively, we propose a diagonal preconditioner to facilitate the Krylov solvers \citep{golubMatrixComputations2013}, which outputs inexact solutions of the Newton systems.
This preconditioner is provably optimal among diagonal preconditioning matrices. These schemes, taking advantage of the curvature information, may also be applicable to first-order price updating mechanisms.

In addition, we develop two interior-point methods to capitalize on these approaches for computing the inexact steps.
These methods can be seen as \enquote{second-order t\^atonnement} processes since they have the same information dependency as the classical t\^atonnement process.
The first method is based on the logarithmic barrier function.  We show that this price updating mechanism converges to the equilibrium prices in at most $\cO(\log(\tfrac{1}{\epsilon}))$ iterations. Under the conditions that warrant self-concordance \citep{nesterovLecturesConvexOptimization2018}, a path-following method without the barrier function has a \emph{globally superlinear rate} of convergence based on a few recent technical developments \citep{dvurechenskyImprovedGlobalPerformance2025}.
Note that neither scheme needs an auxiliary procedure to locate initial approximate centers of the corresponding homotopy paths, which has been a typical computational challenge for traditional IPMs \citep{nesterovInteriorPointPolynomialAlgorithms1994,yeInteriorPointAlgorithms1997,enagyNewInteriorpointAlgorithm2024}.
The convergence rates apply to the entire additive family of utility functions, which in particular includes the CES economy.

Finally, we discuss extensions of our approach to more general demand models. This includes the case when $(1)$ the players' utilities are linear (a limiting case of the above discussion), and $(2)$ the players are restricted by homogeneous linear constraints, such as demand models with aggregations, and the flow market. The complexity rates of the two interior-point schemes remain the same: linear and superlinear, respectively.
Numerical experiments are provided to illustrate the effectiveness of the proposed methods for large-scale problems.

\subsection{Organization}
The rest of this paper is structured as follows.
In \cref{sec.ipm}, we briefly go through our methodology.
In \cref{sec.logump}, we discuss the $\lump$, the calculus, and the SLC property.
After thorough discussions on $\lump$ and the calculus, our approximate schemes for computing the inexact steps are presented in \cref{sec.fast.approx}.
The convergence analyses of the two methods are presented in \cref{sec.conv}.
In \cref{sec.extend}, we discuss the potential of our approach to more general problems.
Preliminary experiments are presented in \cref{sec.numerical} and the conclusion (\Cref{sec.conclusion}) is given thereafter.
For the sake of clarity, we will leave the proofs in the appendix.

 \section{Preliminaries and Overview of the Approaches}\label{sec.ipm}

In this section, we briefly go through the elements centered around the best responses to construct the decentralized interior-point strategies for the market equilibrium.
Throughout the paper, we let $\cH^d$ be the set of continuous and homogeneous mappings of order $d$. Namely, if $u \in \cH^d$, then for any $\x \in \real_+^n, \alpha \in \real_+$, $u(\alpha \x) = \alpha^d u(\x)$; similar definition could be defined for functions of $\p$. For the domain $\cP$, we denote its relative interior by $\cP^\circ$ and its boundary by $\partial \cP$. We let a capital letter be the diagonal matrix with the lower-case letters on the diagonal, e.g., $\bP = \diag(\p)$. The notation $\|\cdot\|$ without any subscript represents the $\ell_2$-norm.

We focus on heterogeneous players with $\cH^{d_1}, ..., \cH^{d_m}$ utilities, and more specifically, the \emph{additively homogeneous utilities} $\cA^{k_i, r_i}$ for some real numbers $k_i, r_i$ (the concept will be made clear later in \Cref{eq.additive.utility}). This ubiquitous family includes the CES utilities. To start with, we introduce the Logarithmic Utility Maximization problem.

\subsection{The logarithmic utility maximization problem}
Instead of using the classical utility maximization problem (UMP), a pair of functions are introduced for each player $i$: the negative logarithmic utility $v_i(\x_i) = - \log(u_i(\x_i))$, and the \emph{dual function} $f_i(\p)$, connected via the Logarithmic Utility Maximization Problem ($\lump{}$):
\begin{subarr}{eq.logump}{rl}
    (\lump{}) \quad f_i(\p) := \max_{\x_i \in \cX} -v_i(\x_i), ~\st~ \inner{\p}{\x_i} \le w_i.
\end{subarr}
We will show in \Cref{lem.equiv.conj.purt} that the best response of $\lump{}$ is the same as that of the classical UMP \eqref{eq.ump}, and the gradient of $f_i(\p)$ is proportional to $\x_i(\p)$.
Using $\lump$, the potential function that we minimize for the heterogeneous case is easily extended from \eqref{eq.indirect.dual.eg} (cf. \cref{corr.potential}):
\begin{equation}\label{eq.potential.general}
    \min_{\p\in\cP} ~\varphi(\p) = \inner{\p}{\1} + \sum_{i\in\cI} \frac{w_i}{d_i} f_i(\p).
\end{equation}
As before, the gradient is the negative excess demand, $\nabla\varphi(\p) = \1 - \sum_{i\in\cI} \x_i(\p)$, and the stationary point of $\varphi$ corresponds to the equilibrium price (note again that each player has a $\cH^{d_i}$ utility).

The rationale of focusing on $(v_i, f_i)$ from $\lump$ instead of $(u_i, \pi_i)$ from UMP is its expressibility.
Since $\x_i(\p)$ is a solution to an optimization problem itself (either the UMP or $\lump{}$), to compute the high-order derivatives, including $\nabla \x_i(\p)$, via $(u_i, \pi_i)$, we would have to differentiate through formulation  \eqref{eq.ump}. This approach provides little economic insight and encounters heavy algebraic manipulations when investigating the functional properties (e.g., Lipschitzness of $\x_i(\p)$).
By comparison, $\lump$ completely relieves such a burden. We show that $(v_i, f_i)$ are Fenchel conjugate dual functions of each other up to some constants \citep{rockafellarVariationalAnalysis2009}.
Hence, both $v_i$ and $f_i$ are $d$-logarithmically homogeneous functions, allowing an arsenal of tools famous for self-scaled cones \citep{gulerBarrierFunctionsInterior1996,nesterovSelfscaledBarriersInteriorpoint1997}.
The calculus can be derived therein without pain.
For the $i^\text{th}$ player, $\x_i(\p) = -\frac{w}{d}\nabla f_i(\p)$, and in addition, computing $\nabla^{p}\x_i(\p)$, $p^{\text{th}}$-order derivatives of $\x_i$, shares the same task of computing $\nabla^{p+1} f_i(\p)$.
In fact, we only need the bidding vector (the distribution of money spent on the bundle):
\begin{equation}\label{eq.bidding.vector}
    \bgamma_i = \frac{\bP\x_i(\p)}{w_i}.
\end{equation}
As shown in \Cref{lem.pd.bestresp}, it encodes all we have to know for the computations of $(1)$ the derivatives of $v_i, f_i$ of any order, and $(2)$ the inverse operator of $\nabla \x_i(\p)$.

Next we discuss the Lipschitzness properties of the best response. Generally, $\x_i(\p)$ cannot be Lipschitz continuous with respect to the Euclidean distances.
Consider the simple example where $\cP = \mathbb{R}_+^n$. If the price $\p$ approaches the boundary of the domain (i.e., the goods become almost free of charge), the best response $\x_i(\p)$ for any player becomes unbounded.  Essentially, $\x_i(\p)$ itself serves as a \enquote{implicit barrier function} for the price space $\cP$ \citep{fiaccoNonlinearProgrammingSequential1968}.   We could summarize the properties of the implicit barrier specifically as follows:
\begin{enumerate}[label=$(\alph*)$, leftmargin=*]
    \item The best response $\x_i(\p)$ (and $\nabla f_i(\p)$) diverges as $\p$ approaches the boundary $\partial\cP$, and so does its Jacobian $\nabla \x_i(\p)$ (and $\nabla^2 f_i(\p)$), cf. \eqref{eq.x.jac};
    \item Under an appropriate affine scaling, $\x_i(\p)$, the Jacobian $\nabla \x_i(\p)$, and the second-order derivatives of $\x_i(\p)$ (and $\nabla^3 f_i(\p)$) remain uniformly bounded, cf. \eqref{eq.bounded.directdev.dual};
    \item The perturbation of $\x_i(\p)$ in $\p$ has the Scaled Lipschitz Continuity (SLC), cf. \eqref{eq.br.localstable}.
\end{enumerate}

Because of the separability, the potential function $\varphi$ in \eqref{eq.potential.general} inherits these properties as an agglomeration of $f_1, ..., f_m$. Based on these facts, we develop interior-point strategies to find
\begin{equation}\label{eq.fosp}
    \epsilon-\textbf{approximate equilibrium price vector}: \quad \p \in \cP, \quad \|\nabla\varphi(\p)\|_\infty \le \epsilon
\end{equation}
for a small enough $\epsilon > 0$.
At some $\p$, we first collect the best responses of all players and compute $$\nabla \varphi = \1 - \sum_{i\in\cI} \x_i(\p), \quad \nabla^2\varphi = -\sum_{i\in\cI} \nabla \x_i(\p)$$ as the summation of the individual ones (simple aggregations of $\bgamma_1, ..., \bgamma_m$).
Each step of the IPMs calculates an update of the price vector $\p$ by solving one Newton-type system, using the Hessian operator $\nabla^2\varphi(\p)$ and the gradient $\nabla\varphi(\p)$.
With these analyses on hand, our development of interior-point strategies consists of two parts as described in the following two subsections.
\subsection{Solving the Newton-type systems efficiently}
The first part involves developing efficient schemes to solve the Newton-type systems (\Cref{sec.fast.approx}). Specifically, we focus on the linear system of $\q \in \real^n$ in the form of
\begin{equation}\label{eq.original.newton}
    (\nabla^2\varphi(\p) + \bG) \q = \r,
\end{equation}
where $\bG \succeq 0$ is some regularization matrix, the right-hand side $\r$ is computed from $\nabla\varphi(\p)$ (it could be just $\nabla\varphi(\p)$ itself). Then, the price vector $\p_+ = \p + \alpha \q$ is updated by some stepsize $\alpha > 0$ while keeping $\p_+ \in \cP^\circ$.
The Hessian operator, if affinely scaled by $\bP$, is uniformly bounded for any $\p \in \cP^\circ$:
\begin{eqalign}{eq.affine.hessian}
    \bH(\p) := \bP\nabla^2 \varphi(\p)\bP.
\end{eqalign}
An equivalent form of \eqref{eq.original.newton} using scaled direction $\q = \bP\d$ is
\begin{equation}\label{eq.inexact.newton}
    (\bH(\p) + \bP\bG\bP)\d = \bP\r,
\end{equation}
For small markets, we can directly use the exact $\bH(\p)$.  Otherwise, we solve \eqref{eq.inexact.newton} inexactly. This task is made easy via the bidding vectors $\bgamma_1, ..., \bgamma_m$ in two different ways.

Because of the uniform boundedness of $\bH(\p)$, we could adopt some approximate Hessian operator $\widetilde{\bH}(\p)$ and require
\begin{equation}\label{eq.def.epsH}
    \|\widetilde{\bH}(\p) - \bH(\p)\| < \epsilon_H,
\end{equation}
where $\epsilon_H \ge 0$ is a predefined accuracy parameter.
This requirement can be satisfied using a variety of techniques in randomized linear algebra \citep{drineasRandNLARandomizedNumerical2016}.
For our purposes, we present a simple invertible DR1 approximation by averaging the ``bidding vectors'' $\bgamma_1, ..., \bgamma_m$. We show under mild conditions in the large-market regimes that the error diminishes in high probability (see \cref{sec.fast.approx.dr1}).

Another possibility is explored for Krylov subspace methods, which starts with an initial guess $\d_0$ and maintains an evolving subspace $\cV$. At each iteration, it updates the direction $\d$ in the subspace spanned by $\cV$, until the residual of \eqref{eq.inexact.newton} is small enough:
\begin{equation}\label{eq.inexact.krylov}
    \left\|(\bH(\p) + \bP\bG\bP)\d - \bP\r \right\| \le \epsilon_K \|\d\|
\end{equation}
Otherwise, it increases $\dim(\cV)$ by one. In this case, $\widetilde{\bH}(\p)$ is defined as the projected matrix of $\bH(\p)$ onto $\cV$.
We could require the approximation to be good enough along the direction $\d$ in the following sense:
\begin{eqalign}{eq.def.epsH.direction}
    \|(\widetilde{\bH}(\p) - \bH(\p))\d\|
    &\le \epsilon_H \|\d\|.
\end{eqalign}
Theoretically, for any $\epsilon_K, \epsilon_H > 0$, we could enlarge the dimension of $\cV$ until \eqref{eq.inexact.krylov}-\eqref{eq.def.epsH.direction} hold, and they hold trivially when $\dim(\cV) = n$.
For example, when using the Conjugate Gradient (CG) method, it terminates finitely while having a linear convergence rate, which depends on the condition number $\kappa_2(\bH(\p)) = \tfrac{\lambda_{\max}(\bH(\p))}{\lambda_{\min}(\bH(\p))}$ \citep{golubMatrixComputations2013}. To reduce the Krylov iterations, many heuristics have been proposed for finding good preconditioners for $\bH(\p)$, i.e., positive definite matrices that re-scale \eqref{eq.inexact.newton} in a different metric and reduce the condition number.
Among them, the diagonal matrices are prevalent because of the simplicity.
However, finding an optimal diagonal preconditioner to achieve the maximum reduction in condition number is a nontrivial quasi-convex optimization problem and can be found by solving a positive semidefinite program \citep{gaoScalableApproximateOptimal2024,quOptimalDiagonalPreconditioning2025}. Perhaps surprisingly, we show that the optimal diagonal preconditioner is analytic and is the row-summation of the Hessian operator, which again, is certain aggregation of $\bgamma_1, ..., \bgamma_m$.

\subsection{The interior-point methods}
The above constructions lay the groundworks for lightweight iterations in second-order t\^atonnement processes, i.e., the IPMs.
Next, we introduce two IPMs using different choices of Newton systems and homotopy paths.

\subsubsection*{Second-order T\^atonnement 1: \logbar{}}\label{sec.logbar}
The most natural approach constructs a sequence of logarithmic barrier models as follows:
\begin{equation}\label{eq.log barrier varphi}
    \varphi_{\mu}(\p) = \varphi(\p) - \mu \inner{\1}{\log(\p)}, \ \mu >  0.
\end{equation}
where $\varphi$ is the potential function described in \eqref{eq.potential.general} and we gradually decrease $\mu$ to zero.

In \Cref{alg.logbarrier}, at the initialization (\cref{line.logbar.ac}), we begin with a sufficiently large $\mu_0 > 0$, and choose $\p_0 \in \cP^\circ$ that is an approximate Analytic Center (AC), a point in the neighborhood of the central path for some constant $Q \in (0, \tfrac{1}{2})$:
\begin{equation}\label{eq.central path}
    \cC(\mu,Q) = \left\{\p > 0: \tfrac{\|\bP \nabla \varphi(\p) - \mu \1\|}{\mu} \le Q \right\}.
\end{equation}
Note that finding the AC is highly nontrivial for general convex optimization \citep{yeInteriorPointAlgorithms1997}. However, for the market equilibrium problem, we can provide an initial approximate AC and the corresponding $\mu_0$ from an economic fact (see \Cref{lem.logbar.init}).
\begin{center}
    \begin{minipage}[h]{0.9\textwidth}
        \begin{algorithm}[H]
        \renewcommand{\algorithmcfname}{Second-order T\^atonnement}
            \caption{The Logarithmic Barrier Method (\logbar{})}\label{alg.logbarrier}
            \small
            \KwIn{$k=0$, $\sigma \in (0, 1)$, $Q \in (0, \tfrac{1}{2})$}
            \tcp{find the approximate analytical center}
            \label{line.logbar.ac} Find some $\mu_0 > 0$ and $\p_0 \in \cP^\circ$, such that
            $\p_0 \in \cC(\mu_0,Q)$\;
            \tcp{main iterates}\label{line.logbar.main}
\For{$k = 0, 1, \ldots, K-1$}
            {
                Update $\mu_{k+1} = \sigma \cdot \mu_k$\;
                Solve $\d_k$ from \eqref{eq.inexact.logbar} at $\p_{k}$ with $\mu = \mu_{k+1}$\;
                Update $\p_{k+1} = \p_{k} + \bP_{k}\d_k$\; \label{line.logbar.newton.main}
            }
            \KwOut{$\p_K$}
        \end{algorithm}
    \end{minipage}
\end{center}
In the main phase (\cref{line.logbar.main}),
any update at some $(\mu,\p)$ (\cref{line.logbar.newton.main}) is computed by a realization of \eqref{eq.inexact.newton} as follows:
\begin{equation}\label{eq.inexact.logbar}
    (\widetilde{\bH}(\p) + \mu \bI)\d = - (\bP\nabla \varphi(\p) - \mu \1).
\end{equation}
We linearly contract $\mu$ via a shrinkage factor $\sigma \in (0, 1)$. After a single inexact Newton step, we can show that the resulting iterate comes back to the neighborhood defined by an updated $\mu$. We show that this scheme has an overall linear rate of convergence (\Cref{thm.complexity.logbarrier}).

\subsubsection*{Second-order T\^atonnement 2: \pathf{}}
The second strategy is an adaptation of the path-following method \cite[Eq. 29]{dvurechenskyImprovedGlobalPerformance2025}, where the requirement for the exact Hessian is relaxed. Briefly, the method constructs a homotopy path of $t$ without the logarithmic barrier:
\begin{equation}\label{eq.pathf.varphi}
    \varphi_t(\p) = \varphi(\p) - t\inner{\nabla\varphi(\p_0)}{\p}, \quad t \in [0,1].
\end{equation}
The measure of progress is based on the Newton decrement,
\begin{equation}
    \lambda(\p) = \left(\nabla^2 \varphi(\p)^\inv[\nabla \varphi(\p), \nabla \varphi(\p)]\right)^\invhalf = \|\nabla \varphi(\p)\|_{\nabla^2 \varphi(\p)}^*.
\end{equation}
Since we do not invert the exact Hessian, the local norm $\|\cdot\|_{\nabla^2 \varphi(\p)}$ is approximated from $\widetilde{\bH}$. That is, for any $\q$ we use,
$\left\|\bP\q\right\|^*_{\widetilde{\bH}(\p)} \approx {\left(\left(\bP\nabla^2 \varphi(\p)\bP\right)^\inv[\bP\q,\bP\q]\right)}^{\invhalf} = \left\|\q\right\|_{\nabla^2 \varphi(\p)}^*.
$
A corresponding concept of AC means that, for some $t$, the iterate $\p$ satisfies:
\begin{equation}\label{eq.pathf.centered}
    \cC(t) = \left\{\p \in \cP^\circ: \|\bP\nabla \varphi(\p) - t\bP\nabla\varphi(\p_0)\|^*_{\widetilde{\bH}(\p)} \le \tfrac{\beta}{C_\varphi}\right\}, ~\beta \in (0,1],
\end{equation}
where $C_\varphi$ is the self-concordance constant of $\varphi(\p)$ specified later.
We present the method in \Cref{alg.pathfollowing}.
Compared to the original scheme (\cite{nesterovLecturesConvexOptimization2018,dvurechenskyImprovedGlobalPerformance2025}), we replace the Newton decrement $\lambda(\p)$ by the inexact local norm $\|\cdot\|_{\widetilde{\bH}(\p)}$. Each step $\d$ (\cref{line.pathf.newton.main}) is computed by the following inexact Newton step:
\begin{equation}\label{eq.inexact.pathf}
    \widetilde{\bH}(\p)\d = - \bP\left(\nabla \varphi_t(\p) - t\nabla \varphi(\p_0)\right).
\end{equation}
And alternatively, letting $\q = \bP\d$, in view of the original Hessian, it is equivalent to say
\begin{equation}\label{eq.inexact.pathf.original}
    \widetilde{\nabla^2\varphi}(\p) \q = - \left(\nabla \varphi_t(\p) - t\nabla \varphi(\p_0)\right).
\end{equation}
The initialization is also simple: pick $t_0 = 1$, and select $\p_0 \in \cP^\circ$ and we automatically have $\p_0 \in \cC(t_0)$.
As long as $\gamma,\beta$ are selected properly, the new iterate $\p \in \cC(t)$ remains centered.
In \eqref{eq.pathf.updatet}, the parameter $t$ is updated until it hits 0, and once it does, the subsequent iterates converge at a quadratic rate. We show that this scheme has an overall superlinear rate of convergence (\Cref{thm.convergence.rate.pathf}).
\begin{center}
    \begin{minipage}[h]{0.9\textwidth}
        \begin{algorithm}[H]
        \renewcommand{\algorithmcfname}{Second-order T\^atonnement}
            \caption{The Path-Following Method (\pathf{})}\label{alg.pathfollowing}
            \small
            \KwIn{$k=0$, $t_0 = 1$, $\p_0 \in \cP^\circ$, $\beta \in (0,1]$, $\gamma \in (0,1).$}
            \tcp{main iterates}\label{line.pathf.main}
            \For{$k = 0, 1, \ldots, K-1$}
            {
                Update $t_{k+1}$ according to
                \begin{equation}\label{eq.pathf.updatet}
                    t_{k+1} = \max\left\{t_k - \frac{\gamma}{C_\varphi\|\bP\nabla \varphi(\p_0)\|_{\widetilde{\bH}(\p_k)}}, 0\right\};
                \end{equation}

                Solve $\d_k$ from \eqref{eq.inexact.pathf} at $\p_{k}$ with $t = t_{k+1}$\;\label{line.pathf.newton.main}
                Update $\p_{k+1} = \p_{k} + \bP_{k}\d_k$\;
            }
            \KwOut{$\p_K$}
        \end{algorithm}
    \end{minipage}
\end{center}

\section{Logarithmic Utility Maximization and the Scaled Lipschitz Continuity}\label{sec.logump}
In this section, we start to formally present our dicussion. First, we discuss the $\lump$ and show how it enables our analysis on the calculus and Lipschitz properties. For the sake of notation, the subscript $i$ is dropped in this section if no confusion arises.
Recall that $v(\x) = -\log(u(\x))$ and
\begin{align*}
    f(\p) := \max_{\x \in \cX} -v(\x)
    , ~\st~ \inner{\p}{\x} \le w.
\end{align*}
Suppose $u \in \cH^d, d > 0$, then for any $\alpha> 0$, $\x \in \cX \subseteq \real_+^n$,
\begin{eqalign}{eq.logump.def}
    v(\alpha\x) & = -\log(u(\alpha\x)) = v(\x) - d\log(\alpha).
\end{eqalign}
This implies that $v$ is a \emph{$d$-logarithmically homogeneous} function.
The canonical logarithmic barrier function will be repetitively used:
\begin{eqalign}{eq.logbarrier}
    \cB(\x) = -\sum_{j\in\cJ} \log(x^{(j)}), \quad \nabla^2 \cB(\x) = \bX^{-2}, \quad \x \in \real_+^n.
\end{eqalign}
We could write $(\nabla^2\cB(\x))^{0.5} = \bX^\inv$. According to standard interior-point nomenclature, define
\begin{subarr}{eq.localnorm}{rl}
    \h \in \real^n, ~ \|\h\|_{\x} &= \|\bX^\inv\h\| = \inner{\h}{\nabla^2 \cB(\x)\h}^{\tfrac{1}{2}},  \\
    \q \in \real^n, ~\|\q\|^*_{\x} &= \|\bX\q\| = \inner{\q}{(\nabla^2 \cB(\x))^\inv\q}^{\tfrac{1}{2}}.
\end{subarr}
as the primal local norm of $\h$ and dual local norm of $\q$ at $\x$, respectively.
The Dikin ellipsoid centered at $\x$ with radius $\rho \in (0,1)$ is defined as
\begin{eqalign}{eq.dikin}
    \cD(\x, \rho) = \left\{ \x' \in \real^n ~\middle|~ \|\x' - \x\|_{\x} \le \rho < 1 \right\}.
\end{eqalign}
Whenever $\x' \in \cD(\x, \rho)$, we would imply that $\x' > 0$. Similar definitions can be made when referring to the prices.
In the sequel, we make the following assumption.
\begin{assumption}\label{ass.logump.unique.bestresp}
    $u$ is differentiable and strictly concave, and the best response $\x(\p)$ is unique.
\end{assumption}
Note that the strict concavity of $u$ guarantees the uniqueness of $\x(\p)$, then the gradient, and further high-order derivatives exist.
Otherwise, only the price, not the allocations, can be determined by the t\^atonnement.\endnote{Just like linear programming, merely providing an optimal dual vector (the price $\p$) shall not make it any easier to find the optimal primal solution (the bundles $\x_1, ..., \x_m$) \citep{megiddoFindingPrimalDualoptimal1991}.}  A relaxation of it is discussed later for linear utility functions (cf. \Cref{sec.extend.linear}), where we add a perturbation to ensure the uniqueness of the best response.
\subsection{Relationships between $\lump$ and UMP}
Consider the Legendre-Fenchel transform of $v$ in the same fashion of \cite{nesterovSelfscaledBarriersInteriorpoint1997}:
\begin{eqalign}{eq.legendrefenchel}
    v_*(\p) &= \max_{\x \in \cX} \left\{ -\inner{\p}{\x} - v(\x) \right\}.
\end{eqalign}
The following lemma connects $v_*$ and $f$ and asserts that the best response $\x(\p)$ of $\lump{}$ coincides with that of the standard UMP.
\begin{lemma}\label{lem.equiv.conj.purt}
    The conjugate function $v_*$ is a $d$-logarithmically homogeneous function. The dual function $f$ satisfies: $f(\p) = v_*(\p) + d + d\log(w) - d\log(d)$.
    That is, $f$ is also $d$-logarithmically homogeneous. Furthermore, suppose that $\x(\p)$ is the best-response function of $\lump{}$ \eqref{eq.logump}, then it is also the best response of UMP \eqref{eq.ump}. Besides,
    \begin{eqalign}{eq.f.gradient}
        \nabla f(\p) = -\frac{d}{w}\x(\p).
    \end{eqalign}
\end{lemma}
As the dual function $f(\p)$ is equivalent to the conjugate function $v_*(\p)$ up to some constant, we already have an extensive array of results from logarithmically homogeneous functions \citep{nesterovInteriorPointPolynomialAlgorithms1994,nesterovSelfscaledBarriersInteriorpoint1997}. These tools were developed in the 1990s for the analysis of self-concordant barriers of self-scaled cones. We collect a few useful results in \Cref{appendix.sec.logump}. However, one should be careful since by far we cannot state that $v, v_*, f$ are self-concordant barriers (e.g., \cite[Definition 2.1.1]{nesterovInteriorPointPolynomialAlgorithms1994}). Coming back to the market equilibrium problem, \cref{lem.equiv.conj.purt} implies that we can extend the Eisenberg-Gale program for hereogeneous $\cH^{d_i}$ functions without the restriction of $d_i=1$ for all $i\in\cI$.
\begin{corollary}\label{corr.potential}
    Suppose that in the market, for each player $i$ the utility function $u_i \in \cH^{d_i}$.
    Let the potential function $\varphi(\p)$ be as follows,
    \begin{equation}
        \varphi(\p) = \inner{\p}{\1} + \sum_{i\in\cI} \frac{w_i}{d_i} f_i(\p),
    \end{equation}
    then the equilibrium price can be found from its first-order stationary points.
\end{corollary}

\subsection{Additive homogeneous utility functions}
Let us consider a common family of utility functions composed of additive elements.
We say a utility function $u$ is additively homogeneous, $u\in\cA^{k,r}$,
if it is the power of the summation of nonnegative $\cH^r$ functions for $k\neq 0, r\neq 0, d=k\times r > 0$, viz.
\begin{equation}\label{eq.additive.utility}
    u(\x) = \left(\smallsum_{j\in\cJ} \theta_j(x^{(j)})\right)^k, \quad\theta_j:\cX^{(j)}\to\real_+,\theta_j \in \cH^{r},\forall j\in\cJ.
\end{equation}
That is, $u(\x) \in \cH^{d}$.
For conciseness, let
$$\theta^{(j)} := \theta_j(x^{(j)}), ~\btheta = [\theta^{(1)}, \cdots, \theta^{(n)}]= [\theta_1(x^{(1)}), \cdots, \theta_n(x^{(n)})], $$ and so $u(\x) = \inner{\btheta}{\1}^k$.
Hence, it is easy to see that the $\cA^{k,r}$ family is a generalization of the constant-elasticity-of-substitution (CES) utility functions (see \Cref{ex.ces}).
Its negative logarithm reads:
$
    v(\x) = -k \log (\sum_{j\in\cJ} \theta^{(j)}).
$
If $k=1$, we could say that they are separable \citep{chenSettlingComplexityArrowDebreu2009}, the maximization for which belongs to a class of entropy optimization problems \citep{potraQuadraticallyConvergentPolynomial1993}.
The following theorem unveils that after affine scaling, the derivatives (of any order) are \enquote{parallel} to a vector $\bgamma$ in the unit simplex $\Delta_n$ (we later show it is the distribution of money spent on the goods). Besides, the concavity of $\cA^{k,r}$ family can be ensured by restricting the spectra of $(k,r)$.
\begin{lemma}\label{lem.additive.diagonal}
    If $u(\x) \in \cA^{k,r}$, then there exists some $\bgamma \in \Delta_n$, such that the following holds:
    \begin{eqalign}{eq.parallel.theta}
        \frac{\bX\nabla \btheta}{r\inner{\btheta}{\1}}
        = \frac{\bX\nabla^2 \btheta[\x, \cdot]}{r(r-1)\inner{\btheta}{\1}}
        = \frac{\bX\nabla^3 \btheta[\x, \x, \cdot]}{r(r-1)(r-2)\inner{\btheta}{\1}}
        = \frac{\btheta}{\inner{\btheta}{\1}} = \bgamma.
    \end{eqalign}
    Furthermore, if
    \begin{eqalign}{eq.concave.u}
        r \in (0,1), k\in (0, \tfrac{1}{r}] \quad \text{or} \quad r < 0, k \in [\tfrac{1}{r}, 0)
    \end{eqalign}
    then $u(\x)$ is concave ($d > 0$).
\end{lemma}
In light of the above result, we write the gradient and Hessian of $v(\x)$ for $\x \in \cX$:
\begin{align*}
    \nabla v(\x) = -k\frac{\nabla \btheta}{\inner{\btheta}{\1}},
    \qquad \nabla^2 v(\x) & = -k\frac{\nabla^2 \btheta}{\inner{\btheta}{\1}}
    + k\left(\frac{\nabla \btheta}{\inner{\btheta}{\1}}\right)\left(\frac{\nabla \btheta}{\inner{\btheta}{\1}}\right)^{\top}.
\end{align*}
Apparently, the primal Hessian operator $\nabla^2 v(\x)$ has the \emph{Diagonal+Rank-One (DR1)} structure: from \eqref{eq.parallel.theta}, it is composed of a diagonal and a rank-one matrix.
At the same time, the vector $\bgamma \in \Delta_n$ represents the \emph{bidding vector}. We present the following result.
\begin{lemma}\label{lem.pd.bestresp}
    Let $\bgamma$ be defined as in \eqref{eq.parallel.theta}, and all quantities be defined at the best response $\x(\p)$, then the following results hold:
    \begin{equation}
        \bgamma = \frac{\bX(\p)\p}{w} \in \Delta_n.
    \end{equation}
    That is, $\bgamma$ represents the distribution of money spent on the goods -- the bidding vector.
    The primal Hessian operator $\nabla^2 v(\x(\p))$ is DR1, viz.
    \begin{equation}
        \nabla^2 v(\x(\p)) = d(1-r)\frac{1}{w^2} \bP \bGamma^\inv\bP +\frac{dr}{w^2} \p\p^{\top}.
    \end{equation}
    The dual Hessian operator $\nabla^2 f(\p)$ is DR1, viz.
    \begin{equation}\label{eq.f.hess}
        \nabla^2 f(\p) = \frac{d^2}{w^2}\left(\nabla^2v(\x(\p))\right)^\inv = \frac{d}{(1-r)} \bP^\inv \left(\bGamma - r \bgamma\bgamma^{\top}\right)\bP^\inv.
    \end{equation}
    And so is the Jacobian of the best-response mapping, viz.
    \begin{equation}\label{eq.x.jac}
        \nabla \x(\p) = - \frac{w}{d}\nabla^2 f(\p) = - \frac{w}{(1-r)} \bP^\inv \left(\bGamma - r \bgamma\bgamma^{\top}\right)\bP^\inv.
    \end{equation}
\end{lemma}
The Hessian operators possess certain symmetry due to homogeneity: $\nabla^2 v(\x)$ (resp. $\nabla^2 f(\p))$ has been \enquote{affine-scaled} by $\bP$ (resp. $\bP^\inv$).
Because $\bgamma$ is always in the unit simplex, it is apparent that $\|\x(\p)\|, \|\nabla \x(\p)\|$ are unbounded if $\p \to 0$.
Note that no extra information is needed for the second-order derivatives of $f(\p)$ and $v(\x)$.  Even if $r$ is unknown, using the responses at two different prices, $r$ could be estimated.

\subsection{Scaled Lipschitz continuity}
For the additive family, \cref{lem.additive.diagonal} states the \enquote{off-diagonal} elements are zero for the high-order derivatives of $\btheta$.
If we take a direction $\h\in\real^n$, we could observe the following properties of $v(\x)$.
\begin{lemma}[Scaled boundedness of directional derivatives, primal space]\label{lem.bounded.directdev}
    If $u(\x) \in \cA^{k,r}$ satisfying \eqref{eq.concave.u}, then for any $\h\in \real^n$, the following holds
    \begin{subarr}{eq.bounded.directdev}{rl}
        \left|\nabla v(\x)[\h]\right| &\le d \|\h\|_{\x}, \\
        \left|\nabla^2v(\x)[\h,\h]\right| &\le \max\{d, d(1-r)\} \|\h\|_{\x}^2, \\
        \left|\nabla^3v(\x)[\h,\h,\h]\right| &\le T_u \|\h\|_{\x}^3,
        \quad T_u := \max\{2d, d(6r^2-6r+2)\}.
    \end{subarr}
    Besides, $\tfrac{1}{d}\nabla^2 v(\x)[\h,\h] \ge \|\tfrac{1}{d}\nabla v(\x)[\h]\|^2.$
\end{lemma}
Using everywhere boundedness of directional derivatives, the Scaled Lipschitz Continuity (SLC) of $v$ can be specified.
\begin{theorem}[SLC in the primal space]\label{thm.taylor.v}
    If $u(\x) \in \cA^{k,r}$ satisfying \eqref{eq.concave.u}, then for any $\h\in\real^n$, and $\rho := \|\h\|_{\x} < 1$, we have
    \begin{eqalign}{eq.taylor.v}
        \left| v(\x+\h) - v(\x) - \nabla v(\x)[\h] - \tfrac{1}{2}\nabla^2 v(\x)[\h,\h] \right|\le -\tfrac{1}{2}T_u \left(\rho +  \tfrac{1}{2}\rho^2 +  \log(1 - \rho)\right),
    \end{eqalign}
    where $T_u$ is defined as in \eqref{eq.bounded.directdev}.
    Hence,
    \begin{eqalign}{eq.taylor.v.more}
        \left|v(\x+\h) -  v(\x) - \nabla v(\x)[\h] - \tfrac{1}{2}\nabla^2 v(\x)[\h,\h]\right|
        \le \tfrac{T_u \rho^3}{6(1-\rho)}.
    \end{eqalign}
\end{theorem}
Similarly, we are interested in the local perturbation of $f(\p)$ in the space of prices. By using the aforementioned arguments, we could conclude a few results symmetric to those in the primal space.
\begin{lemma}[Scaled boundedness of directional derivatives, dual space]\label{lem.bounded.directdev.dual}
    If $u(\x) \in \cA^{k,r}$ satisfying \eqref{eq.concave.u}, then for any $\q\in \real^n$, the following holds for $f(\p)$:
    \begin{subarr}{eq.bounded.directdev.dual}{rl}
        \left|\nabla f(\p)[\q]\right| &\le d \|\q\|_{\p}, \\
        \left|\nabla^2f(\p)[\q,\q]\right| &\le \tfrac{d(1+|r|)}{1-r} \|\q\|_{\p}^2, \\
        \left|\nabla^3f(\p)[\q,\q,\q]\right| &\le T_f \|\q\|_{\p}^3
    \end{subarr}
    where
    \begin{equation}\label{eq.f.tf.def}
        T_f = \max\left\{\tfrac{6d}{(1-r)^2}, 2d\right\}.
    \end{equation}
\end{lemma}
\begin{theorem}[SLC in the dual space]\label{thm.taylor.f}
    Suppose $u(\x) \in \cA^{k,r}$ satisfies \eqref{eq.concave.u}. For any $\p \in \real_+^n, \q \in \real^n$,
    if $\varrho = \|\q\|_{\p} < 1$, then
    \begin{eqalign}{eq.taylor.f}
        |f(\p+\q) - f(\p) - \inner{\nabla f(\p)}{\q} - \tfrac{1}{2}\nabla^2 f(\p)[\q,\q]|
        \le - \tfrac{T_f}{2} \left(\varrho +  \tfrac{1}{2}\varrho^2 +  \log(1 - \varrho)\right).
    \end{eqalign}
    where $T_f$ is defined as in \eqref{eq.f.tf.def}.
    As well,
    \begin{eqalign}{eq.taylor.f.more}
        \left|f(\p+\q) -  f(\p) - \inner{\nabla f(\p)}{\q} - \tfrac{1}{2}\nabla^2 f(\p)[\q,\q]\right|
        \le \tfrac{\varrho^3T_f}{6(1-\varrho)}.
    \end{eqalign}
\end{theorem}
Specifically, since $\nabla f(\p) = -\frac{d}{w}\x(\p)$, we could summarize the following result in the flavor of $\x$.
\begin{theorem}[SLC for the best-response mapping]\label{thm.taylor.x}
    Suppose $u(\x) \in \cA^{k,r}$ satisfies \eqref{eq.concave.u}. For any $\p \in \real_+^n, \q \in \real^n$,
    if $\varrho = \|\q\|_{\p} < 1$, then
    \begin{eqalign}{eq.br.localstable}
        \|\bP(\x(\p+\q) - \x(\p) - \nabla\x(\p)[\q])\| \le \tfrac{w}{d}\tfrac{T_f\varrho^2}{2(1-\varrho)}.
    \end{eqalign}
    where $T_f$ is defined as \eqref{eq.f.tf.def}.
\end{theorem}

While SLC looks quite similar to a few classical conditions in the interior-point literature \citep{monteiroExtensionKarmarkarType1990,kortanekPolynomialBarrierAlgorithm1993,denhertogSufficientConditionSelfconcordance1995,nesterovInteriorPointPolynomialAlgorithms1994,andersenHomogeneousAlgorithmMonotone1999}, this property and its application in the exchange economy appear to be new and may have independent interest.

\subsection{Active-set boundedness and self-concordance}
Generally speaking, $v, f$ induced from the $\cA^{k,r}$ family is not necessarily self-concordant \citep[definition 2.1]{nesterovInteriorPointPolynomialAlgorithms1994}, which requires a bounded ratio between $2^{\text{nd}}$ and $3^{\text{rd}}$-order directional derivatives (see $\cR(\x)$ defined below).
\begin{lemma}[Dominance between $2^{\text{nd}}$ and $3^{\text{rd}}$-order directional derivatives]\label{lem.dom.23}
    If $u(\x) \in \cA^{k,r}$ satisfying \eqref{eq.concave.u}, then for any $\h\in\real^n$, we have
    \begin{eqalign}{eq.dom.23}
        \cR(\x) := \frac{|\tfrac{1}{d}\nabla^3 v(\x)[\h,\h,\h]|}{\left(\tfrac{1}{d}\nabla^2 v(\x)[\h,\h]\right)^{1.5}} \le \cO(\|[\bgamma]^{-0.5}\|_{\infty}).
    \end{eqalign}
    Let $B(\x)$ be the \emph{active set}:
    \begin{equation}\label{eq.activebasis}
        B(\x) = \{j\in\cJ: \gamma^{(j)} > 0\} = \{j\in\cJ: \theta_j(x^{(j)}) > 0\}.
    \end{equation}
    Then $\cR(\x) \le \cO\left(\underset{j\in B(\x)}{\max}\sqrt{\tfrac{1}{\gamma^{(j)}}}\right).$
\end{lemma}
In general, very few words can be said about the smallest element of $\bgamma$, hence, no uniformly finite upper bound can be established.
The boundedness of this ratio means whenever a good is interested, the share of money spent on it should not be too small.

If we assume that the quantity is uniformly bounded in our price-adjusting procedure, then the logarithmic utility $v$ acts as the self-concordant barrier function with $\cX$ as the natural domain. By Fenchel conjugate duality, $f$ is also self-concordant in the dual space. This assumption holds mildly for the CES utility functions (cf. \Cref{ex.ces}).
\begin{assumption}[Active-set boundedness]\label{asm.bounded.theta}
    There exists some $\kappa_u \in (0, \infty)$ such that the following quantities are bounded for $\x \in \cX$:
    $
        \kappa_u(\x) = \underset{j\in B(\x)}{\max}\{(\gamma^{(j)})^\inv\}
        \le \kappa_u.
    $
\end{assumption}
\begin{corollary}[Primal-dual self-concordance under active-set boundedness]\label{cor.selfconcordant.theta}
    If \cref{asm.bounded.theta} holds, then for any $\h\in\mathbb{R}^n$, if $\|\h\|_{\x} \le 1$, then
    \begin{eqalign}{eq.lhscb.param}
        \left|\nabla^3v(\x)[\h,\h,\h]\right| \le C_v \left(\nabla^2 v(\x)[\h,\h]\right)^{\tfrac{3}{2}},
    \end{eqalign}
    $C_v = \kappa_u^3 \max\{2d^{-\tfrac{1}{2}}, d^{-\tfrac{1}{2}}(6r^2-6r+2)\}.$
    Hence, $v(\x)$ is self-concordant on $\cX$.
    Besides, the dual function $f(\p)$ is also self-concordant on $\cP$ with the same coefficient: $C_f = C_v$.
\end{corollary}

\begin{example}[CES economy]\label{ex.ces}
    The utility function is in the family of constant elasticity of substitution (CES) if:
    \begin{equation}
        u(\x) = \left[\smallsum_{j\in\cJ} c^{(j)} (x^{(j)})^{\rho}\right]^{\frac{1}{{\rho}}} = \inner{\c}{\x^{\rho}}^{\frac{1}{{\rho}}}, {\rho} \in (-\infty,1), \c \in \real^n_+.
    \end{equation}
    where we omit linear and Leontief utilities.
Generally, $\delta = \frac{1}{1-{\rho}} \in (0, \infty)$ is called the \emph{elasticity of substitution}.
Clearly, $u(\x) \in \cA^{k,r}$ with $k = \frac{1}{\rho}, r = \rho$ and $d=1$, and hence $u(\x)$ is concave.
    In CES economy, the best response is, $\x(\p) = w\tfrac{\left[\bP^\inv\c\right]^{\frac{1}{1-{\rho}}}}{\inner{\c^{\frac{1}{1-{\rho}}}}{\p^{-\frac{\rho}{1-{\rho}}}}}.$ Indeed,
    \begin{equation}
        \forall j \in \cJ, \theta^{(j)} = (c^{(j)})^\frac{1}{1-\rho}(p^{(j)})^{-\frac{\rho}{1-\rho}} \quad
        \text{and} \quad
        \kappa_u(\x(\p)) = \tfrac{\sum_{j\in B(\x)} (c^{(j)})^\frac{1}{1-\rho}(p^{(j)})^{-\frac{\rho}{1-\rho}}}
        {\underset{j\in B(\x)}{\min}\{(c^{(j)})^\frac{1}{1-\rho}(p^{(j)})^{-\frac{\rho}{1-\rho}}\}}.
    \end{equation}
\end{example}
For simplicity, we assume that $\rho > 0$.
Clearly, $B(\x) \subseteq B = \{j\in\cJ: c_j > 0\}$.
In a realistic price updating mechanism, we usually keep bounded utilities and prices $u(\x) \le D_u, \|\p\|_\infty \le D_{\cP}$ (since the potential function is bounded), then $\forall j, \theta^{(j)} = (c^{(j)})^\frac{1}{1-\rho}(p^{(j)})^{-\frac{\rho}{1-\rho}} \ge (c^{(j)})^\frac{1}{1-\rho}D_{\cP}^{-\frac{\rho}{1-\rho}}$,
and $\sum_{j\in B(\x)} \theta^{(j)} \le D_u^\rho$.
In this view, it holds that
\begin{align*}
    \kappa_u(\x) = \tfrac{D_u^\rho}{\underset{j\in B(\x)}{\min}\{(c^{(j)})^\frac{1}{1-\rho}(p^{(j)})^{-\frac{\rho}{1-\rho}}\}}
    \le \tfrac{D_u^\rho D_{\cP}^{\frac{\rho}{1-\rho}}}{\underset{j\in B}{\min}\{c^{(j)}\}^\frac{1}{1-\rho}}.
\end{align*}
The last term is finite as long as $\underset{j\in B}{\min}\{c^{(j)}\}$ is bounded away from zero. Hence, \cref{asm.bounded.theta} holds.

Finally, with the analysis on $\lump$, we could show the potential function $\varphi$ is SLC with $T_{\varphi}$, and is self-concordant with $C_\varphi$ under \Cref{asm.bounded.theta} where
\begin{equation}\label{eq.ctphi.def.main}
    T_{\varphi} = \sum_{i\in\cI} \tfrac{w_i}{d_i} \max\left\{\tfrac{6d_i}{(1-r_i)^2}, 2d_i\right\}, \quad C_\varphi = \sum_{i\in\cI} \tfrac{\kappa_{u_i}^3}{\sqrt w_i} d_i^{\invhalf} \max\left\{2, (6r_i^2-6r_i+2)\right\}
\end{equation}
These results are derived immediately from the fact that $\varphi$ is composed of $\{f_i\}_{i\in\cI}$ from the players.
For clarity, we leave the detailed verification to the appendix (cf. \Cref{coro.phi inequality}, \Cref{coro.when.selfconcordant}).

\section{Fast Approximations for Newton-Type Equations}\label{sec.fast.approx}
Before presenting the decentralized IPMs to compute the $\epsilon$-approximate equilibrium price vector \eqref{eq.fosp}, we first introduce the methods to solve the Newton-type equations approximately \eqref{eq.inexact.newton}, which will be used to obtain inexact steps in the IPMs.
First, observe that each of $\nabla^2 f_i(\p)$ is DR1. The affine-scaled Hessian operator reads,
\begin{eqalign}{eq.hessian.exact}
    \bH(\p) & = \bP\nabla^2\varphi(\p)\bP = \sum_{i\in\cI} \frac{w_i}{d_i} \bP \nabla^2 f_i(\p) \bP \stackrel{\eqref{eq.x.jac}}{=} \bP\left(-\sum_{i\in\cI} \nabla \x_i(\p)\right)\bP \\
    & \stackrel{\eqref{eq.f.hess}}{=} \sum_{i\in\cI} \tfrac{w_i}{1-r_i} \bGamma_i - \sum_{i\in\cI}{\tfrac{w_ir_i}{1-r_i}}\bgamma_i\bgamma_i^{\top},
\end{eqalign}
which is a sum of diagonal and rank-one matrices. The following two schemes are constructed from the bidding vectors $\bgamma_1,...,\bgamma_m$, which again consume no extra information from the players.
\subsection{An invertible approximation of the Hessian operator}\label{sec.fast.approx.dr1}

We introduce a simple deterministic DR1 approximation based on the first-order averages. Take $\Omega = \sum_{i\in\cI} \tfrac{w_ir_i}{1-r_i}$, $\omega_i = \Omega^{\inv}{\tfrac{w_ir_i}{1-r_i}}$, then
$\Xi = \sum_{i\in\cI}\omega_i\bgamma_i\bgamma_i^{\top}$ can be seen as the second-order moment of the bidding vectors. It is straightforward to approximate it via weighted first-order averages:
\begin{align*}
    \tilde\Xi = \bxi\bxi^{\top}, \quad \bxi = \sum_{i\in\cI}\omega_i\bgamma_i.
\end{align*}
Hence, we compute a \emph{DR1 approximation} of the Hessian operator:
\begin{eqalign}{eq.hessian.approx}
    \widetilde{\bH}(\p) &= \sum_{i\in\cI} \tfrac{w_i}{1-r_i} \bGamma_i - \Omega \cdot \tilde \Xi,
\end{eqalign}
Note that its inverse operator, $(\widetilde{\bH}(\p))^\inv$, is \emph{explicitly} computable via the Sherman-Morrison-Woodbury formula --- also DR1. An inexact Newton step computed from it requires $\cO(n)$ arithmetic operations.

Assume that in large markets, the best responses are independent and identically distributed.
Since $\bgamma_i \in \Delta_n$, it is expected that the variance, $\Sigma = \va(\bgamma_i)$, diminishes as $\lambda_{\max}(\Sigma) = \cO(\frac{1}{n^g})$ for some $g > 0$. In this case, we can expect the DR1 approximation to be accurate.
\begin{theorem}\label{thm.approx.hessian.dr1}
    Suppose that the endowments $w_i = \tfrac{1}{m}, \forall i\in\cI$, and the best responses are independent and identically distributed with diminishing variance $\Sigma = \va(\bgamma_i)$: there exists $g > 0$ such that $\lambda_{\max}(\Sigma) = \cO(\frac{1}{n^g})$. Then if $n = |\cJ| > \left(\frac{2\Omega}{\epsilon_H}\right)^{\frac{1}{g}}$, and if
    $
        m \ge \frac{8\sqrt{2}(\Omega^2\lambda_{\max}(\Sigma) + 0.5\Omega\epsilon_H)}{\epsilon_H^2}\log\left(\frac{2n}{\delta}\right),
    $
    then
    $0 \preceq \widetilde{\bH}(\p) - \bH(\p) \preceq \epsilon_H\bI$
    holds with probability at least $1-\delta$.
\end{theorem}
An example is when $\bgamma_i$ are i.i.d. from a uniform distribution on the unit simplex, in which case $\lambda_{\max}(\Sigma) = \frac{1}{n(n+1)}$ \citep{kotzContinuousMultivariateDistributions2000},
then if $\|\cJ\| = n \ge \sqrt{\tfrac{2}{\epsilon_H}}$,
it holds that $\|\widetilde{\bH}(\p) - \bH(\p)\| \le \epsilon_H$ with probability at least $1-\delta$.
\subsection{A diagonal preconditioner for Krylov methods}\label{sec.precond}
We next introduce a simple diagonal preconditioner, so that \eqref{eq.inexact.newton} can be solved efficiently by Krylov methods. The preconditioner simply uses the row summation of $\bH(\p)$:
\begin{eqalign}{eq.def.precond}
    \Kc = \diag(\kc),~\text{with} \quad\kc &= \bH(\p)\1 \stackrel{\eqref{eq.hessian.exact}}{=} \sum_{i\in\cI} w_i \bgamma_i.
\end{eqalign}
The centered preconditioned matrix of $\bH(\p)$ is
\begin{equation}\label{eq.precond.system}
    \bH_c := \Kc^{\invhalf}\bH(\p)\Kc^{\invhalf}
\end{equation}
The following theorem shows that the condition number of the scaled system is minimized.
\begin{theorem}[The scaled system is optimal]\label{thm.cond.precond}
    Consider the centered preconditioned system \eqref{eq.precond.system} where $\kc,\Kc$ are defined as in \eqref{eq.def.precond}. Suppose that $r_i \equiv r < 1$ for all $i\in\cI$.
    Then the following holds:
    \begin{enumerate}
        \item[$(a)$] The condition number of $\bH_\mathrm{c}$ is bounded by
            \begin{equation*}
                \kappa_2(\bH_{\mathrm c}) = \tfrac{\lambda_{\max}(\bH_{\mathrm c})}{\lambda_{\min}(\bH_{\mathrm c})} \le
                \begin{cases}
                    \frac{1}{1-r}, & 0 \le r < 1, \\
                    1 - r,         & r < 0.
                \end{cases}
            \end{equation*}
        \item[$(b)$] Furthermore, this condition number cannot be improved by any diagonal preconditioner. That is, suppose that $\cK \subseteq \real^{n\times n}_+$ is the space of diagonal matrices, then $\kappa_2(\bK^\invhalf \bH(\p) \bK^\invhalf) \ge \kappa_2(\bH_\mathrm{c})$, $\forall \bK \in \cK.$
\end{enumerate}
\end{theorem}
Hence, if $r > 0$, then \eqref{eq.precond.system} can be solved in $\cO(n^2\sqrt{(1-r)^\inv}\log(\tfrac{1}{\epsilon}))$ time.
For CES utilities, if either $r\to -\infty$ (Leontief) or $r\to 1$ (Linear), then the condition number tends to be unbounded -- actually, $\varphi(\p)$ is not differentiable. In fact, the convergence rates of first-order methods degrade to $\cO(\tfrac{1}{\epsilon})$ from $\cO(\log(\tfrac{1}{\epsilon}))$ for $r \in (-\infty, 1)$ \citep{cheungDynamicsDistributedUpdating2018}.
\begin{figure}[h]
    \centering
    \subcaptionbox{
        {\small Comparison of $\widetilde{\bH}(\p)$ and $\bH(\p)$ using DR1 approximation under iid uniform distribution over $\Delta_n$.}
    }[0.45\textwidth]{
        \includegraphics[width=0.33\textwidth]{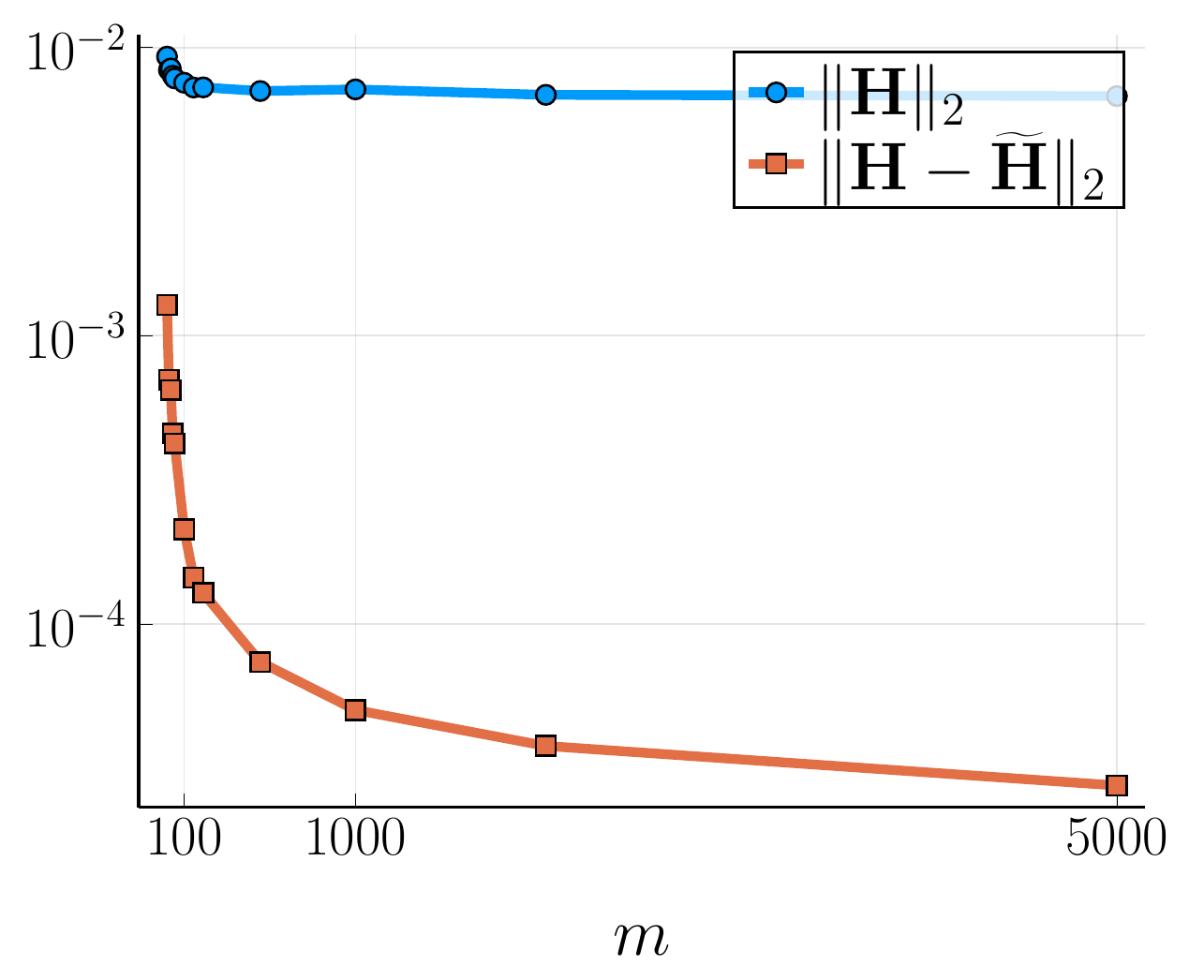}
    }
    ~
    \subcaptionbox{
        {\small Solving the linear system with/without diagonal preconditioner $\bK_c$. }
    }[0.45\textwidth]{
        \includegraphics[width=0.35\textwidth]{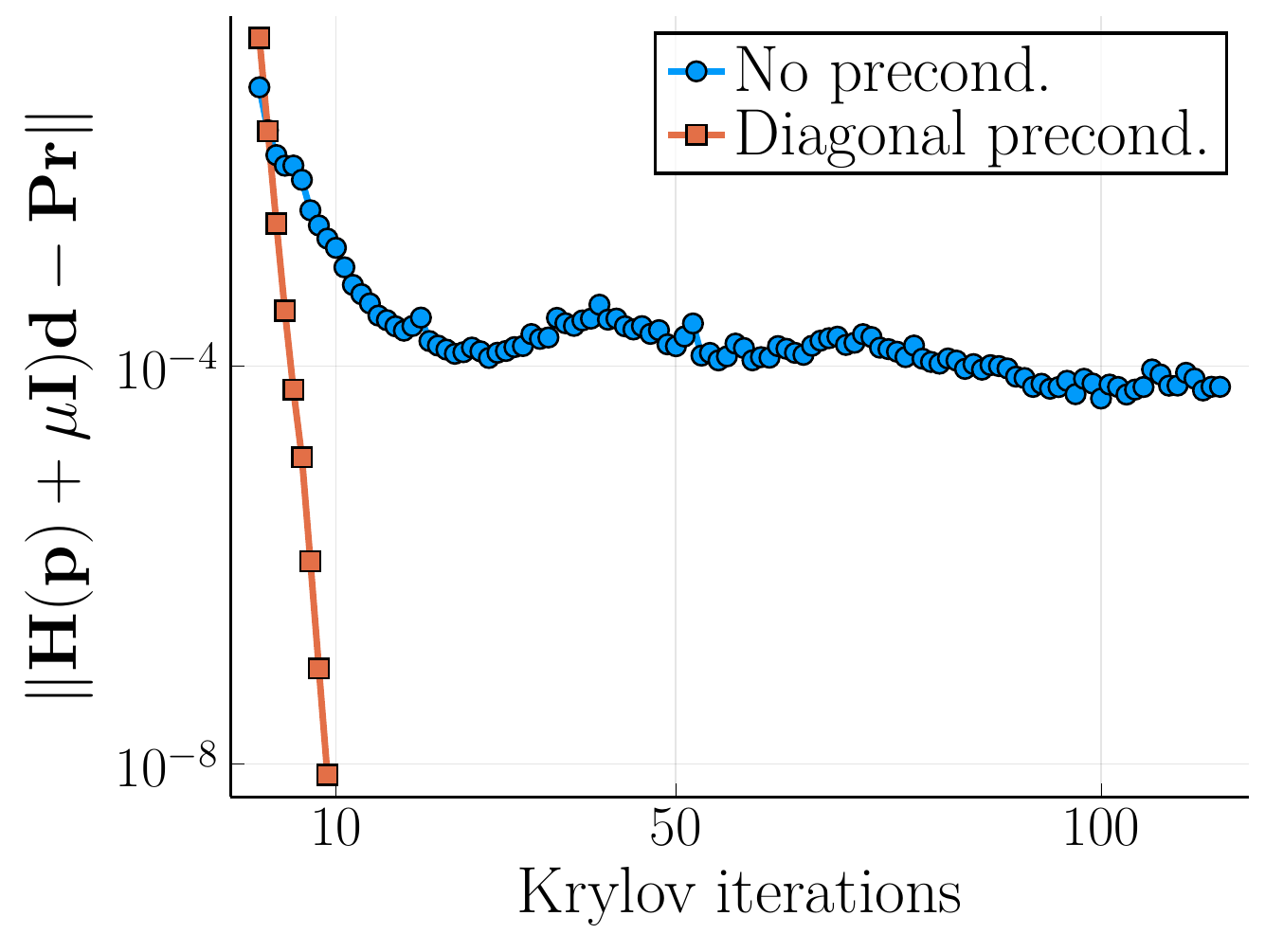}
    }
    \caption{Quality of approximations when solving the Newton equations.}
    \label{fig.approx}
\end{figure}

We illustrate these approximations in \Cref{fig.approx} for the CES economy (cf. \Cref{ex.ces}) with $n=2,000$. In the left figure, we plot the approximation error of DR1 while enlarging the number of players $m$ if the coefficients are generated according to the uniform distribution over $\Delta_n$ (and so the variance diminishes). The operator norm of the scaled Hessian operator is also plotted for comparison. These results show $m\ge 100$ is sufficient for inexact solutions. In the right figure, we consider solving the regularized linear system with matrix $\bH + \mu\bI$ with and without $\Kc$. The $y$-axis shows the residual of the linear equation. By comparing the required Krylov iterations, the results show that the preconditioned system is much easier to solve.

\section{Convergence Analysis of Second-Order T\^atonnement Processes}\label{sec.conv}
In this section, we provide convergence analysis of the second-order t\^atonnement processes (IPMs) using inexact $\widetilde{\bH}(\p)$ in the Newton steps to compute the stationary point of $\varphi$.  Note that $\varphi$ is SLC with parameters $T_{\varphi}$, and is self-concordant with parameter $C_\varphi$ under \Cref{asm.bounded.theta} (cf. \eqref{eq.ctphi.def.main}, \Cref{coro.phi inequality,coro.when.selfconcordant}).
To proceed, we present the following assumption for the convergence analysis.
\begin{assumption}\label{asm.approx.hessian}
    For any $\epsilon_H > 0$, the affine-scaled approximate Hessian operator $\widetilde{\bH}(\p)$ satisfies
    \begin{equation}\label{eq.approx.hessian.cond}
        \widetilde{\bH}(\p) \succeq 0, \quad \|\widetilde{\bH}(\p) - \bH(\p)\| \le \epsilon_H.
    \end{equation}
\end{assumption}
We will assume this general form so that the convergence analysis is not affected by any specific choice of $\widetilde{\bH}(\p)$.
In general, this condition can be satisfied by subsampling or sketching methods \citep{woodruffSketchingToolNumerical2014,pilanciNewtonSketchLineartime2017,drineasRandNLARandomizedNumerical2016}. For example, if using a subsampling method that randomly picks a subset of players, then the size of needed players is about $\Omega(\tfrac{1}{\epsilon_H^2} \log( \frac{n}{\delta}))$ to achieve \eqref{eq.approx.hessian.cond} with probability at least $1-\delta$ \citep{xuNewtontypeMethodsNonconvex2020}.
Alternatively, if we use Krylov methods to solve the linear systems inexactly, generally only conditions like \eqref{eq.def.epsH.direction} can be guaranteed.
For conciseness, we only focus on the analysis, assuming \eqref{eq.approx.hessian.cond}, but the analysis based on Krylov solvers should be migrated easily to the present framework; see, for example, \cite{xuNewtontypeMethodsNonconvex2020,zhangHomogeneousSecondorderDescent2025}.

\subsection{Convergence analysis of \Cref{alg.logbarrier}}
In the \Cref{alg.logbarrier}: \logbar{}, we use the barrier model \eqref{eq.log barrier varphi}, where the potential $\varphi(\cdot)$ is SLC with parameter $T_\varphi$ (\Cref{coro.phi inequality}). For the convenience of the reader, recall the neighborhood of the central path \eqref{eq.central path}:
\begin{equation}\label{eq.logbar.ac.repeat}
    \cC(\mu,Q) = \left\{\p > 0: \tfrac{\|\bP \nabla \varphi(\p) - \mu \1\|}{\mu} \le Q \right\}.
\end{equation}
For a starting $\mu_0$, we need an initialization point $(\mu_0, \p_0)$ such that $\p_0 \in \cC(\mu_0, Q)$. Generally, almost all path-following IPMs require a separate process to find such $\p_0$ (in various kinds). This process is highly nontrivial and has been a major challenge for traditional IPMs \citep{nesterovInteriorPointPolynomialAlgorithms1994,yeInteriorPointAlgorithms1997,enagyNewInteriorpointAlgorithm2024}. For the market equilibrium, we could again benefit from an economic question: what happens if the prices are too large? The answer leads to an initial AC.
\begin{lemma}[Approximate AC]\label{lem.logbar.init}
    Suppose that $\p_0 = \mu_0\1,$ for some $\mu_0 > 0$. Then if $\mu_0 \ge \sqrt{\tfrac{\sum_{i\in\cI} w_i}{Q}}$, we have $\p_0 \in \cC(\mu_0, Q)$ for any $Q \in (0, \tfrac{1}{2})$.
\end{lemma}
With the specific choice of \Cref{lem.logbar.init}, we have
$$\bP_0 \nabla \varphi(\p_0) - \mu_0 \1 = \mu_0\bI \left(\1-\smallsum_{i\in\cI}\x_i(\mu_0\1)\right) - \mu_0 \1 = -\smallsum_{i\in\cI}\x_i(\mu_0\1).$$
The demand would diminish if $\mu_0$ is sufficiently large: $\lim_{\mu_0\to\infty} \sum_{i\in\cI}\x_i(\mu_0\1) = 0$. Hence, \eqref{eq.logbar.ac.repeat} must hold for some $\mu_0$ that is large enough. Specifically, since for each $i$, $\inner{\mu_0\1}{\x_i(\mu_0\1)} = w_i$, we have,
$\left\|\sum_{i\in\cI} \x_i(\mu_0\1)\right\| \le \sum_{i\in\cI} \|\x_i(\mu_0\1)\|_1 = \frac{\sum_{i\in\cI} w_i}{\mu_0}.$
An initial AC \eqref{eq.logbar.ac.repeat} would require $\tfrac{\sum_{i\in\cI} w_i}{\mu_0} \le \mu_0 Q$, that is, $\mu_0 \ge \sqrt{\tfrac{\sum_{i\in\cI} w_i}{Q}}$.

At the following main phase, each step is computed by solving the inexact Newton system \eqref{eq.inexact.logbar}. Each iterate is guaranteed to stay in the neighborhood of the central path.  Finally, when the parameters are properly set, we would conclude the following iteration complexity of \logbar{} to compute an $\epsilon$-approximate equilibrium price vector. The overall dependency of iteration complexity on the choice of $\mu_0$ will be on a logarithmic scale.
\begin{theorem}[Convergence rate of \Cref{alg.logbarrier}]\label{thm.complexity.logbarrier}
    Suppose \Cref{ass.logump.unique.bestresp,asm.approx.hessian} hold. If in \Cref{alg.logbarrier} we choose $\sigma = \frac{Q + \sqrt{n}}{2Q + \sqrt{n}} < 1$ and $Q \le \frac{\epsilon}{14\epsilon + 4T_{\varphi}(\sqrt{n} + 1)},$
    then the number of Newton iterations needed to obtain an $\epsilon$-approximate equilibrium price vector \eqref{eq.fosp} satisfies $K \le \cO(\sqrt{n}\log(\tfrac{\mu_0}{\epsilon})).$
    If DR1 approximation is used (\cref{thm.approx.hessian.dr1} holds), then this estimate holds with probability at least $1- K\delta, \delta \in (0,1)$.
\end{theorem}
For large markets, if DR1 approximation can be applied (\cref{thm.approx.hessian.dr1} holds), each Newton equation is solved in $\cO(n)$ arithmetic operations. Hence,
the iteration complexity implies that the total number of arithmetic operations is about $\cO(n^{1.5}\log(\tfrac{1}{\epsilon}))$ with probability at least $1- K\delta, \delta \in (0,1)$.
In practice (see \cref{sec.numerical}), we find $K$ is about $10^1$ to $10^2$ for most tested instances, even using the simple DR1 approximation.

\subsection{Convergence analysis of \Cref{alg.pathfollowing}}
We now present the convergence analysis of \Cref{alg.pathfollowing}: \pathf{}. Our analysis extends the analysis in \cite{dvurechenskyImprovedGlobalPerformance2025} to account for approximation errors. Under active-set boundedness (\Cref{asm.bounded.theta}), the potential $\varphi(\cdot)$ is $C_\varphi$-self-concordant (see \Cref{coro.when.selfconcordant}). Then the homotopy problem \eqref{eq.pathf.varphi}
$$
    \varphi_t(\p) = \varphi(\p) - t\inner{\nabla\varphi(\p_0)}{\p}, \quad t \in [0,1]
$$
is also $C_\varphi$-concordance for any $t$. To proceed, let us verify the relative error of the DR1 approximation in the presence of the active-set boundedness.
\begin{lemma}[Relative error of the DR1 approximation]\label{lem.relative.error.hess}
    Suppose that \Cref{asm.bounded.theta} holds and the affine-scaled approximate Hessian $\widetilde{\bH}(\p)$ satisfies \eqref{eq.approx.hessian.cond}. If the DR1 approximation is used, defining
    \begin{equation}
        \delta:=\frac{\epsilon_H}{\min_{i\in\cI} \tfrac{d_i}{\kappa_{u_i}}},
    \end{equation}
    then $\widetilde{\bH}$ is also a \emph{$\delta$-spectral approximation} of $\bH$,
    \begin{equation}
        \bH \preceq \widetilde{\bH} \preceq (1+\delta)\bH.
    \end{equation}
\end{lemma}
We could see in the proof that $\p$ is close to the central path and thus bounded, so this assumption is quite reasonable.
By DR1, the absolute error $\epsilon_H$ implies stronger relative deterioration of the second-order derivatives. For Krylov methods, we assume the same since we can efficiently solve the linear system to any desired precision.

In the convergence analysis, we first evaluate the cost of inexactness in terms of the Newton decrement (cf. \Cref{lem.descent.pathf}). This basic fact has several implications. With an inexactness measure $\delta$, a region for quadratic convergence still exists. It happens when the inexact Newton decrement $\widetilde{\lambda}(\p)$ is small enough (cf. \eqref{eq.pathf.quadregion}). Furthermore, once the iterates are centered, it never leaves (cf. \Cref{lem.rules.betagamma}). On the global scale, we could show that the homotopy parameter $t_k$ is decreasing superlinearly. Once $t_k = 0$, we enter the region of quadratic convergence (\Cref{cor.pathf.desc.small}), and it further indicates the global superlinear convergence of the method itself.
\begin{theorem}[Convergence rate of \Cref{alg.pathfollowing}]\label{thm.convergence.rate.pathf}
    Suppose \Cref{ass.logump.unique.bestresp,asm.bounded.theta,asm.approx.hessian} hold. If in \Cref{alg.pathfollowing}, $(\delta, \sigma, \gamma)$ is chosen by \eqref{eq.pathf.choice.bgchoice}, then for all $k \ge 1$,
    \begin{equation}
        t_k \le \left(1 - \frac{\gamma(\gamma - 2\beta)(k+1)}{2C_{\varphi}^2(\varphi(\p_0) - \varphi(\p^{*}))}\right)^{k+1}.
    \end{equation}
    That is, \Cref{alg.pathfollowing} has global superlinear convergence.
\end{theorem}
When the Hessian approximations are good enough, we could show that \Cref{alg.pathfollowing} has a non-asymptotic superlinear rate of convergence.
Note that the CES utilities belong to $\cA^{k,r}$. To the best of our knowledge, no such results exist in the literature for CES utilities even with exact high-order derivatives.
 \section{Further Applications} \label{sec.extend}
We present two applications of our analyses on $\lump$ to accomodate the linear utility funtions and linear constraints in the fisher market model.

\subsection{Linear exchange market problem}\label{sec.extend.linear}
For $\cA^{k,r}$ family, the uniqueness of the best response (\Cref{ass.logump.unique.bestresp}) holds if $-\infty < r < 1$.
Now we consider the linear utility functions that may be seen as a limit case where $k=r=1$ and for each $j$, $\theta_j(x^{(j)}) = c^{(j)}x^{(j)}$, $c^{(j)} \in \real_+$. Specifically, the market equilibrium is given by
\begin{equation}\label{eq.equilibrium.linear}
    \z(\p) :=\sum_{i\in\cI} \x_i - \1 = 0, \quad \x_i \in \arg \max ~ \inner{\c_i}{\x_i}, ~\st ~\inner{\p}{\x_i} \le w_i.
\end{equation}
As before, $\x_i$ also solves the LUMP by dropping index $i$:
\begin{eqalign}{eq.linear.lump.original}
    f(\p) = \underset{\x \in \real_+^n}{\max} ~& -v(\x) := \log(\inner{\c}{\x})\\
    \st ~& \inner{\p}{\x} \le w.
\end{eqalign}
It is well-known that the best-response mapping is set-valued \citep{mas-colellMicroeconomicTheory1995}, which uses up the budget by any combination of goods according to the maximum utility-price ratio: $\arg\max_j\{\frac{c_j}{p_j}\}$. Let us consider the barrier-regularized $\lump$ so that we recover the differentiability of the dual function. Given parameter $\sigma > 0$, we introduce
\begin{eqalign}{eq.linear.lump}
    f^\sigma(\p) = \underset{\x \in \real_+^n}{\max} ~& -v^\sigma(\x) := \log(\inner{\c}{\x}) + \sigma \inner{\log(\x)}{\1} \\
    \st ~& \inner{\p}{\x} \le w.
\end{eqalign}
With a slight abuse of notation, we write $\x(\p)$ as the solution of \eqref{eq.linear.lump}.
The concept of approximate optimal play is widely used for linear utilities (e.g., see \cite{azevedoStrategyproofnessLarge2019}).
Similarly, we know that $v^\sigma$ is logarithmically homogeneous:
$$-v^\sigma(\alpha \x) = \log(\alpha) (1 + \sigma n) + \log(\inner{\c}{\x}) + \sigma \inner{\log(\x)}{\1} = \log(\alpha) (1 + \sigma n) + v^\sigma(\x), \quad \forall \alpha > 0.$$
We could quickly replicate the calculus: $\nabla v^\sigma(\x) = -\frac{\c}{\inner{\c}{\x}} - \sigma \bX^\inv\1$, $\nabla^2 v^\sigma(\x) = \frac{\c\c^\top}{\inner{\c}{\x}^2} + \sigma \bX^{-2}$, and
\begin{align*}
    \nabla^3v^\sigma(\x) & = -2\sigma\bX^{-3} - \tfrac{2}{\inner{\c}{\x}^3}\c\otimes\c\otimes\c.
\end{align*}
Then, similar to \eqref{eq.bounded.directdev} we obtain the scaled boundedness of directional derivatives,
\begin{eqalign}{eq.linear.lump.dir3}
    \left|\nabla^3v^\sigma[\h,\h,\h]\right| = \left|-2\sigma \inner{\1}{[\bX^\inv\h]^3} - 2 \inner{\bgamma}{\bX^\inv\h}^3\right| \le (2\sigma + 2) \|\bX^\inv\h\|^3
\end{eqalign}
Using the first-order condition, the following result is immediate. Notice that the dual variable $\lambda$ to the budget constraint is, as before, a constant.
\begin{lemma}\label{lem.linear.lump.foc}
    If $\sigma > 0$ is a constant, suppose that $\x$ is the best response, and $u = \inner{\c}{\x}$, then the following holds
    $$
        \x = \sigma\left(\frac{1+\sigma n}{w}\p - \frac{1}{u}\c\right)^\inv, \quad u = \inner{\c}{\x};
    $$
    $u$ is the root of the following equation,
    $$
        \psi(u) = \sum_{j\in\cJ}\sigma c^{(j)}\left(\frac{1+\sigma n}{w}p^{(j)} - \frac{1}{u}c^{(j)}\right)^\inv - u = 0, \quad u \in \real_+;
    $$
    and the Hessian of the dual function is
    \begin{equation}
        \nabla^2 f(\p) = \frac{1}{w}\nabla \x = \frac{1}{w} \left(\frac{\lambda}{\sigma} \bX^{2}
        - \frac{\lambda}{\sigma(\sigma + \|\bgamma\|^2)}\bX\bgamma\bgamma^\top\bX\right),
    \end{equation}
    where $\lambda = \tfrac{1+\sigma n}{w}$, and
    $\bgamma = (1+\sigma n)\frac{\bX\p}{w} - \sigma \1= \frac{\bX\c}{\inner{\c}{\x}}$ is the shifted bidding vector.
\end{lemma}
Note that for classical linear UMP, the best-response mapping could be enumerated by a priority rule; hence, it admits an analytic expression.
The only additional step for $\lump$ with $\sigma$ regularization is solving a 1-dimensional root-finding problem. But similarly, the information required here can be collected from the best-response mappings.

\begin{theorem}[SLC of linear $\lump$]\label{thm.linear.lump.slc}
    The negative logarithmic utility, the dual function, and the best-response mapping of linear $\lump$ \eqref{eq.linear.lump} is SLC.
\end{theorem}
We leave the coefficients of the SLC (in $\cO(\sigma^{-3})$) in the appendix.
For linear utility, we could carry out the same machinery as before. To compute the competitive equilibrium, the potential function in this case would be
\begin{eqalign}{eq.linear.ipmeg.phi}
    \varphi^\sigma(\p) = \inner{\p}{\1} + \sum_{i\in\cI} w_i f_i^\sigma(\p).
\end{eqalign}
The only difference is that the $f^\sigma_i$ is the dual function of a barrier-regularized $\lump$. By noticing $\lambda_i = \tfrac{1+\sigma n}{w_i}, \s_i = \sigma\bX_i^\inv\1$, $\forall i\in\cI$, we write the optimal condition as
\begin{subalign}{eq.linear.ipmeg.foc}
    \label{eq.linear.ipmeg.foc.1}    \1 - \smallsum_{i\in \cI} (1+\sigma n) \x_i &= 0, \\
    \label{eq.linear.ipmeg.foc.2}     \inner{\c_i}{\x_i}^\inv\c_i + \s_i - \tfrac{1+\sigma n}{w_i} \p = 0, ~ \bX_i \s_i &= \sigma \1, ~\forall i\in\cI.
\end{subalign}
Clearly, the solution set of \eqref{eq.linear.ipmeg.foc} is a single-valued mapping and is a reasonable approximation: if $\sigma$ converges to $0$, then it converges to the solution set of the linear Fisher market.
The remaining question is how good the approximation is. We summarize in the following remark.
\begin{remark}
    When applying the IPMs, \eqref{eq.linear.ipmeg.foc.2} holds along the way since we apply the best-response mapping of $\lump$. For example, if we use \logbar{} and it is terminated at some tolerance $\epsilon > 0$, we have, $\left\|\1 - \sum_{i\in \cI} (1+\sigma n) \x_i\right\|_\infty \le \epsilon$. So we have
    \begin{align*}
        (1+\sigma n)\left\|\smallsum_{i\in \cI} \x_i - \1\right\|_\infty \le \left\| \smallsum_{i\in \cI} (1+\sigma n) \x_i - \1 - \sigma n\1\right\|_\infty \le \epsilon + \sigma n.
    \end{align*}
    Hence, $\left\|\sum_{i\in \cI} \x_i - \1\right\|_\infty \le \frac{\epsilon + \sigma n}{1+\sigma n}$. If $\sigma = \Theta(\tfrac{\epsilon}{n})$, we have an $\epsilon$-approximated equilibrium as defined \eqref{eq.fosp}.
    By \Cref{thm.complexity.logbarrier}, because we could find a starting AC from \Cref{lem.logbar.init}, a $\epsilon$-approximated equilibrium could be found in $\cO(\sqrt{n}\log(\tfrac{1}{\epsilon}))$ time.
\end{remark}
For linear utility, each iteration requires a linear system in $\mathcal {O} (n^3)$ arithmetic operations if matrix factorizations are used. This could be improved by the techniques in \Cref{sec.fast.approx} and \pathf{}.
When the data of this market is rational, we could find the exact equilibrium using the interior-point rounding procedures like \cite{mehrotraFindingInteriorPoint1993,andersenCombiningInteriorpointPivoting1996}.

\subsection{Allocation in the subspace}\label{sec.extend.network}
We now consider the case where $\cX$ is not full-dimensional. One such possibility is to place linear constraints $\bA\x = \d$ in the allocation space.
As before, the goal is to find some price-allocation pair $(\p, \x_1, ..., \x_m)$ such that the market clears:
\begin{equation}\label{eq.equilibrium.network}
    \z(\p) :=\sum_{i\in\cI} \x_i - \1 = 0, \quad \x_i \in \arg \max ~ u_i(\x_i), ~\st ~\bA_i\x_i = \d_i, ~\inner{\p}{\x_i} \le w_i.
\end{equation}
In general, if for all $i\in\cI$, $\d_i \neq 0$, the equilibrium price may not exist; even it does, the set of equilibria is non-convex and PPAD-hard to acquire \citep{jalotaFisherMarketsLinear2023}.
Luckily, the homogeneous case where $\d_1 = ... = \d_m = 0$ is tractable -- the allocation is in the affine space. We are motivated by the following examples.
\vspace{-1em}
\begin{example}[Characteristics \& aggregates]
    Consider where the utility function is written as $u_i(\bB_i\x_i)$, for some $\bB_i\ge 0$.
    This utility model has been used to describe characteristics of the goods \citep{lancasterNewApproachConsumer1966}, household production \citep{beckerTheoryAllocationTime1965}, two-stage budgeting and nested demand systems \citep{armingtonTheoryDemandProducts1969} if the group aggregation is linear, and demand systems with connected substitutes \citep{berryConnectedSubstitutesInvertibility2013}.
    We can introduce a new variable $\y_i=\bB_i\x_i$ and rewrite the utility as $u_i(\z_i), \z_i = [\y_i; \x_i]$. Hence, the allocation is subject to the linear constraints:
    $
        \bA_i \z_i = 0, \bA_i = [-\bI, \bB_i].
    $
\end{example}
The second example appears in network applications \citep{jainEisenbergGaleMarkets2010}, where each edge or node is regarded as a good.
\begin{example}[Allocation in the networks]
    Consider when any player $i$ allocates a $s_i- t_i$ flow in a network $G = (V,E): \x_i=[x_{i,0}; x_{i,e}, e\in E], s_i \in V, t_i \in V$.  Such a feasible set of allocation is encoded by a homogeneous feasibility system for $x_{i,0}, x_{i,e} \in \real_+, \forall e\in E$ as follows,
    \begin{equation}\label{eq.lump.flow}
        \begin{array}{rl}
            x_{i,0}+   \sum_{e \in \delta^{+}(s)} x_{i,e}-\sum_{e \in \delta^{-}(s)} x_{i,e}=0, &                     \\
            -x_{i,0}+  \sum_{e \in \delta^{+}(t)} x_{i,e}-\sum_{e \in \delta^{-}(t)} x_{i,e}=0, &                     \\
            \sum_{e \in \delta^{-}(v)} x_{i,e}-\sum_{e \in \delta^{+}(v)} x_{i,e}=0,            & v \in V - \{s, t\}.
        \end{array}
    \end{equation}
    where $\delta_-, \delta_+$ are outflow and inflow edges, respectively. Clearly, all the constraints are homogeneous equalities, and so we can write $\bA_i\x_i = 0, \x_i \in \real_+^n$ for short.
\end{example}
Dropping the subscript $i$, $\lump$ is given by,
\begin{equation}\label{eq.lump.affine}
    \begin{array}{rl}
        f(\p) := \underset{\x\in\cX}{\max} & -v(\x) = \log(u(\x))               \\
        \st                                & \inner{\p}{\x} \le w, ~\bA \x = 0.
    \end{array}
\end{equation}
The following lemma shows that $\nabla f(\p) \propto \x(\p)$ and so the convex optimization of the potential function $\varphi(\p)$ is still permitted; the Hessian $\nabla^2 f(\p)$ is still computable. Unfortunately, the DR1 property is lost and it would require a bit more computations as does an individual response itself.
\begin{theorem}\label{thm.lump.affine.hess}
    For any $\p \in \cP$, let $\x$ be the best response that solves \eqref{eq.lump.affine}. Then, the dual variable to the budget constraint is given by $\lambda = \tfrac{d}{w}$.
    Take $\bW = \nabla^2 v(\x)$, then
    the gradient and Hessian of $f(\p)$ are given by
    \begin{equation}
        \nabla f(\p) = -\lambda\x(\p), ~\nabla^2 f(\p) = \frac{d^2}{w^2} (\bW^\inv-\bW^\inv\bA^{\top}(\bA\bW^\inv\bA^{\top})^\inv\bA\bW^\inv).
    \end{equation}
\end{theorem}
Recall that we have shown that $\nabla^2 v(\p)$ is contructed from the best responses, precisely, the bidding vectors ($\bgamma_1, ..., \bgamma_m$), we could say that in this case, a second-order price updating mechanism as before requires no extra information.
For $\cA^{k, r}$ utilities, the scaled Lipschitz continuity is still valid.
\begin{theorem}\label{thm.scaled.lipschitz.network}
    If $u \in \cA^{k, r}$, the dual function of a homogeneously affine-constrained $\lump$ \eqref{eq.lump.affine} is scaled Lipschitz continuous. The coefficient associated with the scaled Lipschitz continuity is no bigger than $T_f$ (cf. \eqref{eq.f.tf.def}) in the full-dimensional case (\cref{thm.taylor.f}).
\end{theorem}

To summarize, we could apply the interior-point methods to solve \eqref{eq.equilibrium.network} by updating the price according to players' best responses; the estimates of iteration complexity remain the same.
 
\section{Numerical Experiments}\label{sec.numerical}
We present a preliminary numerical experiment to show the efficiency of second-order t\^atonnement processes to solve the Fisher model.
All experiments are conducted on a single machine with a 14-core Apple M4 Pro CPU and 48GB LPDDR5 RAM using Julia Programming Language 1.10.
We implement our proposed method, including \logbar{} (\Cref{alg.logbarrier}), \pathf{} (\Cref{alg.pathfollowing}), using the DR$1$ approximation \eqref{eq.hessian.approx}. Besides, we add another method that uses the preconditioned conjugate gradient method in the \Cref{alg.logbarrier}; this method uses the diagonal preconditioner \eqref{eq.precond.system} and is denoted by \logbarpcg{}.

To test our proposed methods, we use the CES utility function (cf. \cref{ex.ces}) with the utility function as follows:
$
    u(\x) = \left[\smallsum_{j\in\cJ} c^{(j)} (x^{(j)})^{\rho}\right]^{\frac{1}{{\rho}}} = \inner{\c}{\x^{\rho}}^{\frac{1}{{\rho}}}, {\rho} \in (-\infty,1), \c \in \real^n_+. $
We compare with a few first-order methods (FOMs) in the CES economy: $(1)$ the Tâtonnement process (\tat{}) as a mirror descent method, based on \cite{birnbaumDistributedAlgorithmsGradient2011,cheungTatonnementGrossSubstitutes2013}; $(2)$ the Proportional Response (\pr{}) implemented according to \cite{cheungDynamicsDistributedUpdating2018}.

\begin{figure}[h]
    \centering
    \subcaptionbox{$\rho=-1.9$}{
        \includegraphics[width=0.35\textwidth]{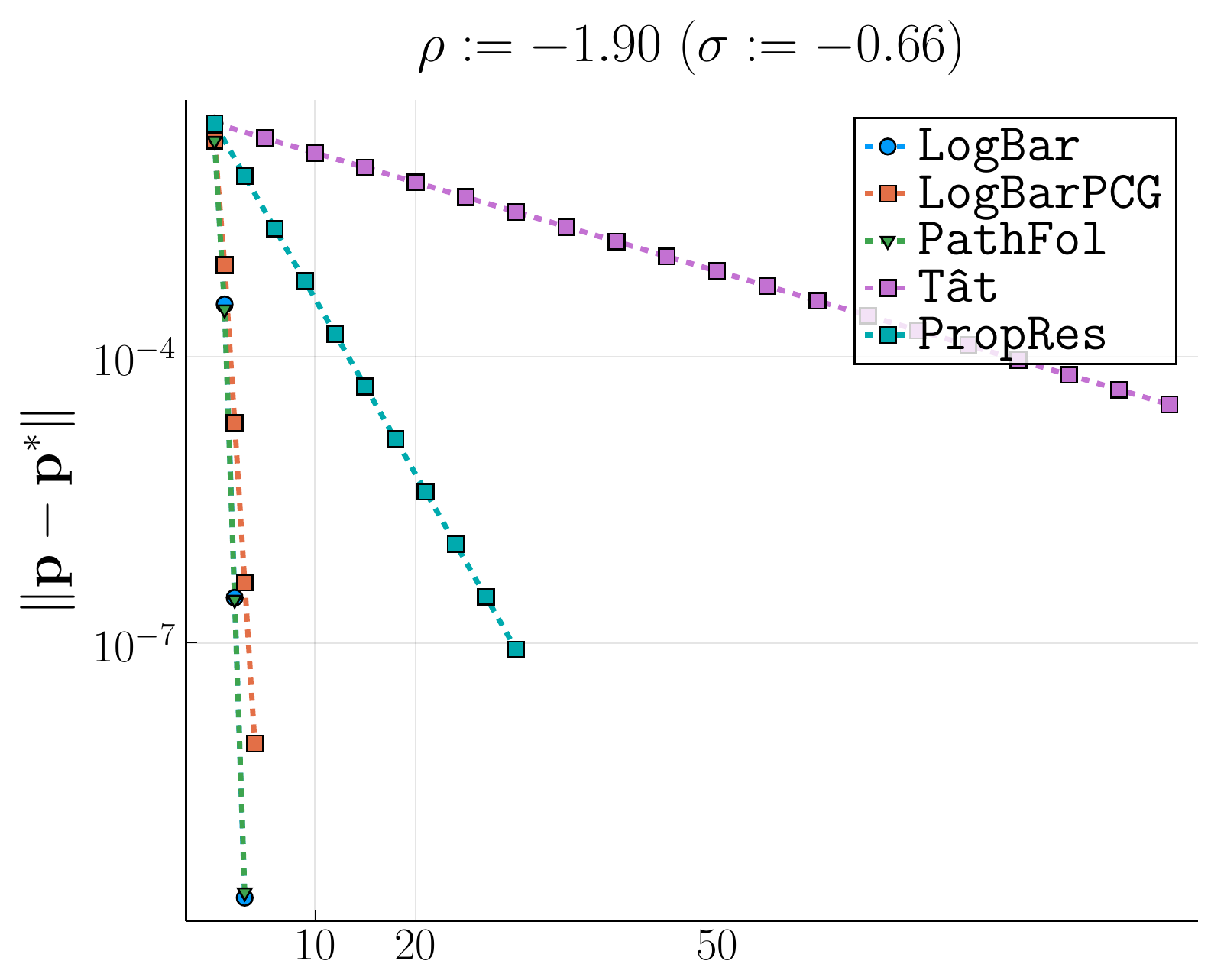}
    }
    \subcaptionbox{$\rho=0.9$}{
        \includegraphics[width=0.35\textwidth]{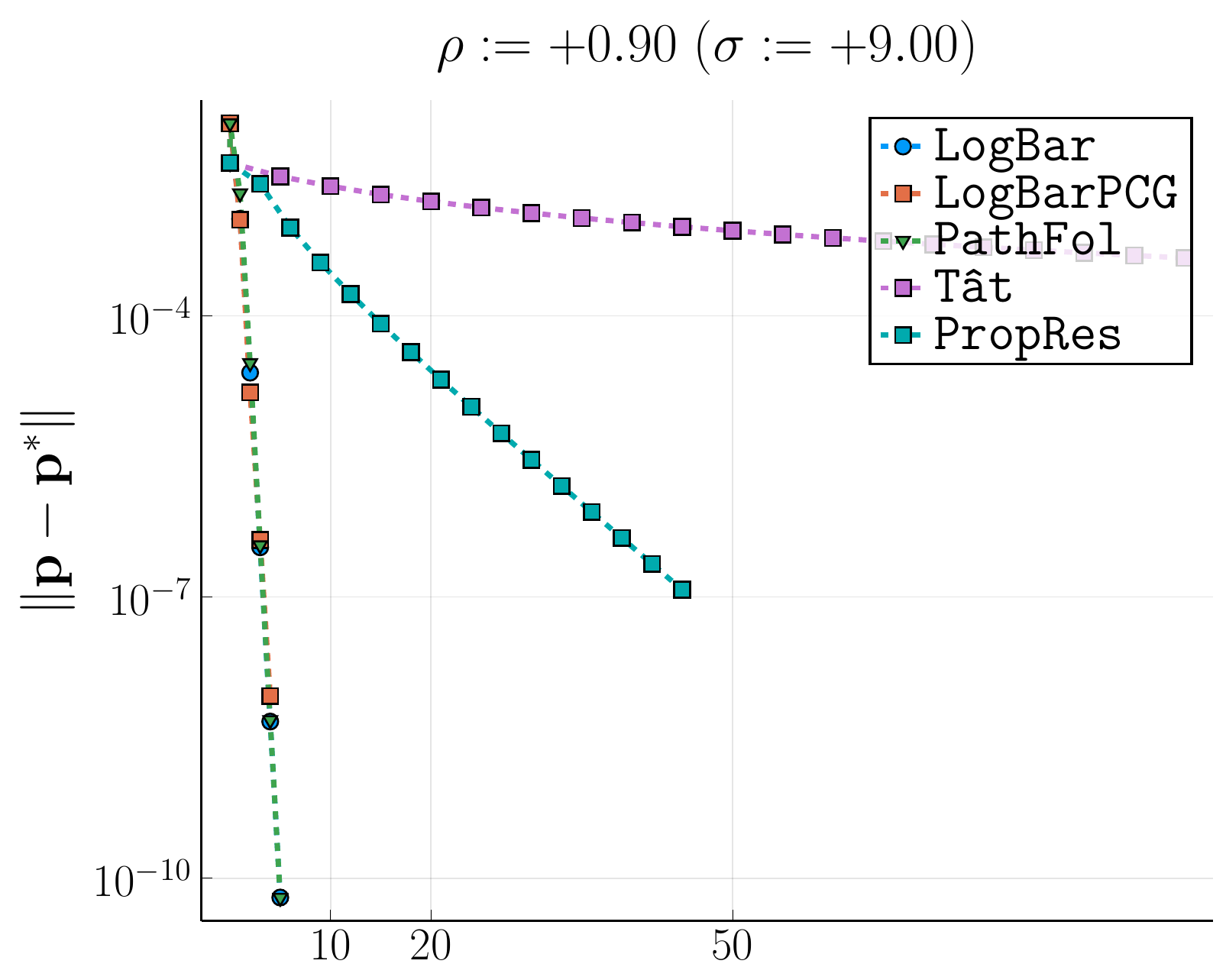}
    }
    \caption{Convergence behavior of Interior-Point Methods and the First-Order Methods in CES economy using the randomly generated dataset. Parameters: $n = 1,000$, $m = 3,000, \tau=0.2$.
        The horizontal axis is the number of iterations.
        The vertical axis is the distance between the current price to the ground truth $\p^*$.
        The parameter $\sigma = \frac{\rho}{1-\rho}$.
    }\label{fig.ces.rho}
\end{figure}

In the randomly generated instances, we generate a random matrix $\c = [\c_1,..., \c_m]$ to represent the coefficients in the utility:
$
    \c = \delta \times \texttt{sprand}(\real_+^{m \times n}, \tau),
$
where $\delta$ is a scaling factor and $\tau$ is the sparsity level; $\texttt{sprand}$ is the function to generate a sparse random matrix in Julia.
The endowments $w^{(1)},...,w^{(m)}$ are generated from the uniform distribution and normalized to unit summation.
In \Cref{fig.ces.rho}, we present the convergence behavior of the proposed methods and FOMs in the substitution regime ($\rho=0.9$) and the complementarity regime ($\rho=-1.9$), respectively.
Note that both our methods and \tat{} are based on best-response mappings, while \pr{} does not guarantee player-optimality in the process. Instead, the market clearance is always guaranteed.
To make a fair comparison, we focus on the distance between the current price to the ground truth $\p^*$, and the stopping criterion is set as: $\|\p - \p^*\| \le 10^{-7}.$

In both cases, \logbar{}, \logbarpcg{}, and \pathf{} behave similarly and converge faster than the other methods.
Now we test larger instances to show the scalability of our proposed methods by enlarging $n, m$ to $10^3-10^5$.  Again, we terminate the methods at an intermediate precision: $\|\p - \p^*\| \le 10^{-5}$ or a time limit of 200 seconds. The results are presented in \Cref{tab.auction.rho.compare}. In all cases, \logbar{} has the best running time and outputs the solution with the highest precision (9-10 digits). \logbarpcg{}  is also quite competitive with the preconditioner; most of the linear equations can be solved in $\le 20$ Krylov iterations. We omit \pathf{} because it has no obvious difference from \logbar{}.
\begin{table}[h]
    \centering
    \caption{Results of finding $\epsilon$-approximate equilibrium price vector for different problem sizes. Parameters: $\tau=0.2$. The time is in seconds. The symbol - marks the case that the method fails to converge, and $^\dagger$ denotes the inaccurate solution.}\label{tab.auction.rho.compare}
    \scriptsize
    \begin{longtable}{lrr|rrrr|rrrr}
        \toprule
                                  &                             &        & \multicolumn{4}{c|}{$t$} & \multicolumn{4}{c}{$\|\p - \p^*\|$}                                                                               \\
        $n$                       & $m$                         & $\rho$ & \logbar{}                & \logbarpcg{}                        & \pr{}  & \tat{} & \logbar{}  & \logbarpcg{} & \pr{}    & \tat{}             \\
\midrule
\multirow[c]{6}{*}{2,000} & \multirow[c]{2}{*}{10,000}  & -0.9   & 1.22                     & 1.91                                & 3.64   & 16.69  & {1.76e-09} & {4.41e-09}   & 9.82e-06 & 9.95e-06           \\
        \addlinespace[2pt]\cline{3-11}\addlinespace[2pt]
                                  &                             & 0.9    & 1.87                     & 2.49                                & 6.81   & 85.97  & {1.99e-10} & {8.98e-11}   & 9.66e-06 & 9.93e-06           \\
        \addlinespace[2pt]\cline{2-11}\addlinespace[2pt]
                                  & \multirow[c]{2}{*}{50,000}  & -0.9   & 5.77                     & 7.92                                & 15.52  & 73.08  & {2.37e-10} & {4.41e-09}   & 9.44e-06 & 9.77e-06           \\
        \addlinespace[2pt]\cline{3-11}\addlinespace[2pt]
                                  &                             & 0.9    & 7.85                     & 7.65                                & 28.68  & -      & {6.00e-11} & {4.40e-09}   & 8.30e-06 & 1.54e-05$^\dagger$ \\
        \addlinespace[2pt]\cline{2-11}\addlinespace[2pt]
                                  & \multirow[c]{2}{*}{100,000} & -0.9   & 11.64                    & 15.98                               & 31.82  & 135.02 & {1.13e-10} & {4.41e-09}   & 6.72e-06 & 9.71e-06           \\
        \addlinespace[2pt]\cline{3-11}\addlinespace[2pt]
                                  &                             & 0.9    & 13.52                    & 15.54                               & 52.87  & -      & {1.12e-09} & {4.41e-09}   & 9.11e-06 & 8.37e-05$^\dagger$ \\
        \addlinespace[2pt]\cline{2-11}\addlinespace[2pt]
        \multirow[c]{6}{*}{5,000} & \multirow[c]{2}{*}{10,000}  & -0.9   & 3.17                     & 4.76                                & 10.40  & 44.68  & {3.89e-10} & {1.76e-09}   & 4.86e-06 & 9.28e-06           \\
        \addlinespace[2pt]\cline{3-11}\addlinespace[2pt]
                                  &                             & 0.9    & 5.03                     & 6.45                                & 16.64  & -      & {3.80e-11} & {3.62e-11}   & 9.46e-06 & 9.49e-05$^\dagger$ \\
        \addlinespace[2pt]\cline{2-11}\addlinespace[2pt]
                                  & \multirow[c]{2}{*}{50,000}  & -0.9   & 17.15                    & 26.36                               & 44.65  & 195.73 & {2.20e-11} & {1.76e-09}   & 9.50e-06 & 9.79e-06           \\
        \addlinespace[2pt]\cline{3-11}\addlinespace[2pt]
                                  &                             & 0.9    & 19.77                    & 25.11                               & 83.69  & -      & {5.19e-10} & {1.76e-09}   & 8.15e-06 & 4.52e-04$^\dagger$ \\
        \addlinespace[2pt]\cline{2-11}\addlinespace[2pt]
                                  & \multirow[c]{2}{*}{100,000} & -0.9   & 32.94                    & 59.79                               & 121.84 & -      & {1.36e-08} & {1.76e-09}   & 6.55e-06 & 5.18e-04$^\dagger$ \\
        \addlinespace[2pt]\cline{3-11}\addlinespace[2pt]
                                  &                             & 0.9    & 49.76                    & 62.41                               & 190.05 & -      & {1.23e-10} & {1.76e-09}   & 1.07e-05 & 7.53e-04$^\dagger$ \\

        \bottomrule
    \end{longtable}
\end{table}
\begin{figure}[h]
    \centering
    \subcaptionbox{$\rho=-0.9$}{
        \includegraphics[width=0.33\textwidth]{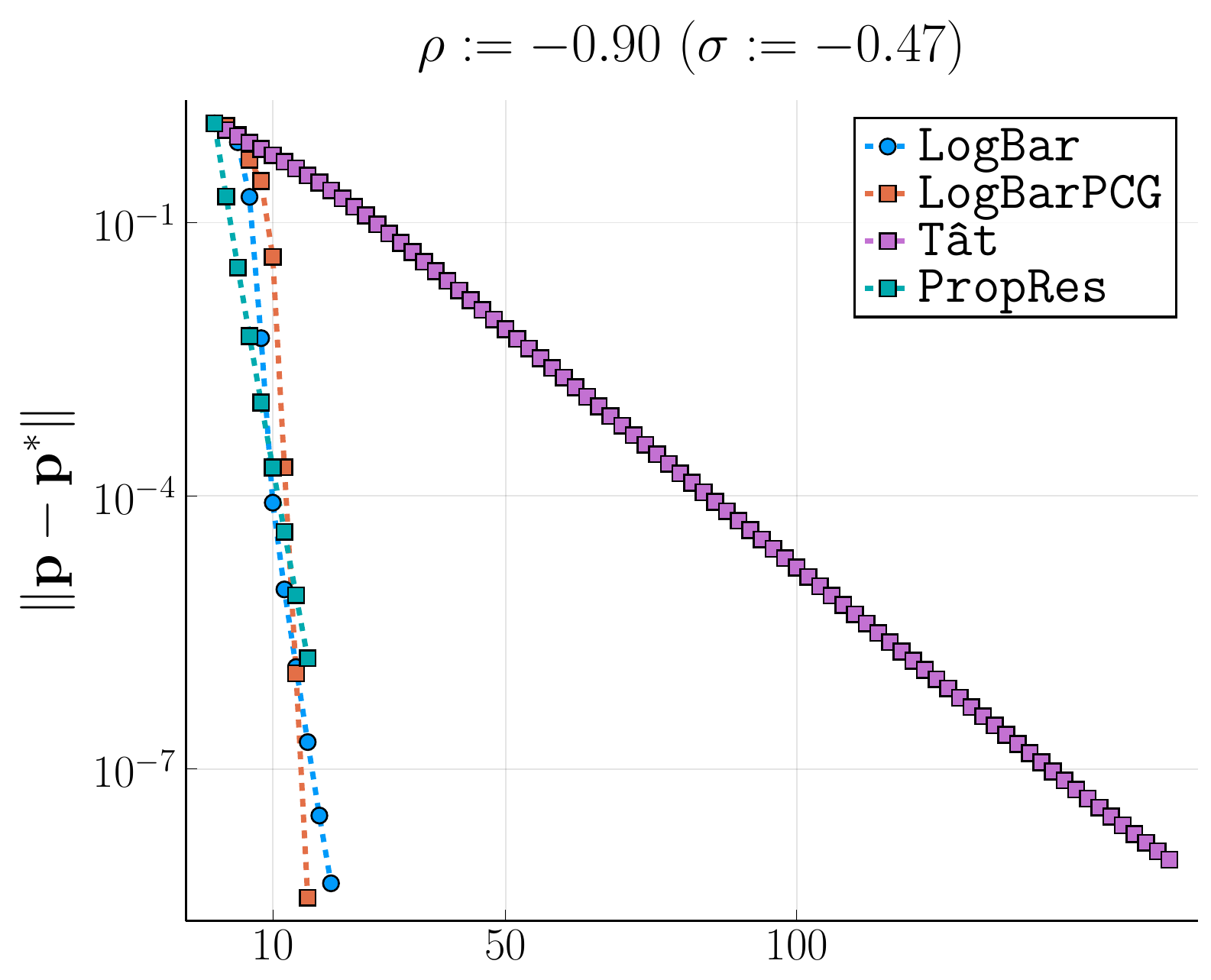}
    }
    \subcaptionbox{$\rho=0.9$}{
        \includegraphics[width=0.33\textwidth]{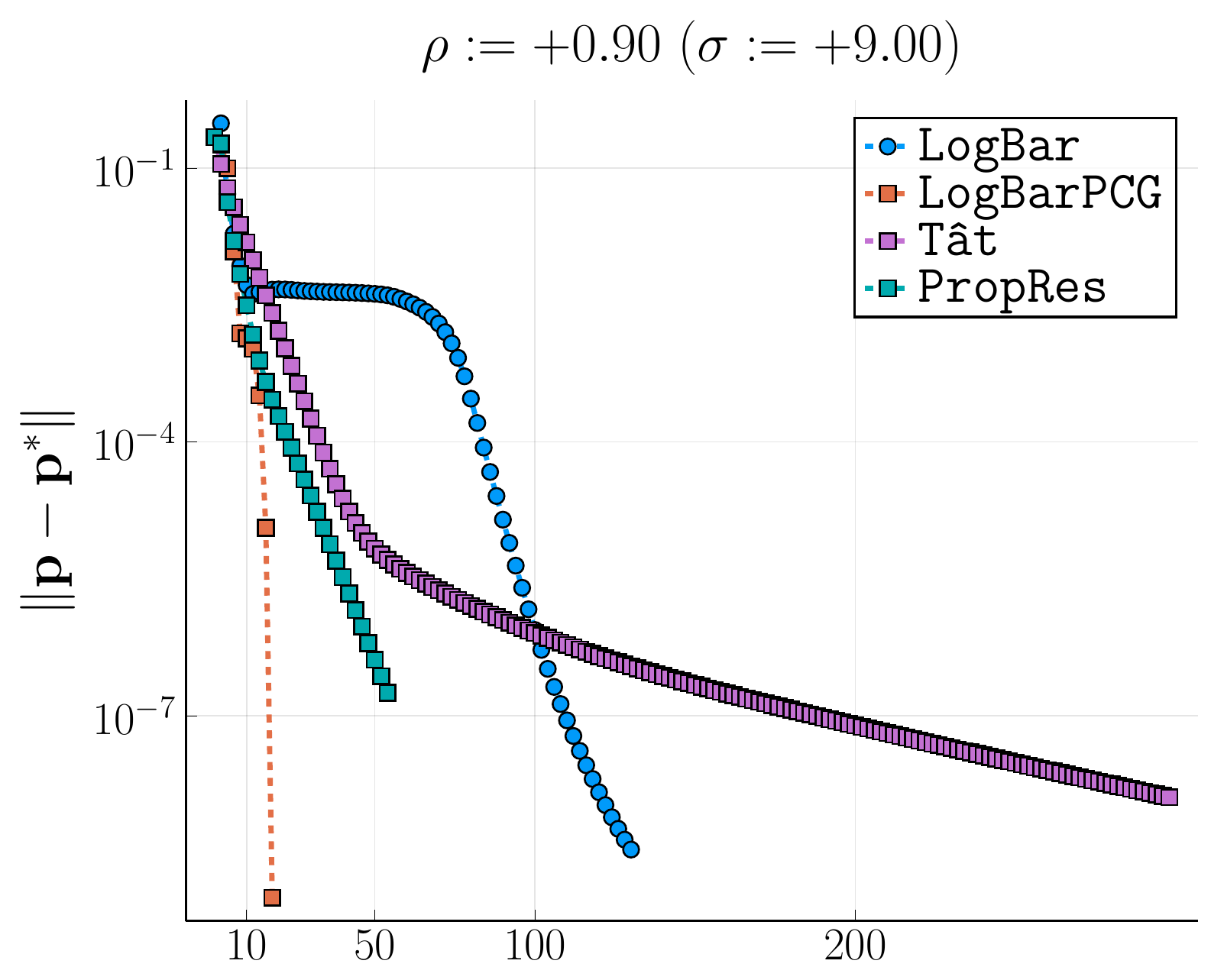}
    }
    \caption{Convergence behavior of the IPMs and the First-Order Methods in the CES economy using the MovieLens 32M dataset. Parameters: $n = 1,000$, $m = 10,000$.
    }\label{fig.ml-32m.compare}
\end{figure}

We next test the proposed methods on the MovieLens 32M dataset\endnote{The dataset is available at \url{https://grouplens.org/datasets/movielens/32m/}.}. Specifically, the movie ratings are used to construct the coefficients, $\c_1, \ldots, \c_m$ in the utility function.
We use the first 10,000 users and 1,000 items in the dataset.
A notable feature of this real-world dataset, compared to the instances generated randomly, is the sparsity. In this part of the experiment, the sparsity level is around $\tau=0.0018$. As a result, our preliminary finding is that the IPMs, \logbar{} (with the DR$1$ approximation), can be worse than \pr{} when in the substitution regime ($\rho > 0$) because of a slow start in the initial phase. In comparison, \logbarpcg{} preserves the fast speed of global and local convergence. For completeness, we present the results in \Cref{fig.ml-32m.compare}.

\section{Conclusion}\label{sec.conclusion}
In this paper, we show that the classical Smale's process can be improved as the decentralized interior-point methods for the equilibrium computation of the Fisher exchange market. The revitalization is due to the analysis of the logarithmic utility maximization, which enables us to design second-order t\^atonnement processes with complexity guarantees.
Overall, these second-order t\^atonnement processes, i.e., the decentralized IPMs, have the same information dependency as the classical t\^atonnement process, which is an improvement upon the informationally demanding Smale's process.
On the other hand, these methods significantly improve the centralized IPMs \citep{yePathArrowDebreu2008} that require $\cO(\sqrt{mn}(m+n)^3\log(\tfrac{1}{\epsilon}))$ arithmetic operations.
Still, we provide a few future directions.

For real-world applications, the revealed preferences are typically very sparse, which can weaken the convergence speed of the IPMs.
Smart segmentations and low-rank approximation schemes such as \cite{kroerComputingLargeMarket2022} may improve the speed of convergence of the IPMs.  Furthermore, it is clear that the current implementation can be adapted to GPU hardware, which has been adopted to solve very large-scale linear economic problems \citep{liuPDHCGScalableFirstorder2025}.

On the theory side, we only present the complexity analysis of \Cref{alg.logbarrier} (\logbar{}) using the SLC conditions, which is essentially a short-step small neighborhood IPM. Whether a predictor-corrector \logbar{} can be designed here could be interesting. The techniques in \cite{mehrotraImplementationPrimaldualInterior1992,mizunoAdaptivestepPrimaldualInteriorpoint1993,gondzioMultipleCentralityCorrections1996} and more recently \cite{aiO$sqrtnL$Iteration2005,nesterovParabolicTargetSpace2008,enagyNewInteriorpointAlgorithm2024} might be helpful. Furthermore, it would be interesting to see whether the similar analyses apply to the Arrow-Debreu model \citep{arrowStabilityCompetitiveEquilibrium1958}.
These exciting questions form the subject of our further research.
 
\section{Code and Data Disclosure}\label{sec:code-disclosure}The code and data to support the numerical experiments in this paper can be found at URL.

\begingroup\parindent 0pt \parskip 0.0ex \def\enotesize{\normalsize} \theendnotes\endgroup

\ACKNOWLEDGMENT{Part of this work was done when Chuwen Zhang was a PhD candidate at the Shanghai University of Finance and Economics. We would like to thank Wanyu Zhang and Prof. Shuzhong Zhang for carefully reading the preliminary version of this paper. We also thank the anonymous reviewers for their helpful comments and suggestions.}

\bibliographystyle{informs2014}

\clearpage
\begin{APPENDICES}
    
\section{Analysis of the Logarithmic Utility Maximization Problem}\label{appendix.sec.logump}
\paragraph{Further notations}
Denote $\bbN,\bbR$ as the set of natural and real numbers, respectively.
For a real number $a\in \real$, we let $[a]_+ = \max\{a, 0\}$ be the non-negative part of $a$.
For a vector $\x = [x^{(1)}, ..., x^{(n)}]$, $[\x]^p$ denotes the element-wise powers: $[\x]^p = [(x^{(1)})^p, ..., (x^{(n)})^p].$ The symbol $\Delta_n$ means the $n$-dimensional simplex: $\Delta_n = \{\x: \inner{\x}{\1}=1,\x\in\real_+^n\}$.
The operator norm for some $\bA: \real_+^n \mapsto \real_+^n$ is denoted by $\|\bA\| := \sup_{\|\p\|\le 1} \|\bA\p\|^*$. For a self-adjoint operator $\bA$, $\lambda_{\min}(\bA)/\lambda_{\max}(\bA)$ refer to its minimum / maximum eigenvalues. The condition number is denoted as $\kappa_2(\bA) = \tfrac{\lambda_{\max}(\bA)}{\lambda_{\min}(\bA)}.$
We use $\|\p\|_\bA$ to denote $\inner{\p}{\bA\p}^{\tfrac{1}{2}}$ and $\|\p\|^*_{\bA}$ for $\inner{\p}{\bA^\inv\p}^{\tfrac{1}{2}}$.
The underlying space $\real_+^n$ is omitted if the context is clear. For a $p$-th order differentiable function $f$, let $\nabla^p f(\p)[\h_1, ..., \h_p]$ be the $p$-th order directional differential along $(\h_1, ..., \h_p)$ at point $\p$; $\nabla f$ and $\nabla^2 f$ refer to the gradient and Hessian of $f$, respectively.
Finally, we use the notation $\otimes$ to denote the tensor product.
\subsection{Basic facts of the logarithmically homogeneous functions}
Suppose that $v(\x)$ is a $d$-logarithmically homogeneous function defined as
\begin{eqalign}{eq.logump.def.appendix}
    v(\alpha\x) & = -\log(u(\alpha\x)) = v(\x) - d\log(\alpha),
\end{eqalign}
then it satisfies the following properties \citep{nesterovSelfscaledBarriersInteriorpoint1997}:
\begin{subalign}{eq.loghomof}
    \label{eq.logump.1} \inner{\nabla v(\x)}{\x}=-d, \\
    \label{eq.logump.2}\nabla v(\alpha\x)=\frac{1}{\alpha} \nabla v(\x), \quad \nabla^2 v(\alpha\x)=\frac{1}{\alpha^2} \nabla^2 v(\x), \\
    \label{eq.logump.3}\nabla^2 v(\x) \x=-\nabla v(\x), \quad \nabla^3v(\x)[\x]=-2 \nabla^2 v(\x), \\
    \label{eq.logump.4}\inner{\nabla^2 v(\x) \x}{\x}=d, \quad\inner{\nabla v(\x)}{(\nabla^2 v(\x))^\inv \nabla v(\x)}=d.
\end{subalign}
Briefly, \eqref{eq.loghomof} can be derived by differentiating \eqref{eq.logump.def.appendix} and \eqref{eq.logump.1} on either $\alpha$ or $\x$. For the conjugate function $v_*$ (and the dual function $f$), note $-\nabla v(\x) \in \cX_*^\circ$ and $-\nabla v_*(\p) \in \cX^\circ$ (here our case is simplified as $\cX_* = \cX = \cP$), we would have:
\begin{subalign}{eq.loghomofdual}
    \label{eq.logump.dual.0} v_*(-\nabla v(\x)) = \inner{\nabla v(\x)}{\x} - v(\x) = -d - v(\x), \\
    \label{eq.logump.dual.1} v(-\nabla v_*(\p)) = -d - v_*(\p), \\
    \label{eq.logump.dual.2}\nabla v_*(-\nabla v(\x)) = -\x, \quad \nabla v(-\nabla v_*(\p)) = -\p, \\
    \label{eq.logump.dual.3}\nabla^2 v(-\nabla f(\p))= (\nabla^2 f(\p))^\inv,
    \quad \nabla^2 f(-\nabla v(\x)) = (\nabla^2 v(\x))^\inv, \\
    \label{eq.logump.dual.4} v(\x) + v_*(\p) \ge -d + d\log(d) - d\log(\inner{\p}{\x}).
\end{subalign}
Note \eqref{eq.logump.dual.0} can be implied by the so-called conjugate (sub)-gradient theorem (e.g., \cite[Theorem 4.20]{beckFirstOrderMethodsOptimization2017}).
Take $\y = -\nabla v(\x)$, then by convexity,
\begin{align*}
                  & v(\z) \ge v(\x) + \inner{-\y}{\z-\x}, \forall \z\in\cX,                      \\
    ~\Rightarrow~ & \inner{-\y}{\x} - v(\x) \ge \inner{-\y}{\z} - v(\z)                          \\
    ~\Rightarrow~ & \inner{-\y}{\x} - v(\x) \ge \max_{\z\in\cX}\inner{-\y}{\z} - v(\z)= v_*(\y).
\end{align*}
By Fenchel's inequality, $v_*(\y) +  v(\x)  \ge \inner{-\y}{\x}$, hence we conclude \eqref{eq.logump.dual.0}, and the same holds for $v_*$ \eqref{eq.logump.dual.1}. Differentiating over $\y$, we have
$\nabla v_*(\y) = -\x$
and similarly \eqref{eq.logump.dual.2} holds. Finally, by $\nabla v_*(-\nabla v(\x)) = -\x$, differentiating over $\x$ yields \eqref{eq.logump.dual.3}:
\begin{align*}
    \nabla^2 v_*(-\nabla v(\x))\nabla^2 v(\x) = -\bI.
\end{align*}

\subsection{Proof of \Cref{lem.equiv.conj.purt}}

By definition, for any $\p \in \cP$, $t > 0$,
\begin{align*}
    v_*(t\p) & = \max_{\x \in \cX} \left\{ - \inner{t\p}{\x} - v(\x) \right\}            \\
             & = \max_{\x \in \cX} \left\{ -\inner{\p}{t\x} - v(t\x) - d\log(t) \right\} \\
             & = v_*(\p) - d\log(t).
\end{align*}
Hence, $v_*$ is a $d$-logarithmically homogeneous function. Consider Lagrangian of problem \eqref{eq.logump}
$$
    L(\x, \lambda) = -v(\x) - \lambda \inner{\p}{\x} + \lambda w,
$$
and the optimality condition indicates that there exists some $(\bar \lambda, \bar \x, \bar\s) \in \real_+\times \cX\times \cX_*$ such that
\begin{subalign}{eq.logump.bestresp}
    \label{eq.br.fo1}     \nabla v(\bar\x) + \bar\lambda \p - \bar \s = 0, ~\bar\bX\bar\s = 0,\\
    \real_+ \ni \bar\lambda \perp w - \inner{\p}{\bar\x} \in \real_+.
\end{subalign}
If we consider the standard UMP, the optimal solution $(\lambda',\x',\s')\in \real_+\times \cX\times \cX_*$ should satisfy
\begin{subalign}{eq.ump.bestresp}
    -\nabla u(\x') + \lambda' \p - \s' = 0, ~\bX'\s' = 0,\\
    \real_+ \ni \lambda' \perp w - \inner{\p}{\x'} \in \real_+.
\end{subalign}
Recall that $\nabla v(\bar\x) = -\tfrac{\nabla u(\bar\x)}{u(\bar\x)}$, $(\bar \lambda, \bar \x, \bar\s)$ solves \eqref{eq.logump.bestresp}, then $(u(\bar \x)\bar\lambda, \bar\x, u(\bar \x)\bar\s)$ also solves \eqref{eq.ump.bestresp}, indicating the equivalence of best-response mappings of UMP and $\lump{}$.
On the other hand, by \eqref{eq.br.fo1} and multiply both sides of the first equation by $\bar \x$, we have that
$$\inner{\bar\x}{\nabla v(\bar\x)} + \bar\lambda \inner{\bar\x}{\p} = 0.$$
By Euler's theorem on homogeneous functions (see, e.g., \cite[Theorem 3.1.21]{nesterovLecturesConvexOptimization2018}),
\begin{align*}
    \inner{\bar\x}{\nabla u(\bar\x)} = d u(\bar\x) \qlq  \inner{\bar\x}{\nabla v(\bar\x)} = -d \qlq -d + \bar\lambda \inner{\bar\x}{\p} = 0.
\end{align*}
Since $d \neq 0$ (indeed, it is only possible that $d > 0$),
the last equation implies that $\bar\lambda > 0$ and $\inner{\bar\x}{\p} = w$.
That is, the budget must be active, and the dual variable $\bar\lambda$ is a constant, $\bar\lambda = \frac{d}{w}$.
By strong duality,
\begin{align*}
    f(\p) & = \max_{\x \in \cX} L(\x, \bar\lambda) =
    \max_{\x \in \cX} \left\{- \inner{\bar\lambda \p}{\x} - v(\x) \right\} + \bar\lambda w \\
          & = v_*(\bar\lambda \p) + d = v_*(\p) + d - d\log(\bar\lambda)                   \\
          & = v_*(\p) + d - d\log(d/w).
\end{align*}
So that $f(\p)$ deviates from $v_*(\p)$ by some constant.
Finally, from Danskin's theorem (e.g. \cite[Theorem 9.26]{shapiroLecturesStochasticProgramming2021}),
\begin{align*}
    \nabla f(\p) = \nabla v_*(\p) = -\bar\lambda \bar\x = -\frac{d}{w}\bar\x.
\end{align*}
This completes the proof.
\hfill\qed

\subsection{Proof of \Cref{lem.additive.diagonal}}

For all $j \in \cJ$, $\gamma^{(j)}(\cdot) \in \cH^{r}$.
By Euler's theorem,
$$
    \theta'_j(x^{(j)}) x^{(j)}  = r\theta_j(x^{(j)}), \quad \Longrightarrow \quad  \bX\nabla \btheta  = r\btheta.
$$
Similarly, since $\theta'_j(\cdot) \in \cH^{r- 1}$ (e.g., see \cite{mas-colellMicroeconomicTheory1995}), we have
$$
    \theta''_j(x^{(j)}) x^{(j)}  = (r- 1)\theta'_j(x^{(j)}) \qlq \bX\nabla^2 \btheta[\x, \cdot]  = r(r-1) \btheta.
$$
By induction, all the quantities below belong to the unit simplex by letting $\gamma^{(j)} = \tfrac{\theta^{(j)}}{\inner{\btheta}{\1}}$:
\begin{eqalign}{eq.parallel.theta2}
    \inner{\bgamma}{\1} = 1, \quad \frac{\bX\nabla \btheta}{r\inner{\btheta}{\1}}
    = \frac{\bX\nabla^2 \btheta[\x, \cdot]}{r(r-1)\inner{\btheta}{\1}}
    = \frac{\bX\nabla^3 \btheta[\x, \x, \cdot]}{r(r-1)(r-2)\inner{\btheta}{\1}}
    = \frac{\btheta}{\inner{\btheta}{\1}} = \bgamma \in \Delta_n.
\end{eqalign}
Now we inspect the utility function $u(\x)$. Note that
\begin{align*}
    \nabla u(\x)   & = k\inner{\btheta}{\1}^{k-1}\nabla \btheta,  \\
    \nabla^2 u(\x) & = k\inner{\btheta}{\1}^{k-1}\nabla^2 \btheta
    + k(k-1)\inner{\btheta}{\1}^{k-2}\nabla \btheta\nabla \btheta^{\top}.
\end{align*}
By \eqref{eq.parallel.theta2}, we write $\nabla \btheta = r\inner{\btheta}{\1} \bX^{-1} \bgamma$, and $\nabla^2\btheta = r(r-1) \inner{\btheta}{\1}\bX^\inv \bGamma\bX^\inv$. Recall $u(\x) = \inner{\btheta}{\1}^k$, and it follows
\begin{align*}
    \nabla u(\x)  & =  k r \cdot u(\x) \bX^\inv\bgamma             \\
    \nabla^2u(\x) & = k r(r-1)\cdot u(\x) \bX^\inv \bGamma\bX^\inv
    + k(k-1)r^2 u(\x)\bX^\inv\bgamma\bgamma^{\top}\bX^\inv.
\end{align*}
For any vector $\v\in\real^n$, observe that
$$
    \inner{\v}{(\bGamma - \bgamma\bgamma^{\top})\v}
    = \sum_{j=1}^n \gamma^{(j)} (v^{(j)})^2 - \left(\sum_{j=1}^n \gamma^{(j)}v^{(j)}\right)^2 \ge 0.
$$
The right-hand side represents the variance of $\v$ under the probability $\bgamma$, and thus the value is nonnegative. Hence, $\bGamma \succeq \bgamma\bgamma^{\top}$.
By $u(\x)>0$, the condition for $u(\x)$ being concave is,
$$k r(r-1) < 0 \quad \text{and} \quad - k r(r-1) > k(k-1)r^2.$$
This is equivalent to saying that:
$r \in (0,1], k\in (0, \tfrac{1}{r}] ~\text{or}~ r < 0, k \in [\tfrac{1}{r}, 0).$
This also gives $d = k\cdot r > 0$.
\hfill\qed

\subsection{Proof of \Cref{lem.pd.bestresp}}
In view of \eqref{eq.f.gradient}, \eqref{eq.loghomof}, and \eqref{eq.loghomofdual} and the fact that $v$ is $d$-logarithmically homogeneous, we have that
\begin{align*}
    -\p & \stackrel{\eqref{eq.loghomofdual}}{=} \nabla v(-\nabla f(\p))
    \stackrel{\eqref{eq.f.gradient}}{=} \nabla v(\tfrac{d}{w}\x(\p))
    \stackrel{\eqref{eq.loghomof}}{=} \frac{w}{d}\nabla v(\x(\p)).
\end{align*}
Since the dual gradient is given by $\nabla v(\x) = -k\frac{\nabla \btheta}{\inner{\btheta}{\1}} = -d\bX^\inv\bgamma$, and we conclude that $\bgamma = \frac{\bX(\p)\p}{w}.$
By \eqref{eq.parallel.theta},
\begin{align*}
    \frac{\nabla^2 \btheta}{\inner{\btheta}{\1}} = r(r-1) \bX^\inv\bGamma\bX^\inv = \frac{r(r-1)}{w^2} \bP \bGamma^\inv\bP,
\end{align*}
so we have
\begin{align*}
    \nabla^2 v(\x(\p)) & = -k\frac{\nabla^2 \btheta}{\inner{\btheta}{\1}}
    + k\left(\frac{\nabla \btheta}{\inner{\btheta}{\1}}\right)\left(\frac{\nabla \btheta}{\inner{\btheta}{\1}}\right)^{\top} \\
                       & = \frac{kr(1-r)}{w^2} \bP \bGamma^\inv\bP +\frac{kr^2}{w^2} \p\p^{\top}.
\end{align*}
Similarly, reusing \eqref{eq.loghomof} and \eqref{eq.loghomofdual},
\begin{align*}
    \nabla^2 f(\p) & = \left[\nabla^2 v(-\nabla f(\p))\right]^\inv
    = \left[\nabla^2 v(\tfrac{d}{w}\x(\p))\right]^\inv = \left[\tfrac{w^2}{d^2}\nabla^2 v(\x(\p))\right]^\inv \\
                   & = \frac{d^2}{w^2}\left[\nabla^2 v(\x(\p))\right]^\inv.
\end{align*}
Finally, using Sherman-Morrison-Woodbury formula and $\bgamma \in \Delta_n$ yields
\begin{align*}
    \nabla^2 f(\p)
     & = \frac{d^2}{w^2}
    \left[kr(1-r)\frac{1}{w^2} \bP \bGamma^\inv\bP +\frac{kr^2}{w^2} \p\p^{\top}\right]^\inv
    = \frac{d^2}{kr^2} \bP^\inv
    \left[\frac{1-r}{r}\bGamma^\inv+\1\1^{\top}\right]^\inv\bP^\inv                                                                                                                          \\
     & = \frac{d^2}{kr^2} \frac{r}{(1-r)} \bP^\inv \left(\bGamma - r \bgamma\bgamma^{\top}\right)\bP^\inv = \frac{d}{(1-r)} \bP^\inv \left(\bGamma - r \bgamma\bgamma^{\top}\right)\bP^\inv.
\end{align*}
This completes the proof.
\hfill\qed

\subsection{Proof of \Cref{lem.bounded.directdev}}
We first notice that the directional derivatives of $v(\x)$ can be written as the products of $\bgamma$ and the affine-scaled direction $\bX^\inv\h$ ($[\cdot]^p$ denotes the element-wise powers). Let us introduce the following notation that will be used throughout the proof:
\begin{equation}\label{eq.tau.def}
    \tau_p := \inner{\bgamma}{\left[\bX^\inv\h\right]^p} = \sum_{j\in\cJ} \gamma^{(j)} (h^{(j)}/x^{(j)})^p, p \in \bbN_+
\end{equation}
Apparently, $|\tau_p| \le \|\h\|_{\x}^p.$
Compute the derivatives of $v(\x)$,
\begin{subarr}{eq.3rd.vx}{cl}
    \frac{1}{k}\nabla v(\x) &= -\frac{\nabla \btheta}{\inner{\btheta}{\1}}, \\
    \frac{1}{k}\nabla^2 v(\x) &= -\frac{\nabla^2 \btheta}{\inner{\btheta}{\1}}
    + \left(\frac{\nabla \btheta}{\inner{\btheta}{\1}}\right)\otimes\left(\frac{\nabla \btheta}{\inner{\btheta}{\1}}\right), \\
    \frac{1}{k}\nabla^3 v(\x) &= -\frac{\nabla^3 \btheta}{\inner{\btheta}{\1}}
    + 3\left(\frac{\nabla^2 \btheta}{\inner{\btheta}{\1}}\right)\otimes\left(\frac{\nabla \btheta}{\inner{\btheta}{\1}}\right) - 2 \left(\frac{\nabla \btheta}{\inner{\btheta}{\1}}\right)^{\otimes 3}.
\end{subarr}
By combining with Lemma \ref{lem.additive.diagonal} and \eqref{eq.tau.def}, we have
\begin{subalign}{eq.uh.defs}
    \label{eq.uh1} H_1(\h) =\frac{1}{d}\nabla v(\x)[\h] & = - \tau_1,\\
    \label{eq.uh2} H_2(\h) =\frac{1}{d}\nabla^2 v(\x)[\h,\h] &= (1-r)\tau_2 + r \tau_1^2,\\
    \label{eq.uh3} H_3(\h) =\frac{1}{d}\nabla^3 v(\x)[\h,\h,\h]
    &= -(r-1)(r-2)\tau_3 + 3r(r-1)\tau_2\tau_1 - 2 r^2\tau_1^3 \\
    &= -(r-1)(r-2)\tau_3 + \tau_1(-2r^2\tau_1^2 + 3r(r-1)\tau_2).
\end{subalign}
Let $\z := \bX^{-1} \h$. Since $\bgamma \in \Delta_n$. Then when $p \ge 2$,
\begin{align*}
    |\tau_p| = \left| \sum_{i=1}^n \gamma_i z_i^p \right|
    \le \sum_{i=1}^n \gamma_i |z_i|^p \le \sum_{i=1}^n |z_i|^p = \|\z\|_p^p \le \|\z\|_2^p = \|\h\|_{\x}^p.
\end{align*}
When $p = 1$, we obviously have
\begin{align*}
    |\tau_1| = \left| \sum_{i=1}^n \gamma_i z_i \right| \le \|\bgamma\|_2 \|\z\|_2 \le \|\z\|_2 = \|\h\|_{\x}.
\end{align*}
Then we have $|\nabla v(\x)[\h]| \le d\|\h\|_{\x}$. By \eqref{eq.uh2}, if $r \in (0,1]$, $k \in [0, \tfrac{1}{r}]$,
\begin{align*}
    H_2(\h) = (1-r)\tau_2+r\tau_1^2 \le (1-r)\|\h\|_{\x}^2 + r\|\h\|_{\x}^2 \le \|\h\|_{\x}^2;
\end{align*}
in the other case, $r < 0, k \in [\tfrac{1}{r}, 0]$, we have
\begin{align*}
    H_2(\h) = (1-r)\tau_2+r\tau_1^2 \le (1-r) \tau_2 \le (1-r) \|\h\|_{\x}^2.
\end{align*}
By \eqref{eq.uh3}, again if $r \in (0,1]$, $k \in [0, \tfrac{1}{r}]$,
\begin{align*}
    H_3(\h) & \le (1-r)(2-r)|\tau_3| + 3 r(1-r)\tau_2|\tau_1| + 2 r^2|\tau_1|^3 \le 2 \|\h\|_{\x}^3.
\end{align*}
otherwise if $r < 0, k \in [\tfrac{1}{r}, 0]$,
\begin{align*}
    H_3(\h) & \le (1-r)(2-r)|\tau_3| + 3 r(r-1)\tau_2|\tau_1| + 2 r^2|\tau_1|^3 \\
            & \le (6r^2-6r+2) \|\h\|_{\x}^3.
\end{align*}
This completes the proof.
\hfill\qed

\subsection{Proof of \Cref{thm.taylor.v}}

First observe that
$\diag(\x + t\h) = (\bI + \diag(t\bX^\inv\h))\bX,$
and thus it follows
\begin{align*}
    \lambda_{\min}(\bI + t\diag(\bX^\inv\h)) \ge 1 - t\|\bX^\inv\h\|_{\infty}, \\
    ~ \diag(\x + t\h)\succeq (1 - t\|\bX^\inv\h\|_{\infty})\bX \succeq (1 - t\|\h\|_{\x})\bX.
\end{align*}
Hence, by \cref{lem.bounded.directdev}, for any $t \in [0, 1]$,
\begin{align*}
    |\nabla^3 v(\x + t\h)[\h,\h,\h]| \le T_u \|\h\|^3_{\x + t\h}
    = T_u \|\diag(\x + t\h)^\inv\h\|^3
    \le T_u \frac{\|\h\|_{\x}^3}{(1 - t\|\h\|_{\x})^3}.
\end{align*}
Let $\rho = \|\h\|_{\x} < 1$, and $\psi(t) = v(\x + t\h)$ be a scalar function. Using the Taylor expansion of $\psi$, we have
\begin{align*}
    \psi(1) - \psi(0) - \psi^\prime(0) - \frac{1}{2}\psi^{\prime\prime}(0) = \frac{1}{2}\int_0^1 \psi^{\prime\prime\prime}(t) (1-t)^2 \diff t,
\end{align*}
where $\psi^\prime(\cdot)$ denotes the derivative. It implies that
\begin{subalign}{}
    \nonumber &\left|  v(\x+\h) - v(\x) - \nabla v(\x)[\h] - \tfrac{1}{2}\nabla^2 v(\x)[\h,\h]\right|
    \nonumber   = \left|\int_0^1 \frac{(1-t)^2}{2} \nabla^3 v(\x+t\h)[\h,\h,\h] \diff t\right| \\
    \nonumber &\le \int_0^1 T_u \frac{(1-t)^2}{2} \frac{\rho^3}{(1 - t\rho)^3} \diff t = T_u \rho^3 \int_0^1 \frac{(1-t)^2}{2(1 - t\rho)^3}  \diff t \\
    \nonumber &= T_u \rho^3 \left(-\frac{1}{4} \cdot \frac{\rho (2 + \rho) + 2 \log(1 - \rho)}{\rho^3}\right) = -\frac{1}{4}T_u \left( \rho (2 + \rho) + 2 \log(1 - \rho)\right) \\
    &= -\frac{1}{2}T_u \left(\rho +  \tfrac{1}{2}\rho^2 +  \log(1 - \rho)\right) \\
    \label{eq.taylor.v.karmarkar} &\le \frac{T_u}{2} \frac{\rho^3}{3(1-\rho)}.
\end{subalign}
The last inequality \eqref{eq.taylor.v.karmarkar} is well-known in the literature of IPMs (e.g. \cite{karmarkarNewPolynomialtimeAlgorithm1984,haeserOptimalityConditionComplexity2019}).
\hfill\qed

\subsection{Proof of \Cref{lem.bounded.directdev.dual}}
Similarly to our construction in the primal space \eqref{eq.tau.def}, for any direction $\q$, let
\begin{equation}\label{eq.a.def.new}
    \f = \bP^\inv\q, \quad a_1 = \inner{\bgamma}{\f}, \quad a_2 = \inner{\bgamma}{[\f]^2}, \quad a_3 = \inner{\bgamma}{[\f]^3}.
\end{equation}
Recall from \eqref{eq.f.gradient} and \eqref{eq.f.hess}, we have
$$\nabla \x(\p) = - \frac{w}{d}\nabla^2 f(\p), \quad \nabla^2 v(\x(\p))^\inv = \frac{w^2}{d^2}\nabla^2 f(\p),$$
Hence, $|\nabla f(\p)[\q]| = |\inner{\q}{-\frac{d}{w}\x}| = |\inner{\bP^\inv\q}{d\bgamma}| = d |a_1| \le d \|\q\|_{\p}.$
Besides,
\begin{eqalign}{}
    \nabla^2 f(\p)[\q,\q] &= \left(\tfrac{d}{1-r}\bP^\inv \left(\bGamma - r \bgamma\bgamma^{\top}\right)\bP^\inv\right)[\q,\q] \\
    & \le \tfrac{d}{1-r}\left\|\bGamma - r \bgamma\bgamma^{\top}\right\| \|\bP^\inv\q\|^2 \le \tfrac{d(1+|r|)}{1-r}\|\q\|_{\p}^2.
\end{eqalign}
Finally, the third-order directional derivative of $f(\p)$ is given by
\begin{eqalign}{eq.f.hesstimesq}
    \nabla^3 f(\p)[\q] & = \lim_{t\to 0} \frac{1}{t} \left(\nabla^2f(\p+t\q) - \nabla^2f(\p)\right)                                                                                        \\
    & \stackrel{\eqref{eq.f.hess}}{=} \frac{d^2}{w^2}\lim_{t\to 0} \frac{1}{t} \left(\nabla^2 v(\x(\p+t\q))^\inv - \nabla^2 v(\x(\p))^\inv\right)                       \\
    & = \frac{d^2}{w^2}\lim_{t\to 0} \frac{1}{t} \nabla^2 v(\x(\p))^\inv\left(\nabla^2 v(\x(\p))-\nabla^2 v(\x(\p+t\q))\right)\nabla^2 v(\x(\p+t\q))^\inv               \\
    & = \frac{d^2}{w^2}\nabla^2 v(\x(\p))^\inv \lim_{t\to 0} \frac{1}{t} \left(\nabla^2 v(\x(\p))-\nabla^2 v(\x(\p+t\q))\right)\lim_{t\to 0}\nabla^2 v(\x(\p+t\q))^\inv \\
    & = \frac{d^2}{w^2}\nabla^2 v(\x(\p))^\inv \lim_{t\to 0} \frac{1}{t} \left(\nabla^2 v(\x(\p))-\nabla^2 v(\x(\p+t\q))\right)\nabla^2 v(\x(\p))^\inv \\
    & = -\frac{d^2}{w^2}\nabla^2 v(\x(\p))^\inv \left(\nabla^3v(\x(\p))[\nabla\x(\p)\q]\right)\nabla^2 v(\x(\p))^\inv \\
    & \stackrel{\eqref{eq.f.hess}}{=} -\frac{d^2}{w^2}\left(\frac{w^2}{d^2}\nabla^2 f(\p)\right) \left(\nabla^3v(\x(\p))[- \frac{w}{d}\nabla^2 f(\p)\q]\right)\left(\frac{w^2}{d^2}\nabla^2 f(\p)\right) \\
    & =\frac{w^3}{d^3}\nabla^2 f(\p)\left(\nabla^3v(\x(\p))[\nabla^2 f(\p)\q]\right)\nabla^2 f(\p).
\end{eqalign}
Let $\h = \nabla^2 f(\p)\q$, since $\bgamma = \frac{1}{w}\bP\x(\p)$,
\begin{align*}
    (\bX(\p))^\inv \h = \bX(\p)^\inv\nabla^2 f(\p)\q
     & \stackrel{\eqref{eq.f.hess}}{=} \frac{d}{(1-r)}\bX(\p)^\inv \bP^\inv \left(\bGamma - r \bgamma\bgamma^{\top}\right)\bP^\inv\q \\
     & = \frac{d}{(1-r)}\frac{1}{w} \left(\bI - r \1\bgamma^{\top}\right)\bP^\inv\q.
\end{align*}
Hence, let us take $\bS = \left(\bI - r \1\bgamma^{\top}\right)$, and we have
\begin{eqalign}{eq.ub3.f3qqq.new}
    \nabla^3 f(\p)[\q,\q,\q]
    &=\tfrac{w^3}{d^3}\nabla^3v(\x(\p))[\nabla^2 f(\p)\q,\nabla^2 f(\p)\q,\nabla^2 f(\p)\q] \\
    &=\tfrac{w^3}{d^3}\nabla^3v(\x)[\bX\tfrac{d}{(1-r)}\tfrac{1}{w} \bS\f,\bX\tfrac{d}{(1-r)}\tfrac{1}{w} \bS\f,\bX\tfrac{d}{(1-r)}\tfrac{1}{w} \bS\f],  \\
    &=\tfrac{1}{(1-r)^3}\nabla^3v(\x)[\bX\bS\f,\bX\bS\f,\bX\bS\f].
\end{eqalign}
In view of \eqref{eq.3rd.vx}, we have,
\begin{eqalign}{eq.ub3.f3qqq.new.2}
    &\nabla^3 f(\p)[\q,\q,\q] = \frac{1}{(1-r)^3}\nabla^3v(\x)[\bX\bS\f,\bX\bS\f,\bX\bS\f] \\
    &\stackrel{\eqref{eq.uh3}}{=}
    \frac{d}{(1-r)^3}\left(
    -(1-r)(2-r)\inner{[\bS\f]^3}{\bgamma}
    -3(1-r)r\inner{[\bS\f]^2}{\bgamma}\inner{\bgamma}{\bS\f}
    -2r^2\inner{\bgamma}{\bS\f}^3
    \right).
\end{eqalign}
Then, since $\bS\f = \f - r\1\inner{\bgamma}{\f} = \f - ra_1\1$. We use following identities for 1-dimensional polynomials,
\begin{subarr}{}{ll}
    (q-ra)^2=q^2 - 2 ar q  + a^2 r^2, \\
    (q-ra)^3=q^3 - 3 ar q^2  + 3 a^2r^2 q  - a^3 r^3.
\end{subarr}
Note that $\bgamma \in \Delta_n$, the following holds,
\begin{subarr}{}{rl}
    \inner{\bS\f}{\bgamma} &= \inner{\f - ra_1\1}{\bgamma} = a_1 - ra_1 = (1-r)a_1 \\
    \inner{[\bS\f]^2}{\bgamma} &= \inner{[\f]^2 - 2(r a_1)\f + (r a_1)^2\1}{\bgamma} =a_2 - 2ra_1^2 + r^2a_1^2 \\
    \inner{[\bS\f]^3}{\bgamma} &= \inner{[\f]^3 - 3[\f]^2(r a_1) + 3\f (r a_1)^2 - (r a_1)^3\1}{\bgamma} = a_3 - 3ra_1a_2 + 3r^2a_1^3 - r^3a_1^3.
\end{subarr}
Hence,
\begin{align*}
    \nabla^3 f(\p)[\q,\q,\q]
    = & \frac{d}{(1-r)^3}\Big[(1-r)
    \left((r-2)a_3 + 3r a_1a_2 - 2r^2 a_1^3\right)\Big]              \\
    = & \frac{d}{(1-r)^2}\Big[(r-2)a_3 + 3r a_1a_2 - 2r^2 a_1^3\Big]
\end{align*}
Since $\gamma \in \Delta_n$, the terms $a_1^3, a_3, a_1a_2$ satisfy $|a_1|^3 \le |a_3|$ and $|a_1a_2| \le |a_3|$. Recall that $r < 1$. For the case $r \in [0,1)$, it holds that $|(r-2)a_3 + 3r a_1a_2 - 2r^2 a_1^3| \le (2-r + 3r + 2r^2)|a_3| = 2(1 + r + r^2)|a_3| \le 6|a_3|$; otherwise, it holds that $|(r-2)a_3 + 3r a_1a_2 - 2r^2 a_1^3| \le (2-r - 3r + 2r^2)|a_3| = 2(1 - 2r + r^2)|a_3| = 2(r-1)^2 |a_3|$ when $r < 0$. Consequently, we conclude
\begin{eqalign}{eq.ub3.f3qqq}
    \left|\nabla^3 f(\p)[\q,\q,\q]\right|
    &\le \max\left\{\frac{6d}{(1-r)^2}, 2d\right\}|a_3| \\
    &\le \max\left\{\frac{6d}{(1-r)^2}, 2d\right\}\|\q\|_{\p}^3.
\end{eqalign}
The proof is complete. \hfill\qed

\subsection{Proof of \Cref{thm.taylor.f}}
Using the same arguments for the primal space,
$\diag(\p + t\q) = (\bI + \diag(t\bP^\inv\q))\bP,$
and thus it follows
\begin{align*}
    \lambda_{\min}(\bI + t\diag(\bP^\inv\q)) \ge 1 - t\|\bP^\inv\q\|_{\infty}, \\
    ~ \diag(\p + t\q)\succeq (1 - t\|\bP^\inv\q\|_{\infty})\bP \succeq (1 - t\|\q\|_{\p})\bP.
\end{align*}
The rest of the proof is the same as that for $v(\x)$, cf. \eqref{eq.taylor.v.karmarkar}, using scaled boundedness of third-order  directional derivatives with coefficient $T_f$.

\subsection{Proof of \Cref{thm.taylor.x}}
Let us include a useful result for multi-linear form.
\begin{corollary}[Proposition 9.1.1 of \cite{nesterovInteriorPointPolynomialAlgorithms1994}]\label{lem.bounded.multilinear}
    Let $\bT[\h_1,...,\h_k]$ be a symmetric $k$-linear form.
    If the following is true for
    some positive-semidefinite quadratic form $\bA[\h_1,\h_2]$,
    $$
        |\bT[\h,...,\h]| \le \bA[\h,\h]^{\tfrac{k}{2}}.
    $$
    then for any $\h_j \in \real_+^n, j = 1, 2, ..., k$,
    \begin{eqalign}{eq.multilinear}
        |\bT[\h_1,...,\h_k]| \le \smallprod_{j=1}^{k}\bA[\h_j,\h_j]^{1/2}.
    \end{eqalign}
\end{corollary}

\subsubsection*{Proof of \Cref{thm.taylor.x}}
Observe that
\begin{align*}
    \x(\p+\q) - \x(\p) - \nabla\x(\p)[\q]
     & = -\frac{w}{d}\left(
    \nabla f(\p+\q) - \nabla f(\p) - \nabla^2 f(\p)[\q]
    \right)                                                           \\
     & = -\frac{w}{d}\int_0^1 (1-t)\nabla^3 f(\p+t\q)[\q,\q] \diff t.
\end{align*}
Hence,
\begin{align*}
     & \|\bP(\x(\p+\q) - \x(\p) - \nabla\x(\p)[\q])\|
    \le \frac{w}{d}\int_0^1 (1-t)\|\bP\nabla^3 f(\p+t\q)[\q,\q]\| \diff t                           \\
     & = \frac{w}{d}\int_0^1 (1-t) \sup_{\|\xi\|=1} \nabla^3 f(\p+t\q)[\q,\q, \bP\xi] \diff t       \\
     & = \frac{w}{d}\int_0^1 (1-t) \sup_{\|\bP^{-1}\q'\|=1} \nabla^3 f(\p+t\q)[\q,\q, \q'] \diff t.
\end{align*}
In view of \eqref{eq.ub3.f3qqq}, and \cref{lem.bounded.multilinear}, we have,
\begin{align*}
     & \|\bP(\x(\p+\q) - \x(\p) - \nabla\x(\p)[\q])\|                                                                                       \\
     & \le
    \frac{w}{d} \int_0^1 \sup_{\|\bP^{-1}\q'\|=1} (1-t)T_f \|\diag(\p+t\q)^\inv\q\|^2\|\diag(\p+t\q)^\inv\q'\|\diff t                       \\
     & \le \frac{w}{d} \int_0^1 \sup_{\|\bP^{-1}\q'\|=1} T_f \frac{(1-t)\varrho^2}{(1-t\varrho)^2}\frac{\|\bP^\inv\q'\|}{1-t\varrho}\diff t \\
     & = \frac{w}{d}  T_f\varrho^2\int_0^1 \frac{1-t}{(1-t\varrho)^3}\diff t = \frac{w}{d} \frac{T_f\varrho^2}{2(1-\varrho)}.
\end{align*}
This completes the proof.
\hfill\qed

\subsection{Proof of \Cref{lem.dom.23}}
Recall that we have the following observation of directional derivatives at some allocation $\x\in\cX$, cf. \eqref{eq.3rd.vx} and \eqref{eq.uh.defs}.
\begin{align*}
    H_3(\h) = \frac{1}{d}\nabla^3 v(\x)[\h,\h,\h]
     & = -(r-1)(r-2)\tau_3 + 3r(r-1)\tau_2\tau_1 - 2 r^2\tau_1^3.
\end{align*}
As a function of $\h$, it is homogeneous of degree 3 and $H_3(-\h) = - H_3(\h)$.
Similarly, the quadratic $H_2$ is also homogeneous:
\begin{align*}
    H_2(\h) & = \frac{1}{d} \nabla^2v(\x)[\h,\h] = (1-r) \tau_2 + r\tau_1^2
    \\&= (1-r)\inner{\bgamma}{[\bX^\inv\h]^2} + r\inner{\bgamma}{\bX^\inv\h}^2.
\end{align*}
Then the ratio (cf. \eqref{eq.dom.23}) reads
$$\cR(\x) =\frac{|\tfrac{1}{d}\nabla^3 v(\x)[\h,\h,\h]|}{\left(\tfrac{1}{d}\nabla^2 v(\x)[\h,\h]\right)^{1.5}} = \frac{H_3}{H_2^{1.5}},$$
Since $H_1$ \eqref{eq.uh1} is dominated by $H_2$, we could focus on $\frac{\tau_3}{H_2^{1.5}}$. We could bound this ratio via the following maximization problem,
\begin{eqalign}{eq.ratio.concordance}
    \cL(\x) &:= \max_{\h \in \real_+^n}\inner{\bgamma}{[\bX^\inv\h]^3}, ~\st H_2(\h) = 1.
\end{eqalign}
Note if we take $\bM = (1-r) \bGamma + r\bgamma\bgamma^{\top}$, then
$$H_2(\h) =  \inner{\bX^\inv\h}{\bM\bX^\inv\h}.$$
Note $\bM \succeq 0$, and so the Cholesky decomposition of $\bM=  \bL\bL^{\top}$ yields,
\begin{align*}
    \bL             & = \bGamma^{1/2}\left[\sqrt{1-r}\, \bI + (1-\sqrt{1-r}) \sqrt{\bgamma}\sqrt{\bgamma}^{\top} \right], \\
    \bF := \bL^\inv & = \tfrac{1}{\sqrt{1-r}}\left[\bGamma^{\invhalf} - (1-\sqrt{1-r})\sqrt{\bgamma}\1^{\top}\right].
\end{align*}
Change the coordinates $(\bL^{\top})^\inv\y = \bX^\inv\h$, \eqref{eq.ratio.concordance} can be rewritten as,
\begin{align*}
    \cL(\x) := \max_{\y \in \real_+^n}\inner{\bgamma}{[\bF^{\top}\y]^3}, ~\st \|\y\| = 1.
\end{align*}
Clearly, the optimal value of this tensor eigenvalue problem is in the order of $\cO(\|\bgamma^{-0.5}\|_\infty)$. With the definition of active set, this ratio is relevant to $\{\gamma^{(j)}\}_{j\in B(\x)}$. Hence, we have proved \Cref{lem.dom.23}.
\hfill\qed

\subsection{Proof of \Cref{cor.selfconcordant.theta}}
We let $\gamma_B, \bX^\inv_B\h_B$ be the vectors corresponding to the elements in $B(\x)$.
Using this new notation, it holds that
$$\tau_p \stackrel{\eqref{eq.tau.def}}{=} \inner{\bgamma}{[\bX^\inv\h]^p} = \inner{\gamma_B}{[\bX^\inv_B\h_B]^p}.$$ Since $\|\h\|_{\x} \le 1$, we have $$\tau_1^2 \le \tau_2, |\tau_1|^3 \le \tau_2^{\tfrac{3}{2}} \le |\tau_3| < 1.$$
We have
$$
    \nabla^2 v(\x)[\h,\h] = d (1-r)(\tau_2 - \tau_1^2) + d\tau_1^2 \ge d\tau_1^2\ge d \underset{j\in B(\x)}{\min}\{\gamma^{(j)}\}^2 \|\bX_B^\inv\h_B\|^2.
$$
And similarly, $|\nabla^3 v(\x)[\h,\h,\h]| \le T_u \|\bX_B^\inv\h_B\|^3,$ and so $\tfrac{|\nabla^3 v(\x)[\h,\h,\h]|}{T_u} \le \left(\tfrac{\kappa_u^2\nabla^2 v(\x)[\h,\h]}{d}\right)^{\tfrac{3}{2}}$. Hence we have proved \eqref{eq.lhscb.param} by expanding $T_u$ in view of \eqref{eq.bounded.directdev}.
By Theorem 5.1.7 of \cite{nesterovLecturesConvexOptimization2018}, we see $f(\p)$ is also a self-concordant function with the same coefficient.
This completes the proof.
\hfill\qed

\subsection{Summary of SLC and self-concordance properties of the potential function}
We summarize the properties for the potential function as the summation of player-wise dual functions.
\begin{corollary}[SLC of the potential function]\label{coro.phi inequality}
    Suppose that in the market, for each player $i$ the utility function $u_i \in \cA^{k_i,r_i}, d_i = k_i\times r_i$, then for any $\p \in \real_+^n$ and $\q \in \real^n$, the following inequalities hold:
    \begin{subalign}{eq.phi.inequality}
        \|\bP(\nabla \varphi(\p + \q) - \nabla \varphi(\p) - \nabla^2 \varphi(\p) \q)\|
        & \le \tfrac{T_{\varphi} \varrho^2}{1-\rho},      \\
        \left|\varphi(\p+\q) -  \varphi(\p) - \inner{\nabla \varphi(\p)}{\q} - \tfrac{1}{2}\nabla^2 \varphi(\p)[\q,\q]\right|
        & \le \tfrac{T_\varphi\varrho ^3}{6(1-\varrho )}, \\
        \bP\nabla^2\varphi(\p)\bP &\preceq M_\varphi\bI.
    \end{subalign}
    where $\varrho = \|\q\|_{\p} < 1$. The coefficients are defined as
    \begin{equation}\label{eq.coeff.varphi}
        T_{\varphi} = \sum_{i\in\cI} \tfrac{w_i}{d_i} T_{f_i},
        \quad T_{f_i} = \max\left\{\tfrac{6d_i}{(1-r_i)^2}, 2d_i\right\},
        \quad M_\varphi = \sum_{i\in\cI}\tfrac{w_i(1+|r_i|)}{1-r_i}.
    \end{equation}
\end{corollary}
\subsubsection*{Proof of \Cref{coro.phi inequality}}
By \cref{corr.potential}, $\varphi(\p) = \inner{\1}{\p} + \sum_{i\in\cI} \frac{w_i}{d_i}f_i(\p)$.
Since for each $i\in\cI$, $f_i \in \cA^{k_i,r_i}$, we refer to \cref{thm.taylor.f,thm.taylor.x} and have
\begin{eqalign}{}
    \|\bP(\nabla \varphi(\p + \q) - \nabla \varphi(\p) - \nabla^2 \varphi(\p) \q)\|
    & \le \sum_{i\in\cI}\frac{w_i}{d_i}\|\bP(\nabla^2f_i(\p + \q) - \nabla^2f_i(\p) - \nabla^2 f_i(\p) \q)\| \\
    & \le \frac{\varrho^2}{2(1-\varrho)} \sum_{i\in\cI}\frac{w_i}{d_i}T_{f_i}.
\end{eqalign}
And similarly,
\begin{eqalign}{}
    & \left|\varphi(\p+\q) -  \varphi(\p) - \inner{\nabla \varphi(\p)}{\q} - \tfrac{1}{2}\nabla^2 \varphi(\p)[\q,\q]\right| \\
    & \qquad \le \sum_{i\in\cI}\frac{w_i}{d_i}\left|f_i(\p + \q) - f_i(\p) - \inner{\nabla\varphi_i(\p)}{\q} - \tfrac{1}{2}\nabla^2 f_i(\p)[\q,\q]\right| \\
    & \qquad \stackrel{\eqref{eq.taylor.f},\eqref{eq.f.tf.def}}{\le} \frac{\varrho^3}{6(1-\varrho)} \sum_{i\in\cI}\frac{w_i}{d_i}T_{f_i}.
\end{eqalign}
Besides,
\begin{eqalign}{}
    \bP\nabla^2\varphi(\p)\bP = \sum_{i\in\cI}\frac{w_i}{d_i}\bP\nabla^2f_i(\p)\bP \stackrel{\eqref{eq.f.hess}}{=}
    \sum_{i\in\cI}\frac{w_i}{d_i}\frac{d_i}{1-r_i}\left(\bGamma_i - r_i\bgamma_i\bgamma_i^{\top}\right)
    \preceq \sum_{i\in\cI}\frac{w_i(1+|r_i|)}{1-r_i}\bI.
\end{eqalign}
Take $M_\varphi = \sum_{i\in\cI}\frac{w_i(1+|r_i|)}{1-r_i}$, and the proof is completed.
\hfill\qed

When \cref{asm.bounded.theta} holds, we know for each $i\in\cI$, $f_i$ is self-concordant. Hence, the potential function is also self-concordant.
\begin{corollary}[Self-concordance of the potential function]\label{coro.when.selfconcordant}
    Suppose that \cref{asm.bounded.theta} holds, then the potential function $\varphi(\p)$ is $C_\varphi$-self-concordant, where $C_\varphi = \max_{i\in\cI}\frac{1}{\sqrt w_i}C_{f_i} = \max_{i\in\cI}\frac{1}{\sqrt w_i}C_{v_i}$.
\end{corollary}
The proof is apparent by \cite[Theorem 5.1.1]{nesterovLecturesConvexOptimization2018}.

\section{Proofs for Approximation of the Hessian}\label{appendix.sec.approx.hessian}

\subsection{Proof of \Cref{thm.approx.hessian.dr1}}
We need the following lemma.
\begin{lemma}\label{lem.finite.sample.bound}
    Let $(\cU, \cF, \tau)$ be a probablity space, consider $\{u_i\}_{i\in\cI}$ a sequence of \emph{independent} measurable functions $u_i \in \cL^0(\cU, \cA^{k_i, r_i})$.
    Let $\bH(\p)$ be the affine-scaled Hessian operator defined in \eqref{eq.affine.hessian}, and $\widetilde{\bH}(\p)$ be defined as in \eqref{eq.hessian.approx}. Assume that
    $$
        \frac{w_1r_1}{1-r_1} = \cdots = \frac{w_mr_m}{1-r_m}
        \quad\text{and}\quad \vbE((\bgamma_i-\vbE(\bgamma_i))(\bgamma_i-\vbE(\bgamma_i))^{\top}) = \Sigma, \forall i\in \cI,$$
    then the approximation error
    \begin{eqalign}{eq.approximate Hp}
        \widetilde{\bH}(\p) - \bH(\p) \succeq 0.
    \end{eqalign}
    and
    \begin{equation}
        \vbP\left(\widetilde{\bH}(\p) - \bH(\p) \succeq \Omega\varepsilon\bI + \Sigma\right)
        \le 2n\exp\left(-\frac{m\varepsilon^2}{2\sqrt{2}(\lambda_{\max}(\Sigma) +\varepsilon)}\right).
    \end{equation}
\end{lemma}
\subsubsection*{Proof of \Cref{lem.finite.sample.bound}}
Since $\frac{w_1r_1}{1-r_1} = \cdots = \frac{w_mr_m}{1-r_m}$, then by definition,
$\omega_i = \tfrac{1}{\Omega}\frac{w_ir_i}{1-r_i} = \frac{1}{m}.$
By \eqref{eq.hessian.approx}, for simplicity, take
\begin{equation*}
    \bxi = \frac{1}{m}\sum_{i\in\cI}\bgamma_i, \quad \y_i = \bgamma_i - \bxi,\quad
    \hat\Sigma = \frac{1}{m}\sum_{i\in\cI}\y_i\y_i^{\top},
\end{equation*}
Hence, $\vbE(\y_i) = 0$, $\|\y_i\|_2 \le \sqrt 2$, and we know the variance is bounded:
\begin{eqalign}{}
    0 \preceq \va(\y_i) &= \va(\bgamma_i) + \va(\bxi) - 2\cov(\bgamma_i, \bxi) =\Sigma - \frac{1}{m}\Sigma \preceq \Sigma.
\end{eqalign}
The approximation error is
\begin{eqalign}{}
    \bA &= \widetilde{\bH}(\p) - \bH(\p) = \Omega \Xi - \Omega \tilde \Xi = \frac{\Omega}{m} \sum_{i\in\cI}\y_i\y_i^{\top} \succeq 0.
\end{eqalign}
We have $\vbE(\bA) = \Omega\Sigma$.
By \cite[Corollary 6.20]{wainwrightHighDimensionalStatisticsNonAsymptotic2019}, since for all $i$, $\lambda_{\max}(\va(\y_i)) \le \lambda_{\max}(\Sigma)$, it holds that
\begin{eqalign}{eq.delta.sigma}
    \vbP(\|\hat\Sigma - \Sigma\| \ge \varepsilon)
    \le 2n\exp\left(-\frac{m\varepsilon^2}{2\sqrt{2}(\lambda_{\max}(\Sigma) +\varepsilon)}\right).
\end{eqalign}
\hfill\qed{}

We are now ready to show \Cref{thm.approx.hessian.dr1}.
Taking $\epsilon_H = \frac{\varepsilon}{2\Omega}$, by \eqref{eq.delta.sigma}, if
$$
    m \ge \frac{8\sqrt{2}(\Omega^2\lambda_{\max}(\Sigma) + 0.5\Omega\epsilon_H)}{\epsilon_H^2}\log\left(\frac{2n}{\delta}\right),
$$
then
$\|\widetilde{\bH}(\p) - \bH(\p) - \Omega\Sigma\| \le \frac{1}{2}\epsilon_H$
holds with probability at least $1-\delta$.
Since $n = |\cJ| > \left(\frac{2\Omega}{\epsilon_H}\right)^{\frac{1}{k}}$, we have
$\Omega\Sigma \preceq \frac{1}{2}\epsilon_H\bI.$
\hfill\qed{}
\subsection{Proof of \Cref{thm.cond.precond}}

By definition, we have,
\begin{eqalign}{eq:Hc}
    \bH_\mathrm{c} = \left(\smallsum_{i\in\cI} \w_i \bGamma_i\right)^\invhalf
    \left(\smallsum_{i\in\cI} \tfrac{w_i}{(1-r_i)} \bGamma_i - \smallsum_{i\in\cI}{\tfrac{w_ir_i}{1-r_i}}\bgamma_i\bgamma_i^{\top}\right)
    \left(\smallsum_{i\in\cI} \w_i \bGamma_i\right)^\invhalf
\end{eqalign}
Let us show part $(a)$.
Since we assume all ratios are identical: $r_i \equiv r < 1$, take
$$\y_i = \Kc^\invhalf \bgamma_i, \quad
    \bS := \sum_{i\in\mathcal I} w_i \y_i \y_i^\top;
$$
then, \eqref{eq:Hc} equivalents as $\bH_{\mathrm c} = \frac{1}{1-r}(\bI - r\bS).$ Since $\|\y_i\|_2 \le 1$ and $\sum_i w_i = 1$,
$$
    0 \le \lambda_{\min}(\bS) \le \lambda_{\max}(\bS) \le 1.
$$
Two cases are possible:
\begin{itemize}[leftmargin=*]
    \item If $r \ge 0$, we see
          $$
              \lambda_{\min}(\bH_c) \ge 1, \quad \lambda_{\max}(\bH_c) \le \frac{1}{1-r};
          $$
    \item otherwise,
          $$
              \lambda_{\min}(\bH_c) \ge \frac{1}{1-r}, \quad \lambda_{\max}(\bH_c) \le 1.
          $$
\end{itemize}
Taking the ratios yields the desired bounds.

\noindent For part $(b)$, consider an instance where $m=n=2, w_1=w_2=\tfrac12.$
Suppose at some $\p \in \real_+^2$,
$$\bgamma_1=\bgamma_2=\begin{bmatrix}1/2 \\1/2\end{bmatrix}.$$
Then, $\bGamma_i=\diag(\bgamma_i)=\tfrac12 \bI$
and $\bgamma_i\bgamma_i^{\!\top}=\tfrac14 \1\1^{\top}$, and we have
$$
    \bH(\p)
    =\sum_{i=1}^{2}\frac{w_i}{1-r}\bGamma_i
    -\sum_{i=1}^{2}\frac{w_i r}{1-r}\,\bgamma_i\bgamma_i^{\!\top}
    =\frac{1}{2(1-r)}\Bigl[\bI-\frac{r}{2}\,\1\1^{\top}\Bigr].
$$
Without loss of generality, we let the preconditioner be $\bK(a)=\diag(1, a^{-1}), a \ge 1$.
Due to the symmetry, $\bK(a)$ spans the space of diagonal preconditioners. Then with a slight abuse of notation, the centered Hessian is
$$
    2\left[\bK(a)\right]^\invhalf \bH(\p) \left[\bK(a)\right]^\invhalf
    =
    2\begin{bmatrix}
        1 & 0 \\ 0 & a
    \end{bmatrix}
    -r
    \begin{bmatrix}
        1 & \sqrt a \\ \sqrt a & a
    \end{bmatrix}
    =
    \begin{bmatrix}
        2 - r & -r\sqrt{a} \\ -r\sqrt{a} & (2-r)a
    \end{bmatrix} =: \widehat\bH(a).
$$
The characteristic polynomial of the matrix on the right-hand side is
\begin{align*}
    \det(\widehat\bH(a) - \lambda \bI) & = s(\lambda) - r^2a=0; \quad s(\lambda) = (2-r-\lambda)((2-r)a -\lambda)
\end{align*}
The solutions of this equation can be illustrated in \Cref{fig.cond.precond}. It is easy to see the gap between the two eigenvalues $\lambda_2 - \lambda_1$ is minimized when $a = 1$ (and so is $\tfrac{\lambda_2}{\lambda_1}$).
Hence, using $a=1$, we know the optimal diagonal preconditioner yields eigenvalues $(2(1-r), 2)$.
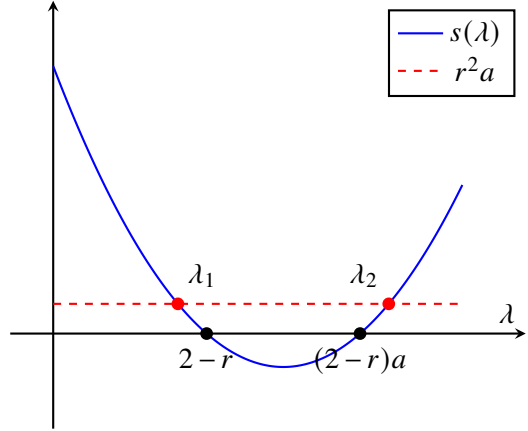
\begin{figure}[h]
    \centering
    \begin{tikzpicture}
        \begin{axis}[
                axis lines=middle,
                xlabel={$\lambda$},
                ylabel={},
                xmin=0, xmax=4.2,
                ymin=-1, ymax=5,
                samples=200,
                domain=0:4,
                thick,
                xtick=\empty,
                ytick=\empty,
                enlargelimits,
                clip=false
            ]

\def\rr{0.5}
            \def\aa{2}

\pgfmathsetmacro{\xone}{2 - \rr}
            \pgfmathsetmacro{\xtwo}{(2 - \rr)*\aa}
            \pgfmathsetmacro{\hline}{\rr*\rr*\aa}

\pgfmathsetmacro{\s}{\xone + \xtwo}
            \pgfmathsetmacro{\d}{(\s)^2 - 4*(\xone*\xtwo - \hline)}
            \pgfmathsetmacro{\lambdaone}{0.5*(\s - sqrt(\d))}
            \pgfmathsetmacro{\lambdatwo}{0.5*(\s + sqrt(\d))}

\addplot[blue, thick] {(\xone - x)*(\xtwo - x)};
            \addlegendentry{$s(\lambda)$}
\addplot[dashed, red, thick, domain=0:4] {\hline};
            \addlegendentry{$r^2 a$}
\addplot[only marks, mark=*, mark size=2pt] coordinates {(\xone, 0) (\xtwo, 0)};
            \node[below] at (axis cs:\xone,0) {\strut $2 - r$};
            \node[below] at (axis cs:\xtwo,0) {\strut $(2 - r)a$};

\addplot[only marks, mark=*, mark size=2pt, red] coordinates {(\lambdaone, \hline) (\lambdatwo, \hline)};
            \node[above right] at (axis cs:\lambdaone,\hline) {\strut $\lambda_1$};
            \node[above left]  at (axis cs:\lambdatwo,\hline) {\strut $\lambda_2$};
        \end{axis}
    \end{tikzpicture}
    \caption{The eigenvalues of the preconditioned Hessian.}\label{fig.cond.precond}
\end{figure}
Hence, if $r < 0$, $\lambda_1 = 2, \lambda_2 = 2(1-r)$; otherwise, $\lambda_1 = 2(1-r), \lambda_2 = 2$.
Both yield the same condition number via $\Kc$ in part $(a)$.
\hfill\qed

\section{Detailed Convergence Analysis of the Interior-Point Strategies}\label{appendix.sec.convergence}

\subsection{Technical lemmas}
The following lemma for the logarithmic barrier function $\cB(\p) = -\sum_{i\in\cJ} \log(p_j)$ is well-known
(see, for example \cite{karmarkarNewPolynomialtimeAlgorithm1984,nesterovInteriorPointPolynomialAlgorithms1994}).
\begin{lemma}
    For any $\p \in \real_+^n$ and $\q \in \real^n$, the following inequality holds:
    \begin{equation}\label{eq.logbar.inequality}
        \cB(\p + \q) - \cB(\p) \le
        \nabla\cB(\p)[\q] + \frac{1}{2}\nabla^2\cB(\p)[\q,\q] + \frac{\varrho^3}{6(1-\varrho)},
    \end{equation}
    where $\varrho = \|\q\|_{\p} < 1$.
\end{lemma}

\subsection{Convergence analysis of the logarithmic barrier method}

In \Cref{alg.logbarrier}, each step is computed by solving the inexact Newton system \eqref{eq.inexact.logbar}.
As indicated by \Cref{lem.logbar.init}, an approximate AC $\p_0$ is obtained by setting $\mu_0$ sufficiently large.
For the rest of the iterations, we introduce the update rule for $\mu$, which ensures that $\p_k \in \cC(\mu_{k+1}, 2Q)$ whenever $\p_k \in \cC(\mu_k, Q)$. Then, we establish the quadratic convergence of the update step in \cref{line.logbar.newton.main}.

\begin{lemma}[Update of $\mu$]\label{lemma.update mu}
    Let $\mu_{k+1} = \sigma \mu_k$, where the shrinkage parameter $\sigma$ is defined by
    \begin{equation}\label{eq.sigma.short}
        \sigma = \frac{Q + \sqrt{n}}{2Q + \sqrt{n}} < 1.
    \end{equation}
    Then, provided that $\p_k \in \cC(\mu_k,Q)$, it follows that $\p_k \in \cC(\mu_{k+1},2Q)$.
\end{lemma}
\subsubsection*{Proof of \Cref{lemma.update mu}}
Notice that
\begin{align*}
    \frac{\|\bP_k \nabla \varphi(\p_k) - \mu_{k+1} \1\|}{\mu_{k+1}}
     & = \frac{\|\bP_k \nabla \varphi(\p_k) - \mu_k \1 + (1-\sigma)\mu_k \1\|}{\sigma\mu_k}                         \\
     & \le \frac{\|\bP_k \nabla \varphi(\p_k) - \mu_k \1\|}{\sigma\mu_k} + \frac{1-\sigma}{\sigma\mu_k}\|\mu_k \1\| \\
     & \le \frac{Q}{\sigma} + \frac{1-\sigma}{\sigma}\sqrt{n},
\end{align*}
where the last inequality holds since $\p_k \in \cC(\mu_k,Q)$ and $\|\1\| = \sqrt{n}$. Substituting the value of $\sigma$ yields $\p_k \in \cC(\mu_{k+1}, 2Q)$.
\hfill\qed

\begin{lemma}[Quadratic convergence near the central path]\label{lem.logbar.quadconv}
    Suppose $\p_k \in \cC(\mu_{k+1}, 2Q)$. Choose $\epsilon_H \le Q^2 / \mu_{k+1}$. Then the update $\p_{k+1} = \p_k + \bP_k \d_k$ satisfies
    \begin{align*}
        \frac{1}{\mu_{k+1}}\|\bP_{k+1} \nabla \varphi(\p_{k+1}) - \mu_{k+1} \1\|
        \le \left(\frac{T_\varphi}{\mu_{k+1}(1- 2Q)} + 3\right)\left(\frac{\|\bP_k \nabla \varphi(\p_k) - \mu_{k+1} \1\|}{\mu_{k+1}}\right)^2
    \end{align*}
    where $\d_k$ is the solution to \eqref{eq.inexact.logbar} with $\p = \p_k$ and $\mu = \mu_{k+1}$.
\end{lemma}
\subsubsection*{Proof of \Cref{lem.logbar.quadconv}}
Let $\bD_k = \operatorname{diag}(\d_k)$, then $\bP_{k+1} = (\bI + \bD_k)\bP_k$. It follows
\begin{align*}
        & \bP_{k+1} \nabla \varphi(\p_{k+1}) - \mu_{k+1} \1                                                                                                                                             \\
    = ~ & \bP_{k+1} \left(\nabla \varphi(\p_{k+1}) - \nabla \varphi(\p_{k}) - \nabla^2 \varphi(\p_{k}) \bP_k \d_k\right) + \bP_{k+1} \left(\nabla \varphi(\p_{k}) - \mu_{k+1} \bP_k^\inv \1\right)      \\
        & + \mu_{k+1}\bP_{k+1}\bP_k^\inv \1 - \mu_{k+1} \1 + \bP_{k+1} \nabla^2 \varphi(\p_{k}) \bP_k \d_k                                                                                              \\
    = ~ & (\bI + \bD_k)\bP_k \left(\nabla \varphi(\p_{k+1}) - \nabla \varphi(\p_{k}) - \nabla^2 \varphi(\p_{k}) \bP_k \d_k\right) + (\bI + \bD_k)\left(\bP_k \nabla \varphi(\p_k) - \mu_{k+1} \1\right) \\
        & + \mu_{k+1}(\bI + \bD_k) \1 - \mu_{k+1} \1 + (\bI + \bD_k) {\bH}(\p_k) \d_k                                                                                                                   \\
    = ~ & (\bI + \bD_k)\bP_k \left(\nabla \varphi(\p_{k+1}) - \nabla \varphi(\p_{k}) - \nabla^2 \varphi(\p_{k}) \bP_k \d_k\right) - (\bI + \bD_k)\left(\widetilde{\bH}(\p_k) + \mu_{k+1} \bI\right)\d_k \\
        & + \mu_{k+1}(\bI + \bD_k) \1 - \mu_{k+1} \1 + (\bI + \bD_k) {\bH}(\p_k) \d_k                                                                                                                   \\
    = ~ & (\bI + \bD_k)\bP_k \left(\nabla \varphi(\p_{k+1}) - \nabla \varphi(\p_{k}) - \nabla^2 \varphi(\p_{k}) \bP_k \d_k\right) + (\bI + \bD_k)\left({\bH}(\p_k) - \widetilde{\bH}(\p_k)\right)\d_k   \\
        & + \mu_{k+1}(\bI + \bD_k) \1 - \mu_{k+1} \1 - (\bI + \bD_k)\mu_{k+1} \d_k                                                                                                                      \\
    = ~ & (\bI + \bD_k)\bP_k \left(\nabla \varphi(\p_{k+1}) - \nabla \varphi(\p_{k}) - \nabla^2 \varphi(\p_{k}) \bP_k \d_k\right) + (\bI + \bD_k)\left({\bH}(\p_k) - \widetilde{\bH}(\p_k)\right)\d_k   \\
        & - \mu_{k+1} \bD_k \d_k.
\end{align*}
On the one hand, by Corollary \ref{coro.phi inequality}, we have
\begin{align*}
     & \|\bP_k (\nabla \varphi(\p_{k+1}) - \nabla \varphi(\p_{k}) - \nabla^2 \varphi(\p_{k}) \bP_k \d_k)\| \le \frac{T_\varphi \|\d_k\|^2}{2(1 - \|\d_k\|)}.
\end{align*}
On the other hand, by Assumption \ref{asm.approx.hessian}, we have
\begin{align*}
     & \|{\bH}(\p_k) - \widetilde{\bH}(\p_k)\| \le \epsilon_H \le \frac{\|\bP_k \nabla \varphi(\p_k) - \mu_{k+1} \1\|^2}{\mu_{k+1}};
\end{align*}
otherwise, we would have $\tfrac{\|\bP_k \nabla \varphi(\p_k) - \mu_{k+1} \1\|}{\mu_{k+1}} \le \sqrt{\epsilon_H \mu_{k+1}} \le Q$, implying that $\p_k$ is already in the neighborhood. Now we analyze the magnitude of the update direction $\d_k$. Notice that the update direction $\d_k$ satisfies
\begin{align*}
    (\widetilde{\bH}(\p_k) + \mu_{k+1} \bI)\d_k = \bP_k \nabla \varphi(\p_k) - \mu_{k+1} \1,
\end{align*}
where the approximate matrix $\tilde\b\nabla^2\varphi \succeq 0$, and $\p_k \in \cC(\mu_{k+1}, 2Q)$. As a result, we obtain
\begin{align*}
    \|\d_k\| \le \frac{\|\bP_k \nabla \varphi(\p_k) - \mu_{k+1} \1\|}{\mu_{k+1}} < 2Q < 1,
\end{align*}
which further implies
\begin{align*}
     & \|\bP_{k+1} \nabla \varphi(\p_{k+1}) - \mu_{k+1} \1\|                                                                                     \\
     & \qquad \le \|\bI + \bD_k\| \cdot \frac{T_\varphi}{2(1-2Q)} \|\d_k\|^2 + \|\bI + \bD_k\| \cdot \epsilon_H + \mu_{k+1}\|\bD_k \d_k\|        \\
     & \qquad \le \frac{T_\varphi}{1- 2Q} \|\d_k\|^2 + \frac{2\|\bP_k \nabla \varphi(\p_k) - \mu_{k+1} \1\|^2}{\mu_{k+1}} + \mu_{k+1}\|\d_k\|^2.
\end{align*}
The last inequality holds because $\|\bD_k \d_k\| \le \|\d_k\|^2$ and $\|\bI + \bD_k\| \le \|\bI\| + \|\bD_k\| \le 2$. Finally, we conclude
\begin{align*}
    \frac{1}{\mu_{k+1}}\|\bP_{k+1} \nabla \varphi(\p_{k+1}) - \mu_{k+1} \1\|
    \le \left(\frac{T_\varphi}{\mu_{k+1}(1- 2Q)} + 3\right)\left(\frac{\|\bP_k \nabla \varphi(\p_k) - \mu_{k+1} \1\|}{\mu_{k+1}}\right)^2
\end{align*}
since $\|\d_k\|^2 \le \|\bP_k \nabla \varphi(\p_k) - \mu_{k+1} \1\|^2 / \mu_{k+1}^2$. The proof is completed.
\hfill\qed

\subsubsection{Proof of \Cref{thm.complexity.logbarrier}}

As indicated by \Cref{lem.logbar.init}, an approximate AC $\p_0$ is obtained by setting $\mu_0$ sufficiently large.
After that, the main phase produces a linearly decreasing sequence $\mu_1, \mu_2, \dots$. Let $\mu_K$ be the first such iterate that $\mu_K \le \frac{\epsilon}{1 + \sqrt{n}} \le \frac{\epsilon}{Q+\sqrt n}$. By standard arguments in short-step IPMs, we know that $K = \cO(\sqrt{n}\log(\tfrac{\mu_0(1 + \sqrt n)}{\epsilon}))$; see, for example, \cite{yeInteriorPointAlgorithms1997}. When it does, $\|\bP_{K} \nabla \varphi(\p_{K}) - \mu_K\1\| \le Q\mu_{K}$. Hence,
$$
    \|\bP_{K} \nabla \varphi(\p_{K})\|_{\infty} \le \|\bP_{K} \nabla \varphi(\p_{K})\|
    \le Q\mu_{K} + \sqrt n \mu_K \le \epsilon.
$$
Besides, $\bP_{K} \nabla \varphi(\p_{K}) \ge -\epsilon \1$.  In other words, $\p_K$ satisfies the following scaled stationary conditions:
\begin{equation}\label{eq.scaled.fosp}
    \p \in \real_+^n, \quad \|\bP\nabla \varphi(\p)\|_{\infty} \le \epsilon, \quad \bP\nabla \varphi(\p) + \epsilon\1 \in \real_+^n.
\end{equation}
By \citep[Proposition 1]{haeserOptimalityConditionComplexity2019}, $\p_K$ is a stationary point in the sense of \eqref{eq.fosp}. We only have to prove that $\p_k \in \cC(\mu_{k}, Q)$ by mathematical induction on $k$.
Suppose the claim holds for some $k \ge 0$. We now investigate the case $k + 1$. By the induction hypothesis and Lemma~\ref{lemma.update mu}, we have
\begin{align*}
    \frac{\|\bP_k \nabla \varphi(\p_k) - \mu_{k+1} \1\|}{\mu_{k+1}} \le 2Q.
\end{align*}
Without loss of generality, we assume that $\mu_{k+1} \ge \frac{\epsilon}{1 + \sqrt{n}}$; otherwise, $\p_k$ is already an $\epsilon$-approximate stationary point. By Lemma~\ref{lem.logbar.quadconv}, the next iterate $\p_{k+1}$ satisfies
\begin{align*}
    \frac{\|\bP_{k+1} \nabla \varphi(\p_{k+1}) - \mu_{k+1} \1\|}{\mu_{k+1}} \le 4Q^2\left(\frac{T_\varphi}{\mu_{k+1}(1-2Q)} + 3\right).
\end{align*}
It remains to prove that the choice of $Q$ guarantees $\p_{k+1} \in \cC(\mu_{k+1}, Q)$, that is,
\begin{align*}
    \left(\frac{T_\varphi}{\mu_{k+1}(1- 2Q)} + 3\right) \cdot (4Q^2) \le Q.
\end{align*}
This is equivalent to
\begin{equation}\label{eq.quadratic inequality of Q}
    24 Q^2 - \left(14 + 4 \cdot \frac{T_{\varphi}}{\mu_{k+1}}\right)Q + 1 \ge 0.
\end{equation}
Let $\ell(Q) = 24Q^2 - \left(14 + 4 \cdot \frac{T_{\varphi}}{\mu_{k+1}}\right)Q + 1$.
It is easy to see that the equation $\ell(Q) = 0$ admits two positive roots, and we denote the smaller one by $Q^*$. Clearly, any $Q \in (0, Q^*)$ satisfies the inequality \eqref{eq.quadratic inequality of Q}. Using $\mu_{k+1} \ge \frac{\epsilon}{1 + \sqrt{n}}$ and quadratic formula yields
\begin{align*}
    48 Q^* = \left(14 + 4 \cdot \frac{T_{\varphi}}{\mu_{k+1}}\right) - \sqrt{\left(14 + 4 \cdot \frac{T_{\varphi}}{\mu_{k+1}}\right)^2 - 96}
    \ge \frac{48\epsilon}{14\epsilon + 4T_{\varphi}(\sqrt{n} + 1)},
\end{align*}
which validates the choice of $Q$. The proof is completed.
\hfill\qed

\subsection{Convergence analysis of the path-following method}
In what follows, we will use a pair of functions
\begin{align*}
    \omega(t) = t - \log(1 + t), \quad\omega_*(t) = -t - \log(1-t).
\end{align*}
We first justify \Cref{lem.relative.error.hess}.
\subsubsection*{Proof of \Cref{lem.relative.error.hess}}
By \Cref{asm.bounded.theta}, recall the bound on the active set,
\begin{align*}
    \min_{j \in B_i(\x_i(\p))} \gamma_i^{(j)} \ge \tfrac{1}{\kappa_{u_i}}, \forall i \in \cI.
\end{align*}
Then, since $\bgamma_i \in \Delta_n$, we have $\bGamma_i \succeq \bgamma_i \bgamma_i^\top$. Hence,
\begin{align*}
    \bH(\p)[\q,\q] & = \bP\nabla^2 \varphi(\p)\bP[\q, \q] \stackrel{\eqref{eq.f.hess}}{=} \sum_{i\in\cI}\frac{d_i}{1-r_i} \left(\bGamma_i - \bgamma_i\bgamma_i^\top + (1-r_i)\bgamma_i\bgamma_i^\top\right)[\q_{\widehat{B}},\q_{\widehat{B}}] \\
                   & \ge \sum_{i\in\cI} d_i \inner{\bgamma_i}{\q_{B_i}}^2 = \sum_{i\in\cI} d_i \inner{\bgamma_i}{\q_{\widehat{B}}}^2                                                                                                           \\
                   & \ge \min_{i\in\cI} \tfrac{d_i}{\kappa_{u_i}} \|\q_{\widehat{B}}\|^2.
\end{align*}
where $\widehat{B} = \bigcup_{i\in\cI} B_i(\x_i(\p))$.
Besides, for DR1 approximation, we have the same sparsity pattern and active set, and thus by
$\|\widetilde{\bH}(\p) - \bH(\p)\| \le \epsilon_H$ (cf. \Cref{thm.approx.hessian.dr1}), we have
\begin{align*}
    \widetilde{\bH}(\p)[\q,\q] - \bH(\p)[\q,\q] = \widetilde{\bH}(\p)[\q_{\widehat{B}},\q_{\widehat{B}}] - \bH(\p)[\q_{\widehat{B}},\q_{\widehat{B}}] & \le \epsilon_H \|\q_{\widehat{B}}\|^2.
\end{align*}
So we have, $\forall \q \in \mathbb{R}^n$,
\begin{equation}
    \widetilde{\bH}(\p)[\q,\q] \le \left(1 + \frac{\epsilon_H}{\min_{i\in\cI} \tfrac{d_i}{\kappa_{u_i}}}\right) \widetilde{\bH}(\p)[\q,\q].
\end{equation}
This completes the proof. \hfill\qed{}

From \Cref{lem.relative.error.hess}, we have
\[
    \tfrac{1}{1+\delta}\bH(\p)^\inv \preceq \widetilde{\bH}(\p)^\inv \preceq \bH(\p)^\inv.
\]
And the same holds for the unscaled operator,
\begin{equation}\label{eq.bounds.invhess}
    \tfrac{1}{1+\delta}\left(\nabla^2\varphi(\p)\right)^\inv \preceq \left(\widetilde{\nabla^2\varphi}(\p)\right)^\inv \preceq \left(\nabla^2\varphi(\p)\right)^\inv.
\end{equation}
Hence, for any $\g$, the dual norm satisfies
$
    \tfrac{1}{\sqrt{1+\delta}}\|\g\|_{\bH(\p)}^{*} \le \|\g\|_{\widetilde{\bH}(\p)}^{*} \le \|\g\|_{\bH(\p)}^{*}.
$
Recall we use $\bH(\p) = \bP\nabla^2\varphi(\p)\bP$ and the Newton decrement satisfies
$
    \lambda(\p) \coloneqq \|\nabla \varphi(\p)\|_{\nabla^2\varphi(\p)}^* = \|\bP\nabla \varphi(\p)\|_{\bH(\p)}^*.
$
We could write the inexact alternative used in the computation
$$
    \widetilde{\lambda}(\p) \coloneqq \|\nabla \varphi(\p)\|_{\widetilde{\nabla^2\varphi}(\p)} = \|\bP\nabla \varphi(\p)\|_{\widetilde{\bH}(\p)}^*.
$$
We are ready to see the following descent property under relative inexactness.
\begin{lemma}\label{lem.descent.pathf}
    For the inexact Newton update \eqref{eq.inexact.pathf}, let $\q = \p_+ - \p = \bP\d$, and suppose that $\tilde\lambda(\p) < \frac{1}{C_{\varphi}}$. We have
    \begin{align*}
        \lambda(\p_+) \le \frac{C_{\varphi} \|\q\|_{\nabla^2\varphi}}{(1 - C_{\varphi} \|\q\|_{\nabla^2\varphi})^2} + \frac{\delta\|\q\|_{\nabla^2\varphi}}{1-C_{\varphi}\|\q\|_{\nabla^2\varphi}}.
    \end{align*}
\end{lemma}
\subsubsection*{Proof of \Cref{lem.descent.pathf}}
By \eqref{eq.inexact.pathf}, $\d = \widetilde{\bH}^\inv (\bP\nabla\varphi_t(\p))$,
Then
\begin{eqalign}{eq.sizeof.q}
    \|\q\|_{\nabla^2\varphi}^2
    & = \nabla^2\varphi [\bP\q, \bP\q] =  \nabla^2\varphi[\bP\widetilde{\bH}^\inv \bP\nabla \varphi(\p), \bP\widetilde{\bH}^\inv\bP\nabla\varphi_t(\p)] \\
    & \le \left(\widetilde{\nabla^2\varphi}\right)^\inv[\nabla \varphi,\nabla\varphi] =  \widetilde{\lambda}^2 \le \frac{1}{C^2_\varphi}.
\end{eqalign}
Hence, the next iterate $\p_+ \in \cP^\circ$ \cite[Theorem 5.1.7]{nesterovLecturesConvexOptimization2018}; it holds that $(1-C_{\varphi} \|\q\|_{\nabla^2\varphi})^2 \nabla^2\varphi(\p) \preceq \nabla^2\varphi(\p_+)$. Then we have
\begin{align*}
    \lambda(\p_+) & = \left(\nabla^2\varphi(\p_+)^\inv[\nabla\varphi(\p_+),\nabla\varphi(\p_+)]\right)^{\frac{1}{2}}                                                  \\
                  & \le \left(\frac{(\nabla^2\varphi)^\inv}{(1-C_{\varphi} \|\q\|_{\nabla^2\varphi})^2}[\nabla\varphi(\p_+),\nabla\varphi(\p_+)]\right)^{\frac{1}{2}} \\
                  & = \frac{1}{1-C_{\varphi} \|\q\|_{\nabla^2\varphi}} \|\nabla\varphi(\p_+)\|_{\nabla^2\varphi}^*.
\end{align*}
Since by \eqref{eq.inexact.pathf.original}
\begin{align*}
    \nabla\varphi(\p_+) = \nabla \varphi(\p) + \int_0^1 \nabla^2\varphi(\p+\tau\q)\q~ \diff\tau = \underbrace{\left(\int_0^1 \nabla^2\varphi(\p+\tau\q) -  \widetilde{\nabla^2\varphi} ~\diff\tau\right)}_{\coloneqq \bG} \q,
\end{align*}
we obtain
\begin{align*}
     & \left(\|\nabla\varphi(\p_+)\|_{\nabla^2\varphi}^*\right)^2 =\left(\bG (\nabla^2\varphi)^\inv \bG\right)[\q, \q]                                       \\
     & = \inner{\q}{\nabla^2\varphi^{\frac{1}{2}} \underbrace{\left(\nabla^2\varphi^{-\frac{1}{2}}\bG \nabla^2\varphi^{-\frac{1}{2}}\right)}_{\coloneqq \bQ}
    \left(\nabla^2\varphi^{-\frac{1}{2}} \bG \nabla^2\varphi^{-\frac{1}{2}}\right) \nabla^2\varphi^{\frac{1}{2}}\q}                                          \\
     & \le  \|\q\|_{\nabla^2\varphi}^2 \cdot \|\bQ\|^2.
\end{align*}
By \cite[Corollary 5.1.5]{nesterovLecturesConvexOptimization2018}, we have
\begin{align*}
    \left(-C_{\varphi}\|\q\|_{\nabla^2\varphi} + \frac{1}{3}C_{\varphi}^2 \|\q\|_{\nabla^2\varphi}^2\right)\nabla^2\varphi \preceq \int_0^1 \nabla^2\varphi(\p+\tau\q) -  \nabla^2\varphi \diff\tau \preceq \frac{C_{\varphi} \|\q\|_{\nabla^2\varphi}}{1 - C_{\varphi} \|\q\|_{\nabla^2\varphi}}\nabla^2\varphi.
\end{align*}
From relative error of the Hessian, we have
$0 \preceq \widetilde{\nabla^2\varphi} - \nabla^2\varphi \preceq \delta \nabla^2\varphi.$
Hence,
\begin{align*}
    \bQ = & \int_0^1 \nabla^2\varphi^{-\frac{1}{2}}\left(\nabla^2\varphi(\p+\tau\q) -  \widetilde{\nabla^2\varphi}\right)\nabla^2\varphi^{-\frac{1}{2}} \diff\tau                                      \\
    =     & \int_0^1 \nabla^2\varphi^{-\frac{1}{2}}\left(\nabla^2\varphi(\p+\tau\q) -  \nabla^2\varphi + \nabla^2\varphi - \widetilde{\nabla^2\varphi}\right)\nabla^2\varphi^{-\frac{1}{2}} \diff\tau.
\end{align*}
Now we can bound $\|\bQ\| \le \frac{C_{\varphi} \|\q\|_{\nabla^2\varphi}}{1 - C_{\varphi} \|\q\|_{\nabla^2\varphi}} + \delta$. Therefore, we conclude
\begin{align*}
    \lambda(\p_+) \le \frac{C_{\varphi} \|\q\|_{\nabla^2\varphi}^2}{(1 - C_{\varphi} \|\q\|_{\nabla^2\varphi})^2} + \frac{\delta\|\q\|_{\nabla^2\varphi}}{1-C_{\varphi}\|\q\|_{\nabla^2\varphi}}.
\end{align*}
The proof is completed. \hfill\qed{}

\subsubsection{The analysis on local convergence}
Two consequences are immediate from the above fact.
Firstly, an analogous region of quadratic convergence also exists for inexact update just as a standard Newton's method. Since $\lambda(\p) \le \sqrt{1+\delta} \cdot \widetilde{\lambda}(\p)$, for some required accuracy $\epsilon > 0$, the following iterate $\p$ shall suffice,
\begin{align*}
    \widetilde{\lambda}(\p) \le \frac{\epsilon}{2} \le \frac{\epsilon}{\sqrt{1+\delta}}.
\end{align*}
Recall we show that $\|\q\|_{\nabla^2\varphi} \le \widetilde{\lambda}(\p)$ (cf. \eqref{eq.sizeof.q}). The above lemma implies that
\begin{align*}
    \widetilde{\lambda}(\p_+) \le \lambda(\p_+) \le \frac{C_{\varphi} \widetilde{\lambda}^2(\p)}{(1 - C_{\varphi} \widetilde{\lambda}(\p))^2} + \frac{\delta\widetilde{\lambda}(\p)}{1-C_{\varphi}\widetilde{\lambda}(\p)}.
\end{align*}
Without loss of generality, we assume $\widetilde{\lambda}(\p) > \tfrac{\epsilon}{2}$; otherwise, the dual norm of the residual is already small. Let $\delta \le \tfrac{C_{\varphi} \epsilon}{2} \le C_{\varphi} \widetilde{\lambda}(\p)$, we have
\begin{align*}
    \widetilde{\lambda}(\p_+) \le \frac{C_{\varphi} \widetilde{\lambda}^2(\p)}{(1 - C_{\varphi} \widetilde{\lambda}(\p))^2} + \frac{C_{\varphi} \widetilde{\lambda}^2(\p)}{1-C_{\varphi}\widetilde{\lambda}(\p)}.
\end{align*}
Therefore, a region of quadratic convergence region can be defined:
\begin{eqalign}{eq.pathf.quadregion}
    \cQ \coloneqq \left\{\p: \widetilde{\lambda}(\p) \le \frac{0.3}{C_{\varphi}}\right\}.
\end{eqalign}
Besides, for any $\p \in \cQ$, the next iterate $\p_+$ satisfies
\begin{align*}
    \widetilde{\lambda}(\p_+) \le \left(\frac{1}{0.49} + \frac{1}{0.7}\right) C_{\varphi} \widetilde{\lambda}^2(\p) \approx \frac{0.312}{C_{\varphi}} \le \frac{1}{2C_{\varphi}} < \frac{1}{C_{\varphi}}.
\end{align*}
This guarantees that \cref{lem.descent.pathf} can be used at iterate $\p_+$, i.e., the inexact Newton's update can be applied at $\p_+$. On the other hand, as for the centering condition in \Cref{alg.pathfollowing},
\begin{equation}
    \mathcal{C}(t) \ = \ \left\{\p:  \|\nabla \varphi(\p) - t \nabla\varphi(\p_0)\|^*_{\widetilde{\nabla^2\varphi}} \le \frac{\beta}{C_{\varphi}}\right\},
\end{equation}
we arrive at the following lemma.
\begin{lemma}\label{lem.rules.betagamma}
    Choose $\beta \in (0,1)$ and $\gamma \in (0,1)$ satisfies $\beta + \gamma < 1$. If
    \begin{equation}\label{eq.rules.betagamma}
        \frac{(1+\delta)(\beta + \gamma)^2}{(1-(\beta + \gamma)\sqrt{1+\delta})^2} + \frac{\delta(\beta + \gamma)\sqrt{1+\delta}}{1-(\beta + \gamma)\sqrt{1+\delta}} \le \beta.
    \end{equation}
    Then $(\p_{k+1}, t_{k+1}) \in \mathcal{C}$ whenever $(\p_k, t_k) \in \mathcal{C}$.
\end{lemma}
\subsubsection*{Proof of \Cref{lem.rules.betagamma}}
Use $\widetilde{\nabla^2\varphi} = \widetilde{\nabla^2\varphi}(\p_k)$, and $\widetilde{\nabla^2\varphi}_{k+1} = \widetilde{\nabla^2\varphi}(\p_{k+1})$ for short. Take,
\[
    \Delta t_k:=\min\left\{\frac{\gamma}{C_{\varphi}
    \|\nabla\varphi(\p_0)\|_{\widetilde{\nabla^2\varphi}}^{*}}, t_k\right\},
\]
then
\begin{align*}
    \nabla\varphi_{t_{k+1}}(\p_k)
     & = \nabla\varphi(\p_k)-t_{k+1}\nabla\varphi(\p_0)
    = \nabla\varphi(\p_k)-t_{k}\nabla\varphi(\p_0) + (t_{k} - t_{k+1})\nabla\varphi(\p_0) \\
     & =\nabla\varphi_{t_{k}}(\p_k)+\Delta t_k\nabla\varphi(\p_0).
\end{align*}
Using $\Delta t_k \le \tfrac{\gamma}{C_{\varphi}\|\nabla\varphi(\p_0)\|_{\widetilde{\nabla^2\varphi}}^{*}}$ yields
\[
    \| \nabla\varphi_{t_{k+1}}(\p_k)\|_{\widetilde{\nabla^2\varphi}}^{*} =  \|\nabla\varphi_{t_{k}}(\p_k)+\Delta t_k\nabla\varphi(\p_0)\|_{\widetilde{\nabla^2\varphi}}^{*}
    \le
    \|\nabla\varphi_{t_{k}}(\p_k)\|_{\widetilde{\nabla^2\varphi}}^{*}+\frac{\gamma}{C_{\varphi}}.
\]
Since $(\p_k, t_k) \in \mathcal{C}$, we have $\|\nabla\varphi_{t_{k}}(\p_k)\|_{\widetilde{\nabla^2\varphi}}^* \le \frac{\beta}{C_{\varphi}}$. Since $\beta + \gamma < 1$, it follows
\begin{eqalign}{eq.lnbd.p.nextt}
    \| \nabla\varphi_{t_{k+1}}(\p_k)\|_{\widetilde{\nabla^2\varphi}}^{*} \le \frac{\beta}{C_{\varphi}} + \frac{\gamma}{C_{\varphi}} < \frac{1}{C_{\varphi}}.
\end{eqalign}
In view of \Cref{lem.descent.pathf}, we have
\begin{align*}
     & \|\nabla\varphi_{t_{k+1}}(\p_{k+1})\|_{\widetilde{\nabla^2\varphi}_{k+1}}^{*} \stackrel{\eqref{eq.bounds.invhess}}{\le} \|\nabla\varphi_{t_{k+1}}(\p_{k+1})\|_{\nabla^2\varphi}^{*}                                                                                     \\
     & \ \le \frac{C_{\varphi} \left(\| \nabla\varphi_{t_{k+1}}(\p_k)\|_{\nabla^2\varphi}^{*}\right)^2}{(1 - C_{\varphi} \| \nabla\varphi_{t_{k+1}}(\p_k)\|_{\nabla^2\varphi}^{*})^2}
    + \frac{\delta\| \nabla\varphi_{t_{k+1}}(\p_k)\|_{\nabla^2\varphi}^{*}}{1-C_{\varphi}\| \nabla\varphi_{t_{k+1}}(\p_k)\|_{\nabla^2\varphi}^{*}}                                                                                                                             \\
     & \stackrel{\eqref{eq.bounds.invhess}}{\le} \frac{C_{\varphi} \left(\sqrt{1+\delta}\| \nabla\varphi_{t_{k+1}}(\p_k)\|_{\widetilde{\nabla^2\varphi}}^{*}\right)^2}{(1 - C_{\varphi}\sqrt{1+\delta}\| \nabla\varphi_{t_{k+1}}(\p_k)\|_{\widetilde{\nabla^2\varphi}}^{*})^2}
    + \frac{\delta\sqrt{1+\delta} \| \nabla\varphi_{t_{k+1}}(\p_k)\|_{\widetilde{\nabla^2\varphi}}^{*}}{1-C_{\varphi}\sqrt{1+\delta}\| \nabla\varphi_{t_{k+1}}(\p_k)\|_{\widetilde{\nabla^2\varphi}}^{*}}                                                                      \\
     & \stackrel{\eqref{eq.rules.betagamma}}{\le}  \frac{(1+\delta)(\beta + \gamma)^2}{C_{\varphi}(1-(\beta + \gamma)\sqrt{1+\delta})^2}
    + \frac{\delta(\beta + \gamma)\sqrt{1+\delta}}{C_{\varphi}(1-(\beta + \gamma)\sqrt{1+\delta})}                                                                                                                                                                             \\
     & \ \ \le \frac{\beta}{C_{\varphi}}.
\end{align*}
It implies that $(\p_{k+1}, t_{k+1}) \in \mathcal{C}$. \hfill\qed{}

Whenever $t_{k+1} = 0$, the update reduces to the standard inexact Newton's update (also recall \eqref{eq.inexact.pathf.original})
\begin{align*}
    \p_{k+1} = \p_k - \widetilde{\nabla^2\varphi}(\p_k)^\inv \nabla\varphi(\p_k).
\end{align*}
Since $\widetilde{\lambda}(\p_{k+1}) \le \tfrac{\beta}{C_{\varphi}}$, as long as $\beta < 0.3$, we could say the next iterate $\p_{k+1}$ is already in the region of quadratic convergence \eqref{eq.pathf.quadregion}.
The method arrives a local region and will terminate in at most $\cO(\log\log(\tfrac{1}{\epsilon}))$ iterations.
Therefore, the non-asymptotic behavior of this scheme is relavant when $t_k \neq 0$.

\begin{figure}[h]
    \centering
    \includegraphics[width=0.5\textwidth]{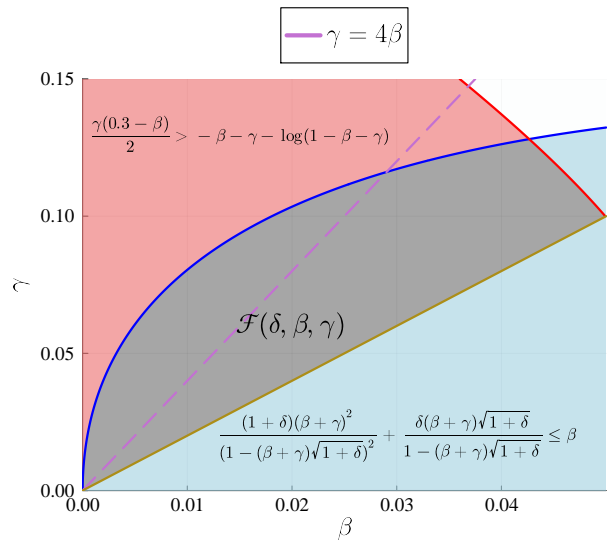}
    \caption{Feasible region for the parameters $(\delta, \beta, \gamma)$ in \eqref{eq.pathf.choice.bgchoice}. The grey region in the middle denotes the valid choice of parameters. The curves denote the boundaries defined by inquality \eqref{eq.pathf.choice.1}, \eqref{eq.pathf.choice.2} and $\gamma \ge 2\beta$ respectively.}
    \label{fig.pathf.choice.bgchoice}
\end{figure}

\subsubsection{Global convergence behavior}
As the previous analysis shows, the path-following method requires a proper choice of parameters. To facilitate the descent property, we choose $(\delta, \beta, \gamma)$ to satisfy the following inequalities:
\begin{subalign}{eq.pathf.choice.bgchoice}
    & \delta \le \tfrac{C_{\varphi} \epsilon}{2},                                                                                       \\
    & \beta + \gamma < 1, \beta < 0.3, \gamma < 1, \gamma > 2\beta,                                                                                       \\
    \label{eq.pathf.choice.1}    & \frac{(1+\delta)(\beta + \gamma)^2}{(1-(\beta + \gamma)\sqrt{1+\delta})^2} + \frac{\delta(\beta + \gamma)\sqrt{1+\delta}}{1-(\beta + \gamma)\sqrt{1+\delta}} \le \beta, \\
    \label{eq.pathf.choice.2}      & \frac{\gamma(0.3 - \beta)}{2} > \omega_*\left(\beta+\gamma\right) = -\beta-\gamma - \log(1-\beta-\gamma).
\end{subalign}
The first inequality can be relaxed to $\delta \le \tilde \lambda(\p)$. In this way, we could set an adaptive accuracy if using the Krylov solvers, or DR1 approximation if the market is sufficiently large.
Note that the above system can be satisfied mildly. For example, if $\delta = 10^{-3}$, the feasible region, $\cF(\delta, \beta, \gamma)$, is shown as the dark grey area in \Cref{fig.pathf.choice.bgchoice}. We could select some $\gamma = 4\beta$, e.g., $\beta = 0.01, \gamma = 0.04$.

In the rest of this section, we assume that $\widetilde{\lambda}(\p_k) > \tfrac{0.3}{C_{\varphi}}$ for all $k = 0, 1, \ldots, K$. Hence,
\begin{align*}
    \|\p_{k+1} - \p_k\|_{\nabla^2\varphi}^2
     & =  \left(\nabla^2\varphi\right)[(\widetilde{\nabla^2\varphi})^\inv\nabla\varphi_{t_{k+1}}(\p_k),(\widetilde{\nabla^2\varphi})^\inv\nabla\varphi_{t_{k+1}}(\p_k)] \\
     & \stackrel{\eqref{eq.bounds.invhess}}{\le}(\widetilde{\nabla^2\varphi})^\inv[\nabla\varphi_{t_{k+1}}(\p_k),\nabla\varphi_{t_{k+1}}(\p_k)]
    = (\|\nabla\varphi_{t_{k+1}}(\p_k)\|_{\widetilde{\nabla^2\varphi}}^*)^2                                                                                             \\
     & \stackrel{\eqref{eq.lnbd.p.nextt}}{\le} \frac{(\beta+\gamma)^2}{C_{\varphi}^2};
\end{align*}
so $\p_{k+1} \in \cP^\circ$ due to that $\beta + \gamma < 1$. Besides, we also have
\begin{eqalign}{}
    (\widetilde{\nabla^2\varphi})^\inv\nabla\varphi_{t_{k+1}}(\p_k)
    &= (\widetilde{\nabla^2\varphi})^\inv(\nabla\varphi(\p_k) - t_{k+1}\nabla\varphi(\p_0))                                                                                                                                                                                                     \\
    &= (\widetilde{\nabla^2\varphi})^\inv\left(\nabla\varphi(\p_k) - \left(t_k - \frac{\gamma}{C_{\varphi}\|\nabla\varphi(\p_0)\|^*_{\widetilde{\nabla^2\varphi}}}\right)\nabla\varphi(\p_0)\right)                                                                                             \\
    &= (\widetilde{\nabla^2\varphi})^\inv\left(\nabla\varphi(\p_k) - t_k \nabla\varphi(\p_0) + \frac{\gamma \nabla\varphi(\p_0)}{C_{\varphi}\|\nabla\varphi(\p_0)\|^*_{\widetilde{\nabla^2\varphi}}}\right).
\end{eqalign}

\begin{lemma}\label{lem.pathf.global.desc}
    Suppose that $\widetilde{\lambda}(\p_k) > \tfrac{0.3}{C_{\varphi}}$ holds, then
    \begin{eqalign}{eq.pathf.global.desc}
        \varphi(\p_k) - \varphi(\p_{k+1}) \ge \frac{\gamma - 2\beta}{2C_{\varphi}}t_k \|\nabla\varphi(\p_0)\|^*_{\widetilde{\nabla^2\varphi}} + \frac{\gamma}{2C_{\varphi}} \left(\widetilde{\lambda}(\p_k) - \frac{\beta}{C_{\varphi}}\right) - \frac{1}{C_{\varphi}^2}\omega_*\left(\beta+\gamma\right).
    \end{eqalign}
    Furthermore, under the parameter choice \eqref{eq.pathf.choice.bgchoice}, we have
    \begin{align*}
        \varphi(\p_k) - \varphi(\p_{k+1}) \ge \frac{\gamma - 2\beta}{2C_{\varphi}}t_k \|\nabla\varphi(\p_0)\|^*_{\widetilde{\nabla^2\varphi}}.
    \end{align*}
\end{lemma}
\subsubsection*{Proof of \Cref{lem.pathf.global.desc}}
Since $\varphi$ is $C_\varphi$-self-concordant, we have
\begin{eqalign}{}
    (\ddagger) = \varphi(\p_k) & - \varphi(\p_{k+1}) + \frac{1}{C_{\varphi}^2}\omega_*\left(C_{\varphi}\|\p_{k+1} - \p_k\|_{\nabla^2\varphi}\right)                                                           \\
    & \ge \inner{\nabla\varphi(\p_k)}{\p_k - \p_{k+1}}                                                                                                                             \\
    & = \inner{\nabla\varphi(\p_k)}{(\widetilde{\nabla^2\varphi})^\inv\nabla\varphi_{t_{k+1}}(\p_k)}                                                                               \\
    & = \inner{\nabla\varphi(\p_k)}{(\widetilde{\nabla^2\varphi})^\inv(\nabla\varphi(\p_k) - t_k \nabla\varphi(\p_0))} \\
    & \qquad + \frac{\gamma}{C_{\varphi}\|\nabla\varphi(\p_0)\|^*_{\widetilde{\nabla^2\varphi}}}\inner{\nabla\varphi(\p_k)}{(\widetilde{\nabla^2\varphi})^\inv\nabla\varphi(\p_0)}.
\end{eqalign}
Reorganizing the terms, we have
\begin{eqalign}{eq.pathf.stepchange}
    (\ddagger)  & = (\|\nabla\varphi(\p_k) - t_k \nabla\varphi(\p_0)\|_{\widetilde{\nabla^2\varphi}}^*)^2
    + \inner{t_k \nabla\varphi(\p_0)}{(\widetilde{\nabla^2\varphi})^\inv(\nabla\varphi(\p_k) - t_k \nabla\varphi(\p_0))}                                                                        \\
    & \qquad + \frac{\gamma}{C_{\varphi}\|\nabla\varphi(\p_0)\|^*_{\widetilde{\nabla^2\varphi}}}\inner{\nabla\varphi(\p_k)}{(\widetilde{\nabla^2\varphi})^\inv\nabla\varphi(\p_0)} \\
    & \ge (\|\nabla\varphi(\p_k) - t_k \nabla\varphi(\p_0)\|_{\widetilde{\nabla^2\varphi}}^*)^2 - t_k\|\nabla\varphi(\p_0)\|_{\widetilde{\nabla^2\varphi}}^* \|\nabla\varphi(\p_k)
    - t_k \nabla\varphi(\p_0)\|_{\widetilde{\nabla^2\varphi}}^*                                                                                                                                  \\
    & \qquad + \frac{\gamma}{C_{\varphi}\|\nabla\varphi(\p_0)\|^*_{\widetilde{\nabla^2\varphi}}}\inner{\nabla\varphi(\p_k)}{(\widetilde{\nabla^2\varphi})^\inv\nabla\varphi(\p_0)} \\
    & \ge - t_k\|\nabla\varphi(\p_0)\|_{\widetilde{\nabla^2\varphi}}^* \|\nabla\varphi(\p_k) - t_k \nabla\varphi(\p_0)\|_{\widetilde{\nabla^2\varphi}}^* \\
    & \qquad + \frac{\gamma}{C_{\varphi}\|\nabla\varphi(\p_0)\|^*_{\widetilde{\nabla^2\varphi}}}\inner{\nabla\varphi(\p_k)}{(\widetilde{\nabla^2\varphi})^\inv\nabla\varphi(\p_0)}.
\end{eqalign}
Now we bound the latter term. Notice that $(\p_k,t_k)\in \mathcal{C}$, then we have
\begin{eqalign}{eq.pathf.att.centering}
    \left(\|\nabla\varphi(\p_k) - t_k \nabla\varphi(\p_0)\|_{\widetilde{\nabla^2\varphi}}^*\right)^2 \le \frac{\beta^2}{C_{\varphi}^2}.
\end{eqalign}
It implies
\begin{align*}
         & \widetilde{\lambda}^2(\p_k) - 2t_k\inner{\nabla\varphi(\p_k)}{(\widetilde{\nabla^2\varphi})^\inv\nabla\varphi(\p_0)}
    + t_k^2 \left(\|\nabla\varphi(\p_0)\|_{\widetilde{\nabla^2\varphi}}^*\right)^2 \le \frac{\beta^2}{C_{\varphi}^2}                              \\
    \qlq & \inner{\nabla\varphi(\p_k)}{(\widetilde{\nabla^2\varphi})^\inv\nabla\varphi(\p_0)} \ge \frac{1}{2t_k}\left(\widetilde{\lambda}^2(\p_k)
    + t_k^2 \left(\|\nabla\varphi(\p_0)\|_{\widetilde{\nabla^2\varphi}}^*\right)^2 - \frac{\beta^2}{C_{\varphi}^2}\right).
\end{align*}
Combining \eqref{eq.pathf.stepchange}, we obtain
\begin{align*}
     & \varphi(\p_k) - \varphi(\p_{k+1}) + \frac{1}{C_{\varphi}^2}\omega_*\left(C_{\varphi}\|\p_{k+1} - \p_k\|_{\nabla^2\varphi}\right)                                                        \\
     & \ge \frac{\gamma}{2C_{\varphi} t_k\|\nabla\varphi(\p_0)\|^*_{\widetilde{\nabla^2\varphi}}}\left(\widetilde{\lambda}^2(\p_k)
    - \frac{\beta^2}{C_{\varphi}^2}\right)
    + \tfrac{\gamma - 2C_{\varphi} \|\nabla\varphi(\p_k) - t_k \nabla\varphi(\p_0)\|_{\widetilde{\nabla^2\varphi}}^*}{2C_{\varphi}}t_k \|\nabla\varphi(\p_0)\|^*_{\widetilde{\nabla^2\varphi}} \\
     & \stackrel{\eqref{eq.pathf.att.centering}}{\ge} \frac{\gamma}{2C_{\varphi} t_k\|\nabla\varphi(\p_0)\|^*_{\widetilde{\nabla^2\varphi}}}\left(\widetilde{\lambda}^2(\p_k)
    - \frac{\beta^2}{C_{\varphi}^2}\right) + \tfrac{\gamma - 2\beta}{2C_{\varphi}}t_k \|\nabla\varphi(\p_0)\|^*_{\widetilde{\nabla^2\varphi}}.
\end{align*}
Recall that
\begin{align*}
    t_k\|\nabla\varphi(\p_0)\|^*_{\widetilde{\nabla^2\varphi}} \le \|\nabla\varphi(\p_k) - t_k \nabla\varphi(\p_0)\|_{\widetilde{\nabla^2\varphi}}^* + \|\nabla\varphi(\p_k)\|_{\widetilde{\nabla^2\varphi}}^* \le \frac{\beta}{C_{\varphi}} + \widetilde{\lambda}(\p_k).
\end{align*}
Then
\begin{align*}
    \frac{\gamma}{2C_{\varphi} t_k\|\nabla\varphi(\p_0)\|^*_{\widetilde{\nabla^2\varphi}}}\left(\widetilde{\lambda}^2(\p_k) - \frac{\beta^2}{C_{\varphi}^2}\right)
    \ge \frac{\gamma}{2C_{\varphi}} \left(\widetilde{\lambda}(\p_k) - \frac{\beta}{C_{\varphi}}\right).
\end{align*}
Also, notice that
\begin{align*}
    \omega_*\left(C_{\varphi}\|\p_{k+1} - \p_k\|_{\nabla^2\varphi}\right) \le \omega_*\left(\beta+\gamma\right).
\end{align*}
Therefore, we conclude
\begin{align*}
    \varphi(\p_k) - \varphi(\p_{k+1}) \ge \frac{\gamma - 2\beta}{2C_{\varphi}}t_k \|\nabla\varphi(\p_0)\|^*_{\widetilde{\nabla^2\varphi}}
    + \frac{\gamma}{2C_{\varphi}} \left(\widetilde{\lambda}(\p_k) - \frac{\beta}{C_{\varphi}}\right) - \frac{1}{C_{\varphi}^2}\omega_*\left(\beta+\gamma\right).
\end{align*}
By combining the above lemma and the parameter choice in \eqref{eq.pathf.choice.bgchoice}, we obtain
\begin{eqalign}{eq.pathf.desc.withparam}
    \varphi(\p_k) - \varphi(\p_{k+1}) & \ge \frac{\gamma - 2\beta}{2C_{\varphi}}t_k \|\nabla\varphi(\p_0)\|^*_{\widetilde{\nabla^2\varphi}} + \frac{\gamma}{2C_{\varphi}} \left(\widetilde{\lambda}(\p_k) - \frac{\beta}{C_{\varphi}}\right) - \frac{1}{C_{\varphi}^2}\omega_*\left(\beta+\gamma\right) \\
    & \ge \frac{\gamma - 2\beta}{2C_{\varphi}}t_k \|\nabla\varphi(\p_0)\|^*_{\widetilde{\nabla^2\varphi}}.
\end{eqalign}
This completes the proof. \hfill\qed

To proceed, we need the following corollary.
\begin{corollary}\label{cor.pathf.desc.small}
    Suppose that $(\delta, \beta, \gamma)$ is selected such that \eqref{eq.pathf.choice.bgchoice} holds. At iteration $k$, if
    \begin{eqalign}{eq.pathf.desc.small}
        \varphi(\p_k) - \varphi(\p_{k+1}) \le \frac{\gamma(\gamma - 2\beta)}{2C_{\varphi}^2}.
    \end{eqalign}
    then $t_{k+1} = 0$, $\p_{k+1} \in \cQ$. Otherwise,
    \begin{eqalign}{eq.pathf.linear.conv.informal}
        t_{k+1} \le t_k \left(1 - \frac{\gamma(\gamma - 2\beta)}{2C_{\varphi}^2(\varphi(\p_k) - \varphi(\p_{k+1}))}\right).
    \end{eqalign}
\end{corollary}
\subsubsection*{Proof of \Cref{cor.pathf.desc.small}}
Suppose \eqref{eq.pathf.desc.small} holds, then by \eqref{eq.pathf.desc.withparam}, we have
\begin{align*}
    \frac{\gamma - 2\beta}{2C_{\varphi}}t_k \|\nabla\varphi(\p_0)\|^*_{\widetilde{\nabla^2\varphi}} \le \frac{\gamma(\gamma - 2\beta)}{2C_{\varphi}^2},
\end{align*}
which further implies $t_k \le \frac{\gamma}{C_{\varphi} \|\nabla\varphi(\p_0)\|^*_{\widetilde{\nabla^2\varphi}}}$.
By the update of $t_{k+1}$, we know $t_{k+1} = 0$, and $\p_{k+1} \in \cQ$.
On the opposite, we would have,
\begin{align*}
    \varphi(\p_k) - \varphi(\p_{k+1}) > \frac{\gamma(\gamma - 2\beta)}{2C_{\varphi}^2},
\end{align*}
and thus
\begin{align*}
    t_{k+1} = t_k - \frac{\gamma}{C_{\varphi}\|\nabla\varphi(\p_0)\|^*_{\widetilde{\nabla^2\varphi}}} = t_k \left(1 - \frac{\gamma}{C_{\varphi} t_k \|\nabla\varphi(\p_0)\|^*_{\widetilde{\nabla^2\varphi}}}\right) \le t_k \left(1 - \frac{\gamma(\gamma - 2\beta)}{2C_{\varphi}^2(\varphi(\p_k) - \varphi(\p_{k+1}))}\right).
\end{align*}
\hfill \qed

\subsubsection{Proof of \Cref{thm.convergence.rate.pathf}}
We could now evaluate how fast $t_k$ decreases. Note that by \eqref{eq.pathf.linear.conv.informal},
\begin{eqalign}{}
    t_{k+1} \le \prod_{j=0}^{k}\left(1 - \frac{\gamma(\gamma - 2\beta)}{2C_{\varphi}^2(\varphi(\p_j) - \varphi(\p_{j+1}))}\right)
\end{eqalign}
By inequality of arithmetic and geometric means one has
\begin{eqalign}{}
    \left(\prod_{j=0}^{k}\left(1 - \frac{\gamma(\gamma - 2\beta)}{2C_{\varphi}^2(\varphi(\p_j) - \varphi(\p_{j+1}))}\right)\right)^{\tfrac{1}{k+1}}
    \le 1 - \frac{1}{k+1}\sum_{j=0}^k \frac{\gamma(\gamma - 2\beta)}{2C_{\varphi}^2(\varphi(\p_j) - \varphi(\p_{j+1}))}.
\end{eqalign}
Using the harmonic mean,
\begin{eqalign}{}
    \frac{1}{k+1}\sum_{j=0}^k \frac{\gamma(\gamma - 2\beta)}{2C_{\varphi}^2(\varphi(\p_j) - \varphi(\p_{j+1}))}
    \ge \frac{k+1}{\sum_{j=0}^k\frac{2C_{\varphi}^2(\varphi(\p_j) - \varphi(\p_{j+1}))}{\gamma(\gamma - 2\beta)}}
    =\frac{\gamma(\gamma - 2\beta)(k+1)}{2C_{\varphi}^2(\varphi(\p_0) - \varphi(\p_{k+1}))}.
\end{eqalign}
Hence,
\begin{align*}
    t_{k+1}
     & \le \left(1 - \frac{1}{k+1}\sum_{j=0}^k \frac{\gamma(\gamma - 2\beta)}{2C_{\varphi}^2(\varphi(\p_j) - \varphi(\p_{j+1}))}\right)^{k+1} \\
     & \le \left(1 - \frac{\gamma(\gamma - 2\beta)(k+1)}{2C_{\varphi}^2(\varphi(\p_0) - \varphi(\p_{k+1}))}\right)^{k+1}.
\end{align*}
This completes the proof. \hfill\qed

\section{Extra Proofs for Further Applications}\label{appendix.sec.extra}

\subsection{Proof of \Cref{lem.linear.lump.foc}}

The first-order condition is,
\begin{subarr}{eq.linear.ipm.foc}{rl}
    \frac{\c}{\inner{\c}{\x}} + \sigma \bX^\inv\1 - \lambda \p = 0 \\
    0 \le \lambda ~\perp ~ w - \inner{\p}{\x} \ge 0.
\end{subarr}
Hence, $1 + \sigma n - \lambda w = 0$, $\lambda = \frac{1 + \sigma n}{w} > 0$
and the budget is tight.
Reusing the first equation, we have,
\begin{align*}
    \x = \sigma\left(\frac{1+\sigma n}{w}\p - \frac{1}{u}\c\right)^\inv, \quad u = \inner{\c}{\x},
\end{align*}
where ``$\inv$" denotes the element-wise inverse. Hence, $u$ is the root to the following equation,
\begin{align*}
    \psi(u) = \sum_{j\in\cJ}\sigma c^{(j)}\left(\frac{1+\sigma n}{w}p^{(j)} - \frac{1}{u}c^{(j)}\right)^\inv - u = 0.
\end{align*}
We use a slightly shifted bidding vector,
\begin{align*}
    \bgamma = (1+\sigma n)\frac{\bX\p}{w} - \sigma \1= \frac{\bX\c}{\inner{\c}{\x}}
\end{align*}
Taking the derivative over $\p$ at \eqref{eq.linear.ipm.foc},
$$
    - \left(\frac{\c\c^\top}{\inner{\c}{\x}^2} + \sigma \bX^{-2}\right)\nabla\x = \lambda \bI.
$$
Solving for $\nabla\x$ completes the proof. \hfill\qed
\subsection{Proof of \cref{thm.linear.lump.slc}}
By similar arguments as before (cf. \cref{thm.taylor.f}), take $\h = \nabla^2f(\p)\q$, note that $\frac{\bX\p}{w} = \frac{1}{1+\sigma n}(\bgamma + \sigma \1)$,
\begin{eqalign}{}
    \bX^\inv\h &= \bX^\inv\nabla^2f(\p)\q
    =\frac{1}{w} \left(\frac{\lambda}{\sigma} \bX\bP
    - \frac{\lambda}{\sigma(\sigma + \|\bgamma\|^2)}\bgamma\bgamma^\top\bX\bP\right)\bP^\inv\q \\
    &= \lambda\left(\frac{1}{\sigma}\frac{1}{1+\sigma n}(\bGamma + \sigma \bI) - \frac{1}{\sigma(\sigma + \|\bgamma\|^2)}\bgamma\bgamma^\top\frac{1}{1+\sigma n}(\bGamma + \sigma \bI)\right)\bP^\inv\q.
\end{eqalign}

We take $\bS = \left(\frac{1}{\sigma}\frac{1}{1+\sigma n}(\bGamma + \sigma \bI) - \frac{1}{\sigma(\sigma + \|\bgamma\|^2)}\bgamma\bgamma^\top\frac{1}{1+\sigma n}(\bGamma + \sigma \bI)\right)$, $\f = \bP^\inv\q$.
Next we observe that
\begin{eqalign}{}
    \|\bS\| &= \left\|\frac{1}{\sigma}\frac{1}{1+\sigma n}(\bGamma + \sigma \bI) - \frac{1}{\sigma(\sigma + \|\bgamma\|^2)}\bgamma\bgamma^\top\frac{1}{1+\sigma n}(\bGamma + \sigma \bI)\right\| \\
    &\le \frac{1}{\sigma(1+\sigma n)}\left\|\bI - \frac{1}{\sigma + \|\bgamma\|^2}\bgamma\bgamma^\top\right\|\|(\bGamma + \sigma \bI)\| \\
    &\le \frac{(1+\sigma)}{\sigma(1+\sigma n)} .
\end{eqalign}
Notice that $\lambda = \frac{1+\sigma n}{w}$, $d=1$ because of linear utility, then,
\begin{eqalign}{}
    \nabla^3 f(\p)[\q,\q,\q] &= w^3\nabla^3v(\x)[\lambda\bX\bS\f,\lambda\bX\bS\f,\lambda\bX\bS\f] \\
    &= (1+\sigma n)^3\left(-2\sigma\bX^{-3} - 2 \frac{1}{\inner{\c}{\x}^3}\c\otimes\c\otimes\c\right)[\bX\bS\f,\bX\bS\f,\bX\bS\f] \\
    &= -2(1+\sigma n)^3\left(\sigma\inner{[\bS \f]^3}{\1} + \inner{\bgamma}{\bS \f}^3\right).
\end{eqalign}
Hence,
\begin{eqalign}{}
    \|\nabla^3 f(\p)[\q,\q,\q]\| \le 2(1+\sigma n)^3 \|\bS\|^3 \|\q\|_{\p}^3 \le \frac{2(1+\sigma)^3}{\sigma^3} \|\q\|_{\p}^3.
\end{eqalign}
In view of \cref{eq.f.hesstimesq,eq.ub3.f3qqq}, the same results can be worked out easily but tediously. \hfill\qed

\subsection{Proof of \Cref{thm.lump.affine.hess}}

By writing out the Lagrangian, we get
\begin{equation*}
    L(\x, \lambda, \y) = -v(\x) - \lambda \inner{\p}{\x} + \lambda w + \inner{\y}{\bA\x}.
\end{equation*}
Then by perturbation theory, the derivative of $f(\p)$ is $\nabla f(\p) = -\lambda\x(\p)$,
and the primal-dual condition is given by
\begin{eqalign}{}
    \nabla v(\x) + \lambda \p + \bA^{\top}\y - \s = 0, ~\bA\x = 0, ~\bX\s = 0,\\
    \real_+ \ni \lambda \perp w - \inner{\p}{\x} \in \real_+.
\end{eqalign}
Multiplying the first equation by $\x$, similar to full-dimensional case \eqref{eq.logump}, we get
\begin{equation}\label{eq.affine.lambda}
    \lambda = \tfrac{d}{w}, ~\inner{\p}{\x} = w.
\end{equation}
To compute the Hessian $\nabla^2 f(\p) = - \lambda \nabla \x = - \lambda \tfrac{\diff}{\diff \p} \x(\p)$,
let us differentiate the first-order conditions,
\begin{equation}
    \begin{bmatrix}
        \nabla^2 v(\x) & \bA^{\top} & - \bI \\
        \bA            & 0          & 0     \\
        \bS            & 0          & \bX
    \end{bmatrix}
    \begin{bmatrix}
        \nabla \x \\
        \nabla \y \\
        \nabla \s
    \end{bmatrix}
    = -\begin{bmatrix}
        \lambda \bI \\
        0           \\
        0
    \end{bmatrix}.
\end{equation}
Because $\bX\s = 0$, solving the above linear system yields,
\begin{align*}
    \nabla \x = - \lambda (\bW^\inv-\bW^\inv\bA^{\top}(\bA\bW^\inv\bA^{\top})^\inv\bA\bW^\inv), ~\bW = \nabla^2 v(\x).
\end{align*}
This completes the proof.
\hfill\qed

\subsection{Proof of \Cref{thm.scaled.lipschitz.network}}
To show this result, let us define the following bordered matrix $\bQ(\p)$ and $\bW$-orthogonal projector $\Pi_{\bW, \bA}$,
\begin{equation}\label{eq.f.hess.border}
    \bQ(\p) = \begin{bmatrix}
        \nabla^2 v(\x(\p)) & \bA^{\top} \\
        \bA                & 0
    \end{bmatrix}, ~\Pi_{\bW, \bA} = \bI-\bW^\inv\bA^{\top}(\bA\bW^\inv\bA^{\top})^\inv\bA.
\end{equation}
Then, we may write
\begin{equation}\label{eq.f.hess.border.inv}
    \bQ^\inv = \begin{bmatrix}
        \bW^\inv- \bW^\inv\bA^\top\left(\bA\bW^\inv\bA^\top\right)^\inv\bA\bW^\inv & \bW^\inv\bA^\top\left(\bA\bW^\inv\bA^\top\right)^\inv \\
        \left(\bA\bW^\inv\bA^\top\right)^\inv\bA\bW^\inv                           & \left(\bA\bW^\inv\bA^\top\right)^\inv
    \end{bmatrix}
\end{equation}
and so $\nabla \x = -\lambda \bQ^\inv_{[1:n,1:n]} = -\lambda \Pi_{\bW, \bA}\bW^\inv$.
Notice that
\begin{eqalign}{eq.f.hesstimesq.af}
    \nabla^3 f(\p)[\q,\q,\q] & = \lim_{t\to 0} \frac{1}{t} \left(\nabla^2f(\p+t\q) - \nabla^2f(\p)\right)[\q,\q]                                                                                       \\
    & \stackrel{\eqref{eq.f.hess.border.inv}}{=} \lambda^2\lim_{t\to 0} \frac{1}{t} \left(\bQ(\p+t\q)^\inv - \bQ(\p)^\inv\right)[[\q; 0], [\q; 0]]                       \\
    & = \lambda^2\lim_{t\to 0} \frac{1}{t} \left(\bQ-\bQ(\p+t\q)\right)[\bQ(\p+t\q)^\inv[\q; 0], \bQ^\inv[\q; 0]]               \\
    & = \lambda^2\lim_{t\to 0} \frac{1}{t} \left(\bQ-\bQ(\p+t\q)\right)[\lim_{t\to 0}\bQ(\p+t\q)^\inv[\q; 0], \bQ^\inv[\q; 0]].
\end{eqalign}
By definition \eqref{eq.f.hess.border},
\begin{align*}
    \lim_{t\to 0} \frac{1}{t} \left(\bQ(\p)-\bQ(\p+t\q)\right)
     & =
    \lim_{t\to 0} \frac{1}{t}
    \begin{bmatrix}
        \nabla^2 v(\x(\p))-\nabla^2 v(\x(\p+t\q)) & 0 \\
        0                                         & 0
    \end{bmatrix}                    \\
     & = \begin{bmatrix}
             - \nabla^3 v(\x(\p))[\nabla \x(\p)\q] & 0 \\
             0                                     & 0
         \end{bmatrix}.
\end{align*}
Hence, by $\bQ^\inv[\q; 0] = [\tfrac{1}{\lambda^2} \nabla^2f(\p)\q; 0]$,
\begin{eqalign}{}
    \nabla^3 f(\p)[\q,\q,\q] & = \lambda^2\begin{bmatrix}
        - \tfrac{1}{\lambda} \nabla^3 v(\x(\p))[\nabla^2f(\p)\q] & 0 \\
        0                                                        & 0
    \end{bmatrix}[\bQ^\inv[\q; 0],\bQ^\inv[\q; 0]] \\
    & = - \tfrac{1}{\lambda^3} \nabla^3 v(\x(\p))[\nabla^2 f(\p)\q,\nabla^2 f(\p)\q,\nabla^2 f(\p)\q].
\end{eqalign}
Besides, take $\h = \nabla^2 f(\p)\q$,
\begin{align*}
    (\bX(\p))^\inv \h & = \bX(\p)^\inv\nabla^2 f(\p)\q = \lambda^2\bX(\p)^\inv\Pi_{\bW, \bA}\bW^\inv\q                                             \\
                      & = \bX(\p)^\inv\Pi_{\bW, \bA}\left(\frac{d}{(1-r)} \bP^\inv \left(\bGamma - r \bgamma\bgamma^{\top}\right)\bP^\inv\right)\q \\
                      & = \bX(\p)^\inv\Pi_{\bW, \bA}\bX(\p)\left(\frac{d}{w(1-r)}\left(\bI - r \1\bgamma^{\top}\right)\bP^\inv\right)\q            \\
                      & = \Pi_{\bW, \bA}\left(\frac{d}{w(1-r)}\left(\bI - r \1\bgamma^{\top}\right)\bP^\inv\right)\q.
\end{align*}
Hence,
$$
    \|\bX(\p)^\inv \h\| = \left\|\Pi_{\bW, \bA}\left(\tfrac{d}{w(1-r)}\left(\bI - r \1\bgamma^{\top}\right)\bP^\inv\right)\q\right\| \le \left\|\Pi_{\bW, \bA}\right\|\left\|\tfrac{d}{w(1-r)}(\bI - r \1\bgamma^{\top})\bP^\inv\q\right\|.
$$
Because $\Pi_{\bW, \bA}$ is a projection matrix, $\|\Pi_{\bW, \bA}\| \le 1$.
In view of \cref{eq.f.hesstimesq,eq.ub3.f3qqq}, the same results can be worked out.
\hfill\qed

 \end{APPENDICES}

\end{document}